\documentclass[12pt, a4paper]{amsart}
\usepackage[utf8]{inputenc}
\usepackage{amsfonts,amssymb,amsxtra,amsthm,amsmath,amscd,mathrsfs,cite}
\usepackage[all]{xy}

\usepackage{tikz}
\usepackage{verbatim}
\usepackage[normalem]{ulem}
\usepackage{soul}
\usepackage{color}

\setstcolor{red}

\makeatletter
\@namedef{subjclassname@2010}{%
	\textup{2010} Mathematics Subject Classification}
\makeatother

\def \N {{\mathbb N}}
\def \P {{\mathbb P}}

\def\Li{\hbox{{\rm Li}}}

\def\leq{\leqslant}
\def\geq{\geqslant}
\def\le{\leqslant}
\def\ge{\geqslant}

\theoremstyle{plain}
\newtheorem{theorem}{Theorem}[section]
\newtheorem{proposition}{Proposition}[section]
\newtheorem{lemma}[proposition]{Lemma}

\theoremstyle{remark}

\numberwithin{equation}{section}
\addtocounter{footnote}{1}

\begin{document}
	
\title[Theorem $(1+1.9)$ on the Goldbach Conjecture]  
{Theorem $\mathbf{(1+1.9)}$ on the Goldbach Conjecture} 
	
\author{Jiamin Li \& Jianya Liu}
	
\address{%
Jiamin Li
		\\
Mathematical Research Center
		\\
Shandong University
		\\
Jinan, Shandong 250100
		\\
P. R. China
	}
\email{lijiamin@mail.sdu.edu.cn}
\address{%
Jianya Liu
		\\
State Key Laboratory of Cryptography and Digital Economy Security \& Mathematical Research Center
		\\
Shandong University
		\\
Jinan, Shandong 250100
		\\
P. R. China
}
\email{jyliu@sdu.edu.cn}
	
\date{\today}
	
\begin{abstract} 
For \( 1 \leq a \leq 2 \), we say Proposition \((1+a)\) holds if every sufficiently large even 
integer \(N\) can be written as  
\[
N = p + rq, \quad r \leq q^{\,a-1},
\]
where \(r\) is either \(1\) or prime, and \(p,q\) are primes.  
Thus Proposition $(1+1)$ is essentially the binary Goldbach Conjecture, and 
Proposition \((1+2)\) is Chen's theorem. 
We prove unconditionally that Proposition \((1+1.9)\) is true. 
Assuming the Elliott--Halberstam Conjecture, the exponent \(1.9\) can be improved to \(1.4\).  
Analogously, Proposition \((1-a)\) is formulated for the Twin Prime Conjecture. 
Unconditionally, we prove Proposition \((1-1.75)\), and under the Elliott--Halberstam Conjecture, Proposition \((1-1.4)\).   

For six decades, a substantial theoretical divide has persisted between Propositions $(1+2)$ and $(1+1)$, 
and likewise between Propositions $(1-2)$ and $(1-1)$. By constructing new weighted sieves and adopting new 
analytic tools, this paper establishes a connecting pathway between them and achieves 
breakthroughs in this line of research. 
\end{abstract}
	
	
\maketitle
	
\section{Introduction}\label{sec1}
\subsection{The binary Goldbach Conjecture}
The binary Goldbach Conjecture states that every even integer greater than $2$ is the sum of two primes. 
For an even integer $N>2$, let $D_{1,1}(N)$ denote the number of solutions to
\begin{equation}\label{BinGol}
N = p + q
\end{equation}
with primes $p$ and $q$.  
Hardy and Littlewood \cite{HL23} proposed the more precise conjecture that, as $N \to \infty$,
\begin{equation}\label{D11}
D_{1,1}(N) \sim 2C(N) \frac{N}{\log^2 N},
\end{equation}
where
\begin{equation}\label{CN}
C(N) = \prod_{p>2} \left(1 - \frac{1}{(p-1)^2}\right) \prod_{\substack{p>2 \\ p \mid N}} \frac{p-1}{p-2}.
\end{equation}
The factor $2C(N)$ is precisely the Hardy--Littlewood singular series associated with the binary Goldbach Conjecture, 
written here in the form customary in sieve theory.

While \eqref{BinGol} remains out of reach, one may relax one variable, say $q$, 
to an almost prime $P_k$, i.e., an integer with at most $k$ prime factors, and study the equation 
\begin{equation}\label{Pro/1+k}
N = p + P_k.
\end{equation}
Thus $P_1$ denotes a prime, and $4$ and $6$ are examples of $P_2$.  
We say Proposition $(1+k)$ holds if \eqref{Pro/1+k} is solvable for all sufficiently large even $N$.  
The strongest result in this direction remains Chen’s theorem of $(1+2)$.  
Chen \cite{Chen66, Chen73} proved that for large even $N$, if $D_{1,2}(N)$ counts solutions to $N = p + P_2$, then
\begin{equation}\label{D12}
D_{1,2}(N) \ge 0.67 \, C(N) \frac{N}{\log^2 N},
\end{equation}
with $C(N)$ as in \eqref{CN}. The constant $0.67$ has since been improved successively to
\[
0.689,\; 0.7544,\; 0.81,\; 0.8285,\; 0.836,\; 0.867,\; 0.899
\]
by Halberstam--Richert \cite{HR74}, Chen \cite{Chen78II}, Chen \cite{Chen78}, Cai--Lu \cite{CaLu02}, Wu \cite{Wu04}, Cai \cite{Cai08a}, and Wu \cite{Wu08}, respectively.

\medskip 

Although Proposition $(1+2)$ was long regarded as the closest result 
to the Goldbach Conjecture, the weighted sieve method has 
advanced considerably in the six decades since Chen’s work \cite{Chen66}. Combined with improved 
levels of distribution for primes, which means certain averages behave as if various Riemann Hypotheses 
hold, these developments open new avenues for attacking problems concerning primes.  
To frame such progress, we introduce a parametrized proposition that may play a role 
analogous to the Goldston--Pintz--Yildirim method for twin primes. 

For six decades, a substantial theoretical divide has persisted between Propositions $(1+2)$ and $(1+1)$. 
By constructing new weighted sieves and adopting new 
analytic tools, this paper establishes a connecting pathway between them and achieves 
breakthroughs in this line of research. 

\medskip
\noindent
{\bf Proposition $(1+a)$.}   
{\it Fix $1 \le a \le 2$. Every sufficiently large even integer $N$ can be written as
\begin{equation}\label{p+rq/r<}
N = p + r q, \quad r \le q^{\,a-1},
\end{equation}
where $r$ is either $1$ or prime, and $p, q$ are primes.} 
\medskip

The definition is consistent: if $a = 1$, then $r = 1$ and we essentially recover the binary Goldbach Conjecture,  
Proposition $(1+1)$. If $a = 2$, we obtain Chen’s theorem $(1+2)$. For $1 \le a < 2$, the statement is stronger than Chen’s theorem.

Our first result is the following.

\begin{theorem}\label{Goldbach1}
Proposition $(1+1.9)$ is true.  
More precisely, for every sufficiently large even $N$, if $D_{1,a}(N)$ denotes the number of solutions to \eqref{p+rq/r<}, then
\[
D_{1,1.9}(N) > 0.0004 \, C(N) \frac{N}{\log^2 N},
\]
where $C(N)$ is as in \eqref{CN}.
\end{theorem}
We remark that the exponent $1.9$ can be improved to $1.894$ by using more intricate parameters and the double sieve method of Wu \cite{Wu08}.

\subsection{The Twin Prime Conjecture}

Let \( x \) be a sufficiently large real number, and let \( \pi_{1,1}(x) \) denote the number of solutions of the equation  
\[
p + 2 = q,
\]  
where \( p, q \) are primes with \( p, q \le x \).  
The Twin Prime Conjecture asserts that this equation admits infinitely many solutions as \( x \to \infty \).  
Hardy and Littlewood~\cite{HL23} conjectured the asymptotic formula  
\begin{equation}\label{pi11}
\pi_{1,1}(x) \sim C \frac{x}{\log^2 x}
\end{equation}  
as \(x \to \infty\), with  
\begin{equation}\label{C2}
C = 2 \prod_{p>2} \left( 1 - \frac{1}{(p-1)^2} \right).
\end{equation}
Let \( \pi_{1,2}(x) \) denote the number of solutions of the equation  
\[
p + 2 = P_2,
\]  
where \( p \) is a prime and \( P_2 \) is an almost-prime with at most two prime factors.  
In the same papers, Chen~\cite{Chen66, Chen73} proved  
\begin{equation}\label{pi12}
\pi_{1,2}(x) \ge 0.335 \, C \frac{x}{\log^2 x} 
\end{equation}  
for sufficiently large $x$. This constant has been improved successively to  
\[
0.3445,\; 0.3772,\; 0.405,\; 0.71,\; 1.015,\; 1.05,\; 1.0974,\; 1.104,\; 1.123,\; 1.13
\]  
by Halberstam \cite{Hal75}, Chen \cite{Chen78II}, Chen\cite{Chen78}, 
Fouvry--Grupp \cite{FG86}, Liu \cite{Liu89}, Wu \cite{Wu90}, Cai \cite{Cai02}, Wu \cite{Wu04}, 
Cai \cite{Cai08a},  Cai \cite{Cai08b} respectively.

We consider a refinement of this problem. 

\medskip

\noindent
{\bf Proposition $(1-a)$.}   
{\it Fix $1 \le a \le 2$.  
There are infinitely many primes \( p \) such that  
\begin{equation}\label{def/tpc/1a}
p + 2 = r q, \quad r \le q^{\,a-1},
\end{equation}
where \( r \) is either prime or \( 1 \), and \( q \) is prime. }  

\medskip

Progress on the Twin Prime Conjecture has outpaced that on the Goldbach Conjecture since 1980s.  
In this century, developments related to the Goldston--Pintz--Yildirim method~\cite{GPY1,GPY2}, Zhang~\cite{Zhang14}, 
Maynard~\cite{Maynard15}, and others reflect the same trend.  
Within our framework, advances concerning the Twin Prime Conjecture are also more substantial.

\begin{theorem}\label{Twin1}
Proposition \((1-1.75)\) is true.  
More precisely, for every sufficiently large $x$, if $\pi_{1,a}(x)$ denotes the number 
of solutions to \eqref{def/tpc/1a} with $p\leq x$, 
then
\[
\pi_{1,1.75}(x) > 0.036 \, C \frac{x}{\log^2 x},
\]  
where \( C \) is as in \eqref{C2}.
\end{theorem}

We note that $1.75$ can be further improved.  
Although Pascadi's mean-value theorem \cite{Pas24} has a higher distribution level, it is not 
applicable to the sieve weight function constructed in this paper and thus not 
adopted herein. For detailed elaboration, see \S \ref{sec3}. 

\subsection{Conditional results}
Chen's theorem $(1+2)$ relies on the Bombieri--Vinogradov theorem
that if 
\[
R(y; q, a)
=\sum_{\substack{p\leq y \\ p \equiv a \bmod q}} 1 
-\frac{\pi(y)}{\varphi(q)},  
\]
then for any $A>0$, there exists $B=B(A)>0$ such that
\begin{equation}\label{1/BV}
\sum_{q\leq x^{1/2}/\log^{B}x} \max_{y\leq x} \max_{(a,q)=1} 
| R(y; q, a) |
\ll \frac{x}{\log^A x}. 
\end{equation}
Pan and Ding \cite{PanDing79} (see also Pan \cite{Pan/HH}) generalized this result to a weighted variant. 
By adopting this weighted form in place of the original estimate \eqref{1/BV}, they substantially simplified the proof of Chen's theorem \cite{Chen66, Chen73}. Now set
\begin{equation}\label{WBV/Err}
R(y; q, a, m)
=\sum_{\substack{ m p\leq y \\ m p \equiv a \bmod q}} 1 
 -\dfrac{\Li(y/m)}{\varphi(q)}. 
\end{equation}
The corresponding weighted mean-value theorem states that for any $A>0$, there exists $B=B(A)>0$ satisfying
\begin{equation}\label{W/BV}
\sum_{q\leq x^{1/2}/\log^{B}x} \max_{y\leq x} \max_{(a,q)=1} 
\bigg| \sum_{m\leq M} f(m) R(y; q, a, m) \bigg|
\ll \frac{x}{\log^A x}, 
\end{equation}
where $f(m)$ is a fairly general weight function and $M$ is a large parameter. A key merit of this weighted approach is that it removes the complicated arguments used originally by Chen in his proof of the $(1+2)$ theorem, which greatly streamlines the overall proof.

As is well known, the Elliott--Halberstam conjecture claims that the distribution level in the Bombieri--Vinogradov theorem can be lifted to $1-\varepsilon$ for any $\varepsilon>0$. We now formulate a weighted version of the Elliott--Halberstam conjecture with distribution level $\theta$, denoted by $\mathrm{WEH}(\theta)$, where $\theta \in (1/2, 1)$ characterizes the equidistribution quality of primes in arithmetic progressions. Roughly speaking, our weighted Elliott--Halberstam conjecture asserts that the above conclusion still holds when equipped with the weight functions introduced by Pan and Ding. We state the conjecture following the formulation in \cite[\S 8]{PanPan92}.

\medskip
\noindent
{\bf Conjecture WEH$(\theta)$.}\label{EH-New}
{\it
Let $\theta\in (1/2, 1)$ be a constant. Suppose $f(m) \ll 1$ and $\alpha \in (0,1]$. Let $r_1(y)$ be a positive function depending on $x$ such that
$$
r_1(y) \ll x^{\alpha}, \quad y \leqslant x.
$$
Let $r_2(m)$ be a positive function depending on $x$ and $y$ such that
$$
mr_2(m) \ll x, \quad m \leqslant x^{\alpha}, \quad y \leqslant x.
$$
Let $R(y;q,a,m)$ be as in \eqref{WBV/Err}. 
Then for every $A > 0$, there exists a constant $B = B(A) > 0$ such that
\[
\sum_{q \leqslant x^{\theta}/\log^{B} x} 
\max_{y \leqslant x} \max_{(a,q)=1} 
\bigg| \sum_{\substack{m \leqslant x^{1-\alpha} \\ (m,q)=1}} f(m) R(y;q,a,m) \bigg| \ll \frac{x}{\log^A x},
\]
\[
\sum_{q \leqslant x^{\theta}/\log^{B}x} 
\max_{y \leqslant x} \max_{(a,q)=1} \bigg| \sum_{\substack{m \leqslant x^{1-\alpha} \\ (m,q)=1}} f(m) 
R(mr_1(y);q,a,m) \bigg| \ll \frac{x}{\log^A x},
\]
and
\[
\sum_{q \leqslant x^{\theta}/\log^{B}x} 
\max_{y \leqslant x} \max_{(a,q)=1} \bigg| \sum_{\substack{m \leqslant x^{1-\alpha} \\ (m,q)=1}} 
f(m) R(mr_2(m);q,a,m) \bigg| \ll \frac{x}{\log^A x}.
\]
}

We note that ${\rm WEH}(\theta)$ is known to hold for $\theta \le \frac{1}{2}$ (see Lemma \ref{BV-New}). 
We point out that the classical Bombieri--Vinogradov theorem cannot handle the error terms arising from the switching principle. In contrast, the Pan--Ding mean value theorem, is not only sufficient to cover those errors, but also encompasses the classical Bombieri–Vinogradov theorem itself.	
We can now state our conditional result. 

\begin{theorem}\label{Goldbach2}
Assume \(\mathrm{WEH}(0.999)\). Then Proposition \((1+1.4)\) holds. 
More precisely, we have
\[
D_{1,1.4}(N) > 0.062\,C(N) \frac{N}{\log^2 N} 
\]
for every sufficiently large even $N$, where \(C(N)\) is as in \eqref{CN}.
\end{theorem}

Since Conjecture $\mathrm{WEH}(\theta)$ applies to the twin-prime problem without modification, 
it yields an analogous result of equivalent strength. The proof, which follows a parallel argument, is omitted in this paper 
for brevity. 

\begin{theorem}\label{Twin2}
Assume \(\mathrm{WEH}(0.999)\). Then Proposition \((1-1.4)\) holds. More precisely, we have 
\[
\pi_{1,1.4}(x) > 0.031\,C \frac{x}{\log^2 x}
\]
for every sufficiently large $x$, where \(C\) is as in \eqref{C2}.
\end{theorem}

We remark that the exponent \(1.4\) can be further improved. 

\subsection{A related result on the Twin Prime Conjecture}
Our results above can be compared with a theorem of \cite{BH96}, which states that there are infinitely many 
primes \(p\) with \(P^{+}(p+a) > p^{0.676}\). Thus, there are infinitely many primes \(p\) such that 
\[
p+2 = rq, \quad r \le p^{0.324} \le q^{0.48},
\]
where \(q\) is prime and \(r\) is an integer, not necessarily a prime.
	
\subsection{Outline of the paper}

Below we outline the structure of the paper. 
In \S \ref{sec2}, we introduce the linear sieve method in two different forms and list lemmas that facilitate calculations.
In \S \ref{sec3}, we study the equidistribution of prime numbers. In Chen’s era, only the Bombieri -- Vinogradov theorem, 
a substitute for the Generalized Riemann Hypothesis (GRH) for Dirichlet \(L\)-functions, 
was available. Developments in recent decades now allow us to 
surpass the Bombieri--Vinogradov theorem in the context of twin prime problems. In particular, we are 
able to avoid using the GRH for zeta-functions of finite fields (in the sense of Weil) and the GRH 
for zeta-functions related to eigenvalues of the Laplace operator (in the sense of Selberg).

In \S \ref{sec4}, we investigate the weighted sieve method. We also observe that for \(a < 1.5\), the combinatorial structure of 
Propositions \((1+a)\) and \((1-a)\) becomes more intricate. This indicates that once the level of distribution increases, the nature of the problem changes.

In \S\S \ref{sec5}-\ref{sec7}, we prove all the theorems stated in the introduction.

In \S \ref{sec8}, we offer some further remarks.

\medskip

\noindent
\textbf{Notation.}  
Throughout, \(0<\varepsilon<10^{-10}\). The letter \(x\) denotes a large real number, and \(N\) a large even integer.  
Primes are denoted by \(p\), with or without subscripts.  
For \(n \ge 2\), \(P^{-}(n)\) and \(P^{+}(n)\) are the least and largest prime factor of \(n\), respectively.   
The offset logarithmic integral is 
\(\Li(x)=\int_{2}^{x} dt/\log t\).  
Notation concerning various sieves is standard and will be introduced in due course.

\section{Preliminary lemmas}\label{sec2}
This section collects several sieve-theoretic estimates that will be used later.

\subsection{Sieve functions and Buchstab’s function}
The first lemma gives an asymptotic formula for the number of integers free of small prime factors, 
together with some explicit bounds for Buchstab’s function $w(u)$.  The asymptotic formula is a classical result 
(see, e.g., \cite[Lemma 12.1]{FI10}); the explicit expressions for $w(u)$ in small ranges come from \cite[Lemma 20]{Jia96}, 
and the upper bound $w(u)\le 0.561522$ for $u\ge 3.5$ is taken from \cite[Lemma 2.10]{Wu04}. 

\begin{lemma}\label{w}
Let
\[
\Psi(x,z)=|\{n\le x:(n,Q(z))=1\}|, \quad Q(z)=\prod_{p<z} p .
\]
Then for $1<u_0\le u\le 100$,
\[
\Psi(x,x^{1/u})\sim u\,w(u)\,\frac{x}{\log x},
\]
where $w(u)$ is Buchstab’s function determined by the differential-difference equation
\[
\begin{cases}
w(u)=u^{-1}, & 1\le u\le 2,\\[2mm]
(u w(u))' = w(u-1), & u>2 .
\end{cases}
\]
Explicitly,
\[
w(u)=\begin{cases}
\displaystyle\frac1u, & 1\le u\le 2,\\[4mm]
\displaystyle\frac{1+\log(u-1)}u, & 2\le u\le 3,\\[4mm]
\displaystyle\frac{1+\log(u-1)}u
      +\frac1u\int_{2}^{u-1}\frac{\log(t-1)}t\,dt, & 3\le u\le 4 .
\end{cases}
\]
Moreover, for all $u\ge 3.5$ one has $w(u)\le 0.561522$.
\end{lemma}

The next lemma records the standard upper- and lower-bound sieve functions $F(s),f(s)$ of the linear sieve.  The formulas below are taken from \cite[(3.11) and (3.12)]{HHR81}; they are valid for the Iwaniec--Rosser sieve with sifting parameter $\kappa=1$.

\begin{lemma}\label{Ff}
Let $f$ and $F$ be the functions in the linear sieve. 
For $F$, we have  
\renewcommand{\arraystretch}{2.4}
\medskip 
\[
\begin{array}{|c|c|c|c|c|}
\hline
s          & 1\le s\le 3 & 3\le s\le 5 & 5\le s\le 7  \\
\hline
F(s) & F_1(s) & F_2(s) & F_3(s) \\
\hline
\end{array}
\] 
where 
\[
F_1(s) = \frac{2e^{\gamma}}s, 
\]
\[
F_2(s)= \frac{2e^{\gamma}}s\bigg(1+\int_{2}^{s-1}\frac{\log(t-1)}t\,dt\bigg), 
\] 
and 
\[ 
F_3(s) = \frac{2e^{\gamma}}s\bigg(1+\int_{2}^{s-1}\frac{\log(t-1)}t\,dt
      +\int_{2}^{s-3}\frac{\log(t-1)}t\,dt\int_{t+2}^{s-1}
        \frac1u\log\frac{u-1}{t+1}\,du\bigg). 
\]
For $f$, we have  
\renewcommand{\arraystretch}{2.4}
\medskip 
\[
\begin{array}{|c|c|c|c|c|}
\hline
s          & 2\le s\le 4 & 4\le s\le 6 \\
\hline
f(s) & f_1(s) & f_2(s) \\
\hline
\end{array}
\]
where 
\[
f_1(s) = \frac{2e^{\gamma}}s\log(s-1), 
\]
\[
f_2(s)= \frac{2e^{\gamma}}s\bigg(\log(s-1)+\int_{3}^{s-1}\frac{dt}t\int_{2}^{t-1} \frac{\log(u-1)}u\,du\bigg), 
\]
and $\gamma=0.577\ldots$ is Euler’s constant.
\end{lemma}

Numerically one finds
\[
F(7)\le 1.000005,\quad f(6)\ge 0.999895,
\]
so both $F(s)$ and $f(s)$ are extremely close to $1$ for these arguments.  
This explains why, in the above lemma, we only show $s \leq 7$ for $F$ and to $s \leq 6$ for $f$. 
Computations involving $F$ and $f$ for larger $s$ would require additional effort while offering only marginal 
improvement to the final results. 

\subsection{The linear sieve}
In this subsection, we introduce the general setup of the linear sieve; later we will study two versions of the linear sieve theorem. In the proofs of Theorems \ref{Goldbach1} and \ref{Twin1} we shall employ the linear sieve with well-factorable coefficients; for Theorems \ref{Goldbach2} and \ref{Twin2} the standard linear sieve will suffice.

Now let $\mathscr{A}$ be a finite sequence of integers and let $\mathscr{P}$ be a set of primes. For $z \ge 2$ we consider the sifting function
\begin{equation}\label{def/SAPz}
S(\mathscr{A},\mathscr{P},z) = | \{ a \in \mathscr{A} : (a, P(z)) = 1 \} |,
\end{equation}
where
\begin{equation}\label{def/Pz}
P(z) = \prod_{\substack{p < z \\ p \in \mathscr{P}}} p.
\end{equation}
If $d$ is a square‑free integer whose prime factors all belong to $\mathscr{P}$, we write $\mathscr{A}_d$ for the set of elements of $\mathscr{A}$ that are divisible by $d$. We assume an approximation of the form
\begin{equation}\label{RIcon1}
|\mathscr{A}_d| = \frac{\omega(d)}{d} X + r(\mathscr{A},d),
\end{equation}
where $X>1$ is independent of $d$, and $\omega$ is a multiplicative function satisfying
\begin{equation}\label{RIcon2}
0 \le \omega(p) < p \quad \text{for } p \in \mathscr{P}.
\end{equation}
We also define
\begin{equation}\label{def/Vz}
V(z) = \prod_{\substack{p < z \\ p \in \mathscr{P}}} \left( 1 - \frac{\omega(p)}{p} \right),
\end{equation}
and we assume that there exists an absolute constant $K > 1$ such that
\begin{equation}\label{RIcon3}
\frac{V(z_1)}{V(z_2)} \le \frac{\log z_2}{\log z_1} \left( 1 + \frac{K}{\log z_1} \right)
\end{equation}
for $z_2 \ge z_1 \ge 2$. This inequality places us in the setting of the linear sieve.

\subsection{Standard linear sieve}
First, we introduce the standard linear sieve\cite{Iwa81}.

\begin{lemma}\label{Sieve-standard}
We adopt the basic notations and assumptions of the linear sieve \eqref{def/SAPz}--\eqref{def/Vz}.
Assume further that
\begin{equation}\label{Iwacon}
\sum_{\substack{z_1 \le p < z_2}} \frac{\omega(p)}{p} = \log \frac{\log z_2}{\log z_1} 
+ O\left(\frac{1}{\log z_1}\right),\quad z_2 > z_1 \ge 2.
\end{equation}
Then, under these assumptions, the function 
$S(\mathscr{A}, \mathscr{P}, z)$ admits the upper and lower bounds 
\[
S(\mathscr{A}, \mathscr{P}, z) \geq XV(z)\left\{f(s) + O\left(\frac{1}{\log^{1/3} D}\right)\right\} - R_D,
\]
\[
S(\mathscr{A}, \mathscr{P}, z) \leq XV(z)\left\{F(s) + O\left(\frac{1}{\log^{1/3} D}\right)\right\} + R_D, 
\]
where 
\[
s = \frac{\log D}{\log z}, \quad R_D = \sum_{\substack{d < D \\ d \mid P(z)}} |r(\mathscr{A},d)|,
\]
and finally, the functions $f(s)$ and $F(s)$ are exactly those introduced in Lemma \ref{Ff}.
\end{lemma}

\subsection{The linear sieve with well‑factorable coefficients}
In this subsection we present the well‑factorable form of the Rosser--Iwaniec sieve, which will replace the standard linear sieve in certain parts of the argument. We begin with the necessary definitions.

Let $k$ be a positive integer. For any positive integer $n$, denote by $\tau_k(n)$ the number of ways to write $n$ as a product of $k$ positive integers. An arithmetic function $\lambda$ is said to be {\it well‑factorable of level $Q$ and order $k$} 
if it satisfies the following conditions:
\begin{equation}\label{WF1}
	\lambda(q) = 0 \quad \text{for all } q \ge Q,
	\end{equation}
	\begin{equation}\label{WF2}
	|\lambda(q)| \le \tau_k(q) \quad \text{for all } q \ge 1, 
	\end{equation}
and if, for every decomposition $Q = Q_1 Q_2$ with $Q_1, Q_2 \ge 1$, there exist two functions $\lambda_1$ and $\lambda_2$ satisfying \eqref{WF1} and \eqref{WF2} with $Q$ replaced by $Q_1$ and $Q_2$ respectively, such that
	$$
	\lambda = \lambda_1 * \lambda_2,
	$$
where $*$ denotes Dirichlet convolution.
The following lemma\cite[Lemma 2.1]{Wu04} is useful.
\begin{lemma}\label{WuWF}
	If $\lambda'$ is an arithmetical function of level $Q' \le Q$ and of order $k'$, then $\lambda*\lambda'$ is well factorable of level $QQ'$ and of order $k+k'$.
\end{lemma}

Now we state 
the well‑factorable version of the linear sieve\cite{Iwa80} as follows.
	
\begin{lemma}\label{Sieve-WF}
Let  
\begin{equation}\label{Lem24/con}
0 < \eta < 1/8,\quad L=\exp(8\eta^{-3}),\quad 2 \leq z \leq Q^{1/2}. 
\end{equation}
Then under the assumptions \eqref{RIcon1}, \eqref{RIcon2} and \eqref{RIcon3}, we have the inequalities
\begin{equation}\label{RI-up}
S(\mathscr{A}, \mathscr{P}, z) \leq XV(z) \left( F\left( \frac{\log Q}{\log z} \right) + E \right) + \sum_{1 \leq l \leq L }\sum_{q \mid P(z) } \lambda_l^+(q) r(\mathscr{A}, q), 
		\end{equation}
		and
		\begin{equation}\label{RI-low}
		S(\mathscr{A}, \mathscr{P}, z) \geq XV(z) \left( f\left( \frac{\log Q}{\log z} \right) + E \right) + \sum_{1 \leq l \leq L }\sum_{q \mid P(z) } \lambda_l^-(q) r(\mathscr{A}, q). 
		\end{equation}
In these formulas,  $\lambda_q^+$ and $\lambda_q^-$ are well-factorable coefficients of order $1$ 
and of level $Q$, and $E$ satisfies
\begin{equation}\label{RI-Error}
E \ll \eta + \eta^{-8} e^K (\log Q)^{-1/3}. 
\end{equation}
\end{lemma}

\section{Equidistribution of primes}\label{sec3} 

This section collects known results on the distribution of primes in arithmetic progressions.  
The first one is a weighted version of the Bombieri--Vinogradov theorem due to Pan and Ding \cite{PanDing79}, 
which implies Bombieri-Vinogradov's theorem. In \S\ref{sec1} we have stated it in the form of \cite[\S 8]{PanPan92}.

\begin{lemma}\label{BV-New}
The weighted Elliott--Halberstam Conjecture {\rm WEH}$(\theta)$ holds for every $\theta \le 1/2$.
\end{lemma}

The next result, due to Bombieri--Friedlander--Iwaniec \cite[Theorem 10]{BFI86}, provides a larger level of distribution for well‑factorable coefficients.

\begin{lemma}\label{BFI}
Let $\lambda$ be a well‑factorable function of order $k$ and level $Q = x^{4/7 - \varepsilon}$.  
For any $\varepsilon > 0$ and any $A>0$, uniformly for $x \ge 3$ and $|a| \le \log^A x$,
		$$
		\sum_{\substack{(q,a)=1}} \lambda(q) \left( \pi(x; q, a) - \frac{\Li(x)}{\varphi(q)} \right) \ll_{\varepsilon, k, A} \frac{x}{ \log^A x},
		$$
		where
		$$ 
		\pi(x; q, a)  =
		\sum_{\substack{p \le x\\n \equiv a \bmod q}} 1.
		$$
\end{lemma}

The following result is due to Pascadi \cite[Theorem 1.2(ii)]{Pas24}. We note, however, that the stated improvement applies only to the upper bound sieve, not the lower bound sieve. 
In our proof, lower bound estimates for \( S(\mathscr{A}, \mathscr{P}, z) \) are required. The terms 
that would necessitate upper bound estimates are those involving \( S(\mathscr{A}_p, \mathscr{P}, z) \), 
which cannot be obtained via Pascadi’s method. Consequently, this lemma is not applicable in our setting.

\begin{lemma}\label{Pascadi}
Let  $a \in \mathbb{Z} \setminus \{0\},  A, \varepsilon > 0,   x \geq 2$, and suppose that $Q \leq x^{3/5 - \varepsilon}$, 
and   $(\lambda_q)$   are the well-factorable upper bound linear sieve weights of level  $Q$. Then one has
$$
\sum_{\substack{q \leq Q \\ (q,a)=1}} \lambda_q \left( \pi(x; q, a) - \frac{\Li(x)}{\varphi(q)} \right) \ll_{\varepsilon, A, a} \frac{x}{\log^A x}. 
$$
\end{lemma}
	
In the proofs of Theorems \ref{Goldbach1} and \ref{Twin1}, 
we shall need some mean value theorems with well-factorable or 
almost well-factorable coefficients. The first one is \cite[Lemma 2.5]{Wu04}, which is a corollary 
that combines \cite[Corollary 2]{F87}, \cite[Lemma 6]{FG86} 
and \cite[Proposition]{Iwa93}.
		
\begin{lemma}\label{Wu-level}
Let  $M \geq 1$, $N \geq 1$ and $X := MN$. Let $\{ \alpha_m \}$ and $\{ \beta_n \}$ be 
two sequences of order $k$ supported in $[M, 2M]$ and $[N, 2N]$ respectively. We also 
suppose that for all constant $B>0$, the asymptotic formula 
\begin{equation}
\sum_{\substack{n \equiv n_0 \bmod{k} \\ (n, d) = 1}} 
\beta_n = \frac{1}{\varphi(k)} \sum_{(n, dk) = 1} \beta_n + O_{B,k} \! \left( \frac{N \tau_k(d)}{(\log 2N)^B} \right)
\end{equation}
holds for $d \geq 1 , k \geq 1$ and $(k, n_0) = 1$. If $n$ has a prime factor $p$ 
with $p < \exp\{ (\log \log n)^2 \}$, then $\beta_n = 0$. 
Then for any $A>0$ and for any $\varepsilon>0$ we have
$$
\sum_{\substack{(q,a)=1}}\lambda(q)
\bigg(\sum_{\substack{mn\equiv a\bmod{q}}}\alpha_m\beta_n-\frac{1}{\varphi(q)}\sum_{\substack{(mn,q)=1}}\alpha_m\beta_n\bigg)
\ll_{\varepsilon,A}\frac{X}{\log^A X}
		$$
		uniformly for $|a|\le \log^{A} X$. Here $\lambda(q)$ is a well-factorable function of order $1$ and of level 
$Q:=X^{\theta(\nu)-\varepsilon}$,  where $\theta(\nu)$ is given by 
\renewcommand{\arraystretch}{2.4}
\medskip 
\[
\begin{array}{|c|c|c|c|c|}
\hline
\nu & \varepsilon \le \nu \le \frac{1}{15} & \frac{1}{15} \le \nu \le \frac{1}{10} & \frac{1}{10} \le \nu \le \frac{3}{14} & \frac{3}{14} \le \nu \le \frac{1}{4} \\
\hline
\theta(\nu) & \dfrac{6-5\nu}{10} & \dfrac{1}{2}+\nu & \dfrac{5-2\nu}{8} & \dfrac{3+2\nu}{6} \\
\hline
\nu & \frac{1}{4} \le \nu \le \frac{2}{7} & \frac{2}{7} \le \nu \le \frac{2}{5} & \frac{2}{5} \le \nu \le \frac{1}{2} & \frac{1}{2} \le \nu \le 1-\varepsilon \\
\hline
\theta(\nu) & \dfrac{2-\nu}{3} & \dfrac{2+\nu}{4} & 1-\nu & \dfrac{1}{2} \\
\hline
\end{array}
\] 
\end{lemma}

\medskip 
		
The following lemma, sourced from \cite[Corollary 2 (i)]{F87}, was originally established for 
investigations into the Twin Prime Conjecture. Wu \cite{Wu04} subsequently adapted this lemma to the 
Goldbach Conjecture in 2004. This is the first result, which is valid uniformly for 
$|a|\le X$ and has the level of distribution $>\frac{1}{2}$.
		
\begin{lemma}\label{Fouvry-new}
Let  $M \geq 1$, $N \geq 1$ and $X := MN$. Let $\{ \alpha_m \}$ and $\{ \beta_n \}$ be two sequences of order $k$ supported in $[M, 2M]$ and $[N, 2N]$ respectively. We also suppose that for all $B$, the equality
\begin{equation}
\sum_{\substack{n \equiv n_0 \bmod{k} \\ (n, d) = 1}} 
\beta_n = \frac{1}{\varphi(k)} \sum_{(n, dk) = 1} \beta_n + O_{B,k} \! \left( \frac{N \tau_k(d)}{(\log 2N)^B} \right)
\end{equation}
holds for $d \geq 1 , k \geq 1$ and $(k, n_0) = 1$. 
If $n$ has a prime factor $p$ with $p < \exp\bigl\{ (\log \log n)^2 \bigr\}$, then $\beta_n = 0$. 
Then for any $A>0$ and for any $\varepsilon>0$ we have
		$$
		\sum_{\substack{(q,a)=1}}\lambda(q)\biggl(\sum_{\substack{mn\equiv a\bmod{q}}}\alpha_{m}\beta_{n}-\frac{1}{\varphi(q)}\sum_{\substack{(mn,q)=1}}\alpha_{m}\beta_{n}\biggr)\ll_{\varepsilon,A}\frac{X}{\log^A X}
		$$
		uniformly for $|a|\le X$ 
Here $\lambda(q)$ is a well-factorable function of order $1$ and of level $Q:=X^{5(1-\nu)/9-\varepsilon}$, where
$$
\varepsilon\le\nu:=\frac{\log N}{\log X} 
\le\frac{1}{10}.
$$ 
\end{lemma} 

The lemma above is essential for establishing Proposition $(1+1.9)$. Historically, the constant in 
Proposition $(1+2)$ was successively refined from $0.67$ to $0.81$ by Chen.  
Since a lower bound of at least $0.8429$ is required for Proposition $(1+1.9)$, 
this constant was far from being met in Chen’s time, so he could not have 
proved Proposition $(1+1.9)$ then. 
More importantly, even in the period after Chen, without employing the aforementioned lemma, 
the constant could at best reach $0.8285$ as in Cai and Lu \cite{CaLu02}.  
This still falls short of the $0.8429$ threshold. Thus, without this lemma, 
subsequent advances 
would likewise be insufficient to support the proof of Proposition $(1+1.9)$. 

\medskip 

When we use the weighted inequality, some coefficients are merely almost well-factorable. So we need the following result due to  
Fouvry--Grupp \cite[Theorem 2]{FG89}.
	
\begin{lemma}\label{Fouvry-level}
Let $\lambda$ be a well-factorable function of level $Q_1$ and of order $k$, $\xi$ an arithmetical function 
satisfying the conditions 
$
|\xi(q_2)| \leq \log x
$
and
$\xi(q_2) = 0
$
for
$
q_2 > Q_2.
$
Let $\Lambda$ be the von Mangoldt function. Then we have for any integer $a$, any $\varepsilon > 0$ and any $A > 0$,  
$$
\sum_{\substack{(q_1q_2, a)=1}} \lambda(q_1)\xi(q_2) \left( \pi(x; q_1 q_2, a) - \frac{\Li(x)}{\varphi(q_1 q_2)} \right) \ll_{a, \varepsilon, k, A} 
\frac{x}{\log^A x},
$$  
as long as one of the following four conditions is true:  
\begin{itemize}
\item[(C.1)] $Q_2 \leq Q_1, Q_1 Q_2 \leq x^{4/7 - \varepsilon}$,  
\item[(C.2)] $Q_2 \geq Q_1, Q_1 Q_2^6 \leq x^{2 - \varepsilon}$,  
\item[(C.3)] $\xi(q) = \Lambda(q), Q_1 Q_2 \leq x^{11/20 - \varepsilon}, Q_2 \leq x^{1/3 - \varepsilon}$.
\end{itemize}  
\end{lemma}
	
\section{Weighted sieves}\label{sec4}

In this section, we develop weighted sieve methods to establish the theorems 
stated in the introduction.  
We adopt the standard notation of the linear sieve. Throughout the rest of the paper, the parameter \(z\) satisfies  
\begin{equation}\label{Cond/z}
\log z\asymp \log x, 
\end{equation}
and in both the Goldbach and the twin prime problems we set  
\begin{equation}\label{Def/P}
\mathscr{P}=\{p: p>2\}. 
\end{equation}
We will also require information on the auxiliary sets  
\begin{equation}\label{Def/Pq}
\mathscr{P}(q)=\{p \in \mathscr{P}: (p, q)=1\}. 
\end{equation} 

\subsection{The Goldbach Conjecture (I)}
In this subsection, we develop the weighted sieve method needed for Theorem \ref{Goldbach1}
which deals with the Goldbach problem unconditionally.  
Recall that \( D_{1,a}(N) \) denotes the number of solutions to \eqref{p+rq/r<}, i.e.,  
\begin{equation}\label{D1a}
D_{1,a}(N) = |\{ p \le N : N - p = r q,\; r \le q^{a-1} \}|,
\end{equation}
where \( r \) is either a prime or \( 1 \), and \( p, q \) are primes.  
We now define the set  
\begin{equation}\label{Def/A=Goldbach}
\mathscr{A} = \{ N - p : p < (1-\varepsilon)N \}.
\end{equation}
First, we establish the following fundamental weighted inequality.

\begin{proposition}\label{Goldbachbasic}
Let 
$$
1+3\varepsilon<a<2,\quad \tau=\dfrac{a-1}{a}-\varepsilon, \quad 0<\alpha<\tau<\frac{1}{2}. 
$$
If $N>\varepsilon^{\frac{1-a}{a\varepsilon}}$ then 
\begin{equation}\label{D1a>}
D_{1,a}(N)\ge\sum\limits_{\substack{n\in \mathscr{A}\\P^{-}(n)\ge N^\alpha}}1
-\sum\limits_{\substack{n\in \mathscr{A}\\P^{-}(n)\ge N^\tau\\\Omega(n)\ge 2}}1
+\sum\limits_{\substack{n\in \mathscr{A}\\P^{-}(n)\ge N^\tau\\\Omega(n)\ge 3}}1
-\sum\limits_{\substack{n\in \mathscr{A}\\P^{-}(n)\ge N^\alpha\\\Omega(n)\ge 3}}1. 
\end{equation}
Moreover, the third term on the right vanishes when $a>1.5+3\varepsilon$.
\end{proposition}
	
\begin{proof} 
Fix an even integer \(N > \varepsilon^{\frac{1-a}{a\varepsilon}} \) and define the weight function  
\[
w(n)=\mathbf{1}_{P^{-}(n)\ge N^\alpha}-\mathbf{1}_{\Omega(n)=2}\mathbf{1}_{P^{-}(n)\ge N^\tau}
+\mathbf{1}_{\Omega(n)\ge 3}\mathbf{1}_{P^{-}(n)\ge N^\tau}
-\mathbf{1}_{\Omega(n)\ge 3}\mathbf{1}_{P^{-}(n)\ge N^\alpha}. 
\]
It suffices to show  
\begin{equation}\label{D1a>wn}
D_{1,a}(n)\ge\sum\limits_{n\in \mathscr{A}}w(n).
\end{equation}

Recall from \eqref{Def/A=Goldbach} that if \( n \in \mathscr{A} \), then \( n > \varepsilon N \).  
We examine \( w(n) \) according to the size of \( P^{-}(n) \). 

		Case 1: $P^{-}(n)<N^{\alpha}$.
		It is easy to see that $w(n)=0$ regardless of the value of $\Omega(n)$.
		
		Case 2: $N^{\alpha}\le P^{-}(n)<N^{\tau}$.
		In this case, $n$ cannot be a prime. When $\Omega(n)=2$, we have $w(n)=1-0+0-0=1$; and when $\Omega(n)\ge 3$, we have $w(n)=1-0+0-1=0$.
		
		Case 3: $P^{-}(n)\ge N^{\tau}$.
		If $n$ is a prime, then $w(n)=1-0+0-0=1$; if $\Omega(n)=2$, we have $w(n)=1-1+0-0=0$; and if $\Omega(n)\ge 3$, we have $w(n)=1-1+1-1=0$.

Therefore, $w(n)=1$ if and only if $n$ is prime or $n=rq$ where $r$ and $q$ are primes satisfying
$$
N^{\alpha}\le r<N^{\tau}, \quad q> \dfrac{\varepsilon N}{N^{\tau}}=\varepsilon N^{1-\tau}. 
$$
In the latter case we have
$r\le q^{a-1}$ when $N>\varepsilon^{\frac{1-a}{a\varepsilon}}$, and \eqref{D1a>wn} follows. 
\end{proof}

Next, we prove the following simple weighted inequality, which will be refined later. 
This is a variant of \cite[Lemma 2.1]{Wu08}; it has a long history of corrections and refinements
 (see for example \cite{CaLu02},\cite{Wu04},\cite{Cai08a},\cite{Wu08}), though originally stated without proof by Chen \cite{Chen78}. 
	
\begin{proposition}\label{Goldbachbig}
Let 
$$
1.5+3\varepsilon<a<2,\quad \tau=\dfrac{a-1}{a}-\varepsilon, \quad \frac{1}{21}<\kappa<\sigma\le \frac{1}{3}. 
$$ 
Then we have 
\begin{equation}\label{Goldbachbigw2}
\begin{split}
2D_{1,a}(N)
\ge\ & 2S_1(\kappa)-2S_2(\tau)-S_3(\kappa,\sigma)-2S_4(\sigma)\\
&-S_5(\kappa,\sigma)+S_6(\kappa,\sigma)+O(N^{1-\kappa}),
\end{split}
\end{equation}
where
\[
S_1(\kappa) := S(\mathscr{A}; \mathscr{P}(N), N^\kappa), 
\]
\[
S_2(\kappa) :=\sum\limits_{\substack{N^\kappa\le p\le N^{1/2} \\(p, N)=1}}S(\mathscr{A}_p ; \mathscr{P}(N), p), 
\]
\[
S_3(\kappa, \sigma) := \sum_{\substack{N^\kappa \le p \le N^\sigma \\ (p, N) = 1}} S(\mathscr{A}_p; \mathscr{P}(N), N^\kappa),
\] 
\[
S_4(\sigma) := \sum_{\substack{N^\sigma \le p_1 \le p_2 \le (N/p_1)^{1/2} \\ (p_1 p_2, N) = 1}} S(\mathscr{A}_{p_1 p_2}; \mathscr{P}(N p_1), p_2), 
\]
\[
S_5(\kappa, \sigma) := \sum_{\substack{N^\kappa \le p_1 \le N^\sigma \le p_2 \le (N/p_1)^{1/2} \\ (p_1 p_2, N) = 1}}  S(\mathscr{A}_{p_1 p_2}; \mathscr{P}(N p_1), p_2),
\]
and 
\[
S_6(\kappa, \sigma) := \sum_{\substack{N^\kappa \le p_1 \le p_2 \le p_3 \le N^\sigma \\ (p_1 p_2 p_3, N) = 1}} S(\mathscr{A}_{p_1 p_2 p_3};  \mathscr{P}(N p_1), p_2).
\]
\end{proposition}
	
\begin{proof}
When 	
$$
1.5+3\varepsilon<a<2,\quad \tau=\dfrac{a-1}{a}-\varepsilon>\frac{1}{3}, \quad \frac{1}{21}<\kappa\le \frac{1}{3}, 
$$
we derive by Proposition \ref{Goldbachbasic} the basic weighted inequality
\begin{equation}\label{Goldbachw1}
\begin{split}
D_{1,a}(N)\ge\ & \sum\limits_{\substack{n\in \mathscr{A}\\P^{-}(n)\ge N^\kappa}}1-\sum\limits_{\substack{n\in \mathscr{A}\\P^{-}(n)\ge N^\tau\\ \Omega(n)\ge 2}}1-\sum\limits_{\substack{n\in \mathscr{A}\\P^{-}(n)\ge N^\kappa\\\Omega(n)\ge 3}}1\\
\ge\ & S(\mathscr{A};\mathscr{P}(N),N^\kappa)-\sum\limits_{\substack{N^\tau\le p\le N^{1/2} \\(p, N)=1}}S(\mathscr{A}_p ; \mathscr{P}(N), p)\\
&-\sum_{\substack{N^\kappa \le p_1 \le p_2 \le (N/p_1)^{1/2} \\ (p_1 p_2, N) = 1}}  
S(\mathscr{A}_{p_1 p_2}; \mathscr{P}(N p_1), p_2)+O(N^{1-\kappa})\\
=\ & S_1(\kappa)-S_2(\tau)-S_4(\kappa)+O(N^{1-\kappa}),
\end{split}	
\end{equation}
where we have used the trivial bound 
\begin{equation}\label{Buchstab1}
\sum_{\substack{N^\kappa \le p \le N^{1/2} \\ (p, N) = 1}} S(\mathscr{A}_{p^2} ; \mathscr{P}(N), p) 
\ll N^{1-\kappa} 
\end{equation}
as well as the trivial bound 
\begin{equation}\label{Buchstab0} 
S(\mathscr{A}_{p} ; \mathscr{P}(N), p) \ll N^{1-\kappa}
\end{equation}
which holds for $N^{\kappa}\le p\le  N^{1/2}$. 
		
Similarly, when $1/21<\sigma\le 1/3$, we have
\begin{equation}\label{Goldbachw2}
\begin{split}
D_{1,a}(N)
\ge\ 
& S(\mathscr{A};\mathscr{P}(N),N^\sigma)-\sum\limits_{\substack{N^\tau\le p\le N^{1/2} \\(p, N)=1}}S(\mathscr{A}_p ; \mathscr{P}(N), p)\\
&-\sum_{\substack{N^\sigma \le p_1 \le p_2 \le (N/p_1)^{1/2} \\ (p_1 p_2, N) = 1}}  S(\mathscr{A}_{p_1 p_2}; \mathscr{P}(N p_1), p_2)+O(N^{1-\sigma})\\
=\ & S_1(\sigma)-S_2(\tau)-S_4(\sigma)+O(N^{1-\sigma}),
		\end{split}	
		\end{equation}
		where we have used the trivial bound
		\begin{equation}\label{Buchstab2}
		\sum_{\substack{N^\sigma \le p \le N^{1/2} \\ (p, N) = 1}} S\left(\mathscr{A}_{p^2} ; \mathscr{P}(N), p\right) \ll N^{1-\sigma}.
		\end{equation}
Noting that $\kappa<\sigma$, we have by \eqref{Goldbachw1} and \eqref{Goldbachw2} that 
\begin{equation}\label{Goldbachw0}
2D_{1,a}(N)\ge S_1(\kappa)+S_1(\sigma)-2S_2(\tau)-S_4(\kappa)-S_4(\sigma)+O(N^{1-\kappa}).
\end{equation}

In the rest of this proof, bounds similar to \eqref{Buchstab1}, \eqref{Buchstab0} and \eqref{Buchstab2} will be applied 
several times. Now $\kappa<\sigma$, we can rewrite $S_4(\kappa)$ as 
\begin{equation}\label{Goldbachw5}
\begin{split}
S_4(\kappa)
=\ &\sum_{\substack{N^\kappa \le p_1 \le p_2 \le (N/p_1)^{1/2} \\ (p_1 p_2, N) = 1}}  S(\mathscr{A}_{p_1 p_2}; \mathscr{P}(N p_1), p_2)\\ 
=\ &\sum_{\substack{N^\kappa \le p_1 \le p_2 \le N^\sigma \\ (p_1 p_2, N) = 1}}  S(\mathscr{A}_{p_1 p_2}; \mathscr{P}(N p_1), p_2)\\
&+\sum_{\substack{N^\kappa \le p_1 \le N^\sigma\le p_2 \le (N/p_1)^{1/2} \\ (p_1 p_2, N) = 1}}  S(\mathscr{A}_{p_1 p_2}; \mathscr{P}(N p_1), p_2)\\
&+\sum_{\substack{N^\sigma \le p_1 \le p_2 \le (N/p_1)^{1/2} \\ (p_1 p_2, N) = 1}}  S(\mathscr{A}_{p_1 p_2}; \mathscr{P}(N p_1), p_2)+O(N^{1-\kappa})\\
=\ &\sum_{\substack{N^\kappa \le p_1 \le p_2 \le N^\sigma \\ (p_1 p_2, N) = 1}}  S(\mathscr{A}_{p_1 p_2}; \mathscr{P}(N p_1), p_2)+S_5(\kappa,\sigma)+S_4(\sigma)+O(N^{1-\kappa}). 
\end{split}
\end{equation}
On the other hand, by Buchstab's identity we can transform $S_1(\sigma)$ into the form 
\begin{equation}\label{Goldbachw3}
\begin{split}
S_1(\sigma)=\ &S(\mathscr{A};\mathscr{P}(N),N^\sigma)\\
		=\ & S(\mathscr{A};\mathscr{P}(N),N^\kappa)-\sum\limits_{\substack{N^\kappa\le p\le N^{\sigma} \\(p, N)=1}}S(\mathscr{A}_p ; \mathscr{P}(N), p)\\
		=\ & S(\mathscr{A};\mathscr{P}(N),N^\kappa)-\sum\limits_{\substack{N^\kappa\le p\le N^{\sigma} \\(p, N)=1}}S(\mathscr{A}_p ; \mathscr{P}(N), N^\kappa)\\
		& +\sum_{\substack{N^\kappa \le p_1 \le p_2 \le N^\sigma \\ (p_1 p_2, N) = 1}}  S(\mathscr{A}_{p_1 p_2}; \mathscr{P}(N), p_1)+O(N^{1-\kappa})\\
		=\ &S_1(\kappa)-S_3(\kappa,\sigma)+\sum_{\substack{N^\kappa \le p_1 \le p_2 \le N^\sigma \\ (p_1 p_2, N) = 1}}  S(\mathscr{A}_{p_1 p_2}; \mathscr{P}(N), p_1)+O(N^{1-\kappa}). 
\end{split}
\end{equation}
We can compute the difference between the two double sums over $p_1$ and $p_2$ in \eqref{Goldbachw5} and in 
\eqref{Goldbachw3}, to get  
\begin{equation}\label{Goldbachw4}
\begin{split}
&\sum_{\substack{N^\kappa \le p_1 \le p_2 \le N^\sigma \\ (p_1 p_2, N) = 1}}  S(\mathscr{A}_{p_1 p_2}; \mathscr{P}(N), p_1)-\sum_{\substack{N^\kappa \le p_1 \le p_2 \le N^\sigma \\ (p_1 p_2, N) = 1}}  S(\mathscr{A}_{p_1 p_2}; \mathscr{P}(Np_1), p_2)\\
& =  \sum_{\substack{N^\kappa \le p_1 \le p_2 \le p_3 \le N^\sigma \\ (p_1 p_2 p_3, N) = 1}} S(\mathscr{A}_{p_1 p_2 p_3}; \mathscr{P}(N p_1), p_2)+O(N^{1-\kappa})\\
& =  S_6(\kappa,\sigma)+O(N^{1-\kappa}). 
\end{split}
\end{equation}
Collecting \eqref{Goldbachw5}, \eqref{Goldbachw3} and \eqref{Goldbachw4} into \eqref{Goldbachw0},  
we get \eqref{Goldbachbigw2} finally.  
\end{proof}
	
Now, we begin by nesting the preceding weighted inequality to obtain a more refined version.
	
\begin{proposition}\label{Goldbachweight1}
Let 
$$
1.5+3\varepsilon<a<2,\quad \tau=\dfrac{a-1}{a}-\varepsilon, \quad \dfrac{1}{18}<\alpha<\beta<\dfrac{1-3\beta}{3}<\gamma< \dfrac{1}{3}<\tau. 
$$ 
Then 
\begin{equation}\label{Goldbachbigw4}
\begin{split}
4D_{1,a}(N) \geq\ & 3G_1 + G_2 -4G_3- G_4- G_5 + G_6 + G_7 \\
&- 2G_8 - G_9 - G_{10} - G_{11} - G_{12} + O(N^{1 - \alpha}),
\end{split}
\end{equation}
where
\[
G_1 := S(\mathscr{A}; \mathscr{P}(N), N^{\alpha}),
\quad 
G_2 := S(\mathscr{A}; \mathscr{P}(N), N^{\beta}),
\]
\[
G_3 :=\sum\limits_{\substack{N^\tau\le p\le N^{1/2} \\(p, N)=1}}S(\mathscr{A}_p ; \mathscr{P}(N), p), 
\quad 
G_4 := \sum_{\substack{N^{\alpha} \leq p \le N^{1/3} \\ (p, N) = 1}} S(\mathscr{A}_p; \mathscr{P}(N), N^{\alpha}), 
\]
\[
G_5 := \sum_{\substack{N^{\alpha} \leq p \le N^{\gamma} \\ (p, N) = 1}} S(\mathscr{A}_p; \mathscr{P}(N), N^{\alpha}),
\quad 
G_6 := \sum_{\substack{N^{\alpha} \leq p_1 \le p_2 \le N^{\beta} \\ (p_1 p_2, N) = 1}} S(\mathscr{A}_{p_1 p_2}; \mathscr{P}(N), N^{\alpha}),
\]
\[
G_7 := \sum_{\substack{N^{\alpha} \leq p_1 \le N^{\beta} \le p_2 \le N^{\gamma} \\ (p_1 p_2, N) = 1}} S(\mathscr{A}_{p_1 p_2}; \mathscr{P}(N), N^{\alpha}),
\]
\[ 
G_8 := \sum_{\substack{N^{\gamma} \leq p_1 \le p_2 \le (N/p_1)^{1/2} \\ (p_1 p_2, N) = 1}} S(\mathscr{A}_{p_1 p_2}; \mathscr{P}(N p_1), p_2),
\]
\[
G_9 := \sum_{\substack{N^{\alpha} \leq p_1 \le N^{1/3} \leq p_2 \le (N/p_1)^{1/2} \\ (p_1 p_2, N) = 1}} S(\mathscr{A}_{p_1 p_2}; \mathscr{P}(N p_1), p_2),
\]
\[
G_{10} := \sum_{\substack{N^{\beta} \leq p_1 \le N^{\gamma} \leq p_2 \le (N/p_1)^{1/2} \\ (p_1 p_2, N) = 1}} S(\mathscr{A}_{p_1 p_2}; \mathscr{P}(N p_1), (N/p_1 p_2)^{1/2}), 
\]
\[
G_{11} := \sum_{\substack{N^{\alpha} \leq p_1 \le p_2 \le p_3 \le p_4 \le N^{\beta} \\ (p_1 p_2 p_3 p_4, N) = 1}} S(\mathscr{A}_{p_1 p_2 p_3 p_4}; \mathscr{P}(Np_1), p_2),
\]
and 
\[
G_{12} := \sum_{\substack{N^{\alpha} \leq p_1 \le p_2 \le p_3 \le N^{\beta} \leq p_4 \le N^{\gamma}\\ (p_1 p_2 p_3 p_4, N) = 1}} S(\mathscr{A}_{p_1 p_2 p_3 p_4}; \mathscr{P}(Np_1), p_2).
\]
\end{proposition}
	
\begin{proof}
First, we apply Proposition \ref{Goldbachbig} with $(\kappa,\sigma,a)=(\alpha,1/3,a)$, to obtain
\begin{equation}\label{GoldbachW1}
\begin{split}
2D_{1,a}(N)
\ge\ & 2S_1(\alpha)-2S_2(\tau)-S_3(\alpha,1/3)-2S_4(1/3)\\
&-S_5(\alpha,1/3)+S_6(\alpha,1/3)+O(N^{1-\alpha})\\
=\ & 2G_1-2G_3-G_4-2S_4(1/3)-G_9+S_6(\alpha,1/3)+O(N^{1-\alpha}).  
\end{split}
\end{equation}
Then we apply Proposition \ref{Goldbachbig} with $(\kappa,\sigma,a)=(\beta,\gamma,a)$, getting 
\begin{equation}\label{GoldbachW2}
\begin{split}
2D_{1,a}(N)
\ge&\  2S_1(\beta)-2S_2(\tau)-S_3(\beta,\gamma)-2S_4(\gamma)\\
&-S_5(\beta,\gamma)+S_6(\beta,\gamma)+O(N^{1-\beta})\\
=&\ 2G_2-2G_3-S_3(\beta,\gamma)-2G_8-S_5(\beta,\gamma)+S_6(\beta,\gamma)+O(N^{1-\beta}). 
\end{split}
\end{equation} 
Adding \eqref{GoldbachW1} and \eqref{GoldbachW2}, 
we have
\begin{equation}\label{GoldbachW3}
\begin{split}
4D_{1,a}(N)
\ge\ & 2G_1+2G_2-4G_3-G_4-2G_8-G_9\\
&-S_3(\beta,\gamma)-S_5(\beta,\gamma)+S_6(\alpha,1/3)+O(N^{1-\alpha}), 
\end{split}
\end{equation}
where we note that $S_4(1/3)=0$ and $S_6(\beta,\gamma)$ is dropped by non-negativity. 

Now we transform the $G_2$ on the right-hand side of \eqref{GoldbachW3}.  
Applying repeatedly Buchstab's identity as well as bounds similar to \eqref{Buchstab1} and \eqref{Buchstab0}, 
we have
\begin{equation*} 
\begin{split}
G_2=\ &S(\mathscr{A};\mathscr{P}(N),N^\beta)\\
		=\ & S(\mathscr{A};\mathscr{P}(N),N^\alpha)-\sum\limits_{\substack{N^\alpha\le p\le N^{\beta} \\(p, N)=1}}S(\mathscr{A}_p ; \mathscr{P}(N), p)\\
=\ & S(\mathscr{A};\mathscr{P}(N),N^\alpha)-\sum\limits_{\substack{N^\alpha\le p\le N^{\beta} \\(p, N)=1}}S(\mathscr{A}_p ; \mathscr{P}(N), N^\alpha)\\
&+\sum_{\substack{N^\alpha \le p_1 \le p_2 \le N^\beta \\ (p_1 p_2, N) = 1}}  S(\mathscr{A}_{p_1 p_2}; \mathscr{P}(N), p_1)+O(N^{1-\alpha}), 
\end{split}
\end{equation*}
and consequently,  
\begin{equation}\label{GoldbachW4}
\begin{split}
G_2
=\ & S(\mathscr{A};\mathscr{P}(N),N^\alpha)-\sum\limits_{\substack{N^\alpha\le p\le N^{\beta} \\(p, N)=1}}S(\mathscr{A}_p ; \mathscr{P}(N), N^\alpha)\\
		&+\sum_{\substack{N^\alpha \le p_1 \le p_2 \le N^\beta \\ (p_1 p_2, N) = 1}}  S(\mathscr{A}_{p_1 p_2}; \mathscr{P}(N), N^\alpha)\\
		&-\sum_{\substack{N^\alpha \le p_1 \le p_2 \le p_3\le N^\beta \\ (p_1 p_2 p_3, N) = 1}}  S(\mathscr{A}_{p_1 p_2 p_3}; \mathscr{P}(N), p_1)+O(N^{1-\alpha})\\
=\ &G_1-G_{13}+G_6-G_{14}+O(N^{1-\alpha}),
		\end{split}
		\end{equation}
where we have written 
		\begin{equation}
		\begin{split}
		G_{13}:=&\sum\limits_{\substack{N^\alpha\le p\le N^{\beta} \\(p, N)=1}}
S(\mathscr{A}_p ; \mathscr{P}(N), N^\alpha),\\
		G_{14}:=&\sum_{\substack{N^\alpha \le p_1 \le p_2 \le p_3\le N^\beta \\ (p_1 p_2 p_3, N) = 1}}  S(\mathscr{A}_{p_1 p_2 p_3}; \mathscr{P}(N), p_1). 
		\end{split}
		\end{equation}
We use \eqref{GoldbachW4} for one of the $G_2$'s on the right-hand side of \eqref{GoldbachW3}, getting 
\begin{equation}\label{GoldbachW5}
\begin{split}
4D_{1,a}(N)
\ge\ & 3G_1+G_2-4G_3-G_4+G_6-2G_8-G_9\\
&-G_{13}-G_{14}-S_3(\beta,\gamma)-S_5(\beta,\gamma)+S_6(\alpha,1/3)+O(N^{1-\alpha}). 
\end{split}
\end{equation}
		
Next we transform the $S_3(\beta,\gamma)$ and $S_5(\beta,\gamma)$ on the right-hand side of \eqref{GoldbachW3}.  
Applying repeatedly Buchstab's identity as well as bounds similar to \eqref{Buchstab1} and \eqref{Buchstab0}, 
we have
\begin{equation*}
\begin{split}
S_3(\beta,\gamma)=\ &\sum_{\substack{N^\beta \le p \le N^\gamma \\ (p, N) = 1}} S(\mathscr{A}_p; \mathscr{P}(N), N^\beta)\\
=\ & \sum_{\substack{N^\beta \le p \le N^\gamma \\ (p, N) = 1}} S(\mathscr{A}_p; \mathscr{P}(N), N^\alpha)-\sum_{\substack{N^{\alpha} \leq p_1 \le N^{\beta} \le p_2 \le N^{\gamma} \\ (p_1 p_2, N) = 1}} S(\mathscr{A}_{p_1 p_2}; \mathscr{P}(N), p_1), 
\end{split}
\end{equation*}
and consequently,  
\begin{equation}\label{GoldbachW6}
\begin{split}
S_3(\beta,\gamma)
=\ & \sum_{\substack{N^\alpha \le p \le N^\gamma \\ (p, N) = 1}} S(\mathscr{A}_p; \mathscr{P}(N), N^\alpha)-\sum_{\substack{N^\alpha \le p \le N^\beta \\ (p, N) = 1}} S(\mathscr{A}_p; \mathscr{P}(N), N^\alpha)\\
&-\sum_{\substack{N^{\alpha} \leq p_1 \le N^{\beta} \le p_2 \le N^{\gamma} \\ (p_1 p_2, N) = 1}} S(\mathscr{A}_{p_1 p_2}; \mathscr{P}(N), N^{\alpha})\\
&+\sum_{\substack{N^{\alpha} \leq p_1 \le p_2 \le N^{\beta} \leq p_3 \le N^{\gamma}\\ (p_1 p_2 p_3, N) = 1}} S(\mathscr{A}_{p_1 p_2 p_3}; \mathscr{P}(N), p_1)+O(N^{1-\alpha})\\
=\ &G_5-G_{13}-G_{7}+G_{15}+O(N^{1-\alpha}),
\end{split}
\end{equation}
where we have written 
\begin{equation}
G_{15}:=\sum_{\substack{N^{\alpha} \leq p_1 \le p_2 \le N^{\beta} \leq p_3 \le N^{\gamma}\\ (p_1 p_2 p_3, N) = 1}} S(\mathscr{A}_{p_1 p_2 p_3}; \mathscr{P}(N), p_1). 
\end{equation}
Inserting \eqref{GoldbachW6} into \eqref{GoldbachW5}, we have
\begin{equation}\label{GoldbachW7}
\begin{split}
4D_{1,a}(N)\ge\ & 3G_1+G_2-4G_3-G_4-G_5+G_6+G_7-2G_8-G_9\\
&-G_{14}-G_{15}-S_5(\beta,\gamma)+S_6(\alpha,1/3)+O(N^{1-\alpha}).
\end{split}
\end{equation}
Using similar arguments to $S_5(\beta,\gamma)$ and noting that $\alpha<\beta$, we have 
\begin{equation}\label{GoldbachW8}
\begin{split}
S_5(\beta,\gamma)=&\sum_{\substack{N^\beta \le p_1 \le N^\gamma \le p_2 \le (N/p_1)^{1/2} \\ (p_1 p_2, N) = 1}}  
S(\mathscr{A}_{p_1 p_2}; \mathscr{P}(N p_1), p_2)\\
=&\ \sum_{\substack{N^{\beta} \leq p_1 \le N^{\gamma} 
\leq p_2 \le (N/p_1)^{1/2} \\ (p_1 p_2, N) = 1}} S(\mathscr{A}_{p_1 p_2}; \mathscr{P}(N p_1), (N/p_1 p_2)^{1/2})\\
&+\sum_{\substack{N^{\beta} \leq p_1 \le N^{\gamma} \leq p_2 \le p_3\le (N/p_1p_2)^{1/2} \\ (p_1 p_2 p_3, N) = 1}} S(\mathscr{A}_{p_1 p_2}; \mathscr{P}(N p_1), p_3)+O(N^{1-\beta})\\
=&\ G_{10}+G_{16}+O(N^{1-\alpha}),
\end{split}
\end{equation}
where 
\begin{equation}\label{def/G16}
G_{16}:=\sum_{\substack{N^{\beta} \leq p_1 \le N^{\gamma} \leq p_2 \le p_3\le (N/p_1p_2)^{1/2} \\ (p_1 p_2 p_3, N) = 1}} S(\mathscr{A}_{p_1 p_2}; \mathscr{P}(N p_1), p_3). 
		\end{equation}
Collecting \eqref{GoldbachW8} into \eqref{GoldbachW7} yields 
\begin{equation}\label{GoldbachW9}
\begin{split}
		4D_{1,a}(N)\ge\ & 3G_1+G_2-4G_3-G_4-G_5+G_6+G_7-2G_8-G_9-G_{10}\\
		&-G_{14}-G_{15}-G_{16}+S_6(\alpha,1/3)+O(N^{1-\alpha}).
\end{split}
\end{equation}

Finally we want to prove  
\begin{equation}\label{GoldbachW10}
S_6(\alpha,1/3)\ge G_{14}+G_{15}+G_{16}-G_{11}-G_{12},
\end{equation}
from which together with \eqref{GoldbachW9} and \eqref{GoldbachW10}, we deduce the desired 
result \eqref{Goldbachbigw4}. To prove \eqref{GoldbachW10}, we apply Buchstab's identity to get 
\begin{equation}\label{GoldbachW11}
\begin{split}
G_{14}-G_{11}
=&\ \sum_{\substack{N^\alpha \le p_1 \le p_2 \le p_3\le N^\beta \\ (p_1 p_2 p_3, N) = 1}}  S(\mathscr{A}_{p_1 p_2 p_3}; \mathscr{P}(N), p_1)\\&-\sum_{\substack{N^\alpha \le p_1 \le p_2 \le p_3\le p_4\le N^\beta \\ (p_1 p_2 p_3 p_4, N) = 1}}  S(\mathscr{A}_{p_1 p_2 p_3}; \mathscr{P}(Np_1), p_2)\\
=&\ \sum_{\substack{N^\alpha \le p_1 \le p_2 \le p_3\le N^\beta \\ (p_1 p_2 p_3, N) = 1}}  S(\mathscr{A}_{p_1 p_2 p_3}; \mathscr{P}(Np_1), p_2)+O(N^{1-\alpha}),
\end{split}
\end{equation}
and also 
\begin{equation}\label{GoldbachW12}
\begin{split}
G_{15}-G_{12}
=&\ \sum_{\substack{N^{\alpha} \leq p_1 \le p_2 \le N^{\beta} \leq p_3 \le N^{\gamma}\\ (p_1 p_2 p_3, N) = 1}} S(\mathscr{A}_{p_1 p_2 p_3}; \mathscr{P}(N), p_1)\\&-\sum_{\substack{N^{\alpha} \leq p_1 \le p_2 \le p_3 \le N^{\beta} \leq p_4 \le N^{\gamma} \\ (p_1 p_2 p_3 p_4, N) = 1}} S(\mathscr{A}_{p_1 p_2 p_3 p_4}; \mathscr{P}(Np_1), p_2)\\
=&\ \sum_{\substack{N^{\alpha} \leq p_1 \le p_2 \le N^{\beta} \leq p_3 \le N^{\gamma}\\ (p_1 p_2 p_3, N) = 1}} S(\mathscr{A}_{p_1 p_2 p_3}; \mathscr{P}(Np_1), p_2)+O(N^{1-\alpha}). 
\end{split}
\end{equation}
The $G_{16}$ on the right-hand side of \eqref{GoldbachW10} can be treated trivially; 
thus by definition \eqref{def/G16} we have 
\begin{equation}\label{GoldbachW13}
G_{16}
\le \sum_{\substack{N^{\beta} \leq p_1 \le N^{\gamma} \leq p_2 \le p_3\le (N/p_1p_2)^{1/2} \\ (p_1 p_2 p_3, N) = 1}} S(\mathscr{A}_{p_1 p_2}; \mathscr{P}(N p_1), p_2). 
\end{equation}
The $S_6$ on the left-hand side of \eqref{GoldbachW10} can be split as 
\begin{equation}\label{GoldbachW14}
\begin{split}
S_6(\alpha, 1/3) 
=&\ \sum_{\substack{N^\alpha \le p_1 \le p_2 \le p_3 \le N^{1/3} \\ (p_1 p_2 p_3, N) = 1}} S(\mathscr{A}_{p_1 p_2 p_3}; \mathscr{P}(N p_1), p_2)\\
		=&\ \sum_{\substack{N^\alpha \le p_1 \le p_2 \le p_3\le N^\beta \\ (p_1 p_2 p_3, N) = 1}}  S(\mathscr{A}_{p_1 p_2 p_3}; \mathscr{P}(Np_1), p_2)\\
		&+\sum_{\substack{N^\alpha \le p_1 \le p_2\le N^\beta \le p_3\le N^{1/3} \\ (p_1 p_2 p_3, N) = 1}}  S(\mathscr{A}_{p_1 p_2 p_3}; \mathscr{P}(Np_1), p_2)\\
		&+\sum_{\substack{N^\alpha \le p_1 \le N^\beta\le p_2\le p_3\le N^{1/3} \\ (p_1 p_2 p_3, N) = 1}}  S(\mathscr{A}_{p_1 p_2 p_3}; \mathscr{P}(Np_1), p_2)\\
		&+\sum_{\substack{N^\beta \le p_1 \le p_2\le p_3\le N^{1/3} \\ (p_1 p_2 p_3, N) = 1}}  S(\mathscr{A}_{p_1 p_2 p_3}; \mathscr{P}(Np_1), p_2),
\end{split}
\end{equation}
where, since $\beta+\gamma>1/3$, the last term satisfies 
\begin{equation}\label{GoldbachW15}
		\begin{split}
		& \sum_{\substack{N^\beta \le p_1 \le p_2\le p_3\le N^{1/3} \\ (p_1 p_2 p_3, N) = 1}}  S(\mathscr{A}_{p_1 p_2 p_3}; \mathscr{P}(Np_1), p_2)\\
		& \ge \sum_{\substack{N^\beta \le p_1\le N^\gamma \le p_2\le p_3\le N^{1/3} \\ (p_1 p_2 p_3, N) = 1}}  S(\mathscr{A}_{p_1 p_2 p_3}; \mathscr{P}(Np_1), p_2)\\
		& \ge \sum_{\substack{N^\beta \le p_1\le N^\gamma \le p_2\le p_3\le (N/p_1p_2)^{1/2} \\ (p_1 p_2 p_3, N) = 1}}  S(\mathscr{A}_{p_1 p_2 p_3}; \mathscr{P}(Np_1), p_2). 
\end{split}
\end{equation}
In conclusion \eqref{GoldbachW10} follows \eqref{GoldbachW11}, \eqref{GoldbachW12}, \eqref{GoldbachW13}, \eqref{GoldbachW14} and \eqref{GoldbachW15}.
\end{proof}
	
\subsection{The Goldbach Conjecture (II)}
In this subsection, we present the weighted sieve method required for Theorem \ref{Goldbach2}, which concerns the conditional version of the Goldbach problem. Now we work with the set
\begin{equation}\label{Def/A=Goldbach2}
\mathscr{A}=\{N-p:p<(1-\varepsilon)N\}, 
\end{equation}
where \(p\) runs over primes. We then prove the following simple weighted inequalities; these will be refined later.

\begin{proposition}\label{Goldbachsmall}
Let 
$$
1.05+3\varepsilon<a<1.5,\quad \tau=\dfrac{a-1}{a}-\varepsilon, \quad 
\frac{1}{21}<\kappa<\sigma\le \tau< \frac{1}{3}. 
$$ 
Then we have
		\begin{equation}\label{Goldbachsmallw2}
		\begin{split}
		2D_{1,a}(N)\ge\ & 2S_1(\kappa)-2S_2(\tau)-S_3(\kappa,\sigma)-S_5(\kappa,\sigma)+S_6(\kappa,\sigma)\\
		&-2S_7(\sigma,\tau)-2S_8(\sigma,\tau)+O(N^{1-\kappa}),
		\end{split}
		\end{equation}
where $S_1, \ldots, S_6$ are as in Proposition~\ref{Goldbachbig}, and 
\begin{equation}
\begin{split}
S_{7}(\sigma,\tau):=&\sum_{\substack{N^\sigma \le p_1 \le p_2 \le N^\tau \\ (p_1 p_2, N) = 1}}  S(\mathscr{A}_{p_1 p_2}; \mathscr{P}(N p_1), p_2), \\
S_{8}(\sigma,\tau):=&\sum_{\substack{N^\sigma \le p_1 \le N^\tau\le p_2 \le (N/p_1)^{1/2} \\ (p_1 p_2, N) = 1}}  S(\mathscr{A}_{p_1 p_2}; \mathscr{P}(N p_1), p_2). 
\end{split}
\end{equation}
\end{proposition}
	
\begin{proof}
When 	
$$
1.05+3\varepsilon<a<1.5,\quad \tau=\dfrac{a-1}{a}-\varepsilon<\frac{1}{3}, 
\quad \frac{1}{21}<\kappa<\sigma\le\tau< \frac{1}{3},
$$
we derive by Proposition \ref{Goldbachbasic} the basic weighted inequality
\begin{equation}\label{Goldbachw7}
\begin{split}
D_{1,a}(N)\ge\ & \sum\limits_{\substack{n\in \mathscr{A}\\P^{-}(n)\ge N^\kappa}}1-\sum\limits_{\substack{n\in \mathscr{A}\\P^{-}(n)\ge N^\tau\\\Omega(n)\ge 2}}1+\sum\limits_{\substack{n\in \mathscr{A}\\P^{-}(n)\ge N^\tau\\\Omega(n)\ge 3}}1-\sum\limits_{\substack{n\in \mathscr{A}\\P^{-}(n)\ge N^\kappa\\\Omega(n)\ge 3}}1\\
\ge\ & S(\mathscr{A};\mathscr{P}(N),N^\kappa)-\sum\limits_{\substack{N^\tau\le p\le N^{1/2} \\(p, N)=1}}S(\mathscr{A}_p ; \mathscr{P}(N), p)\\
&+\sum_{\substack{N^\tau \le p_1 \le p_2 \le (N/p_1)^{1/2} \\ (p_1 p_2, N) = 1}}  S(\mathscr{A}_{p_1 p_2}; \mathscr{P}(N p_1), p_2)\\
&-\sum_{\substack{N^\kappa \le p_1 \le p_2 \le (N/p_1)^{1/2} \\ (p_1 p_2, N) = 1}}  S(\mathscr{A}_{p_1 p_2}; \mathscr{P}(N p_1), p_2)+O(N^{1-\kappa})\\
=\ &S_1(\kappa)-S_2(\tau)+S_4(\tau)-S_4(\kappa)+O(N^{1-\kappa}),
\end{split}	
\end{equation}
where bounds similar to \eqref{Buchstab1} and \eqref{Buchstab0} have been applied. 
Similarly, when $1/21<\sigma\le 1/3$, we have
\begin{equation}\label{Goldbachw8}
\begin{split}
D_{1,a}(N)
\ge\ & S(\mathscr{A};\mathscr{P}(N),N^\sigma)-\sum\limits_{\substack{N^\tau\le p\le N^{1/2} \\(p, N)=1}}S(\mathscr{A}_p ; \mathscr{P}(N), p)\\
&+\sum_{\substack{N^\tau \le p_1 \le p_2 \le (N/p_1)^{1/2} \\ (p_1 p_2, N) = 1}}  S(\mathscr{A}_{p_1 p_2}; \mathscr{P}(N p_1), p_2)\\
&-\sum_{\substack{N^\sigma \le p_1 \le p_2 \le (N/p_1)^{1/2} \\ (p_1 p_2, N) = 1}}  S(\mathscr{A}_{p_1 p_2}; \mathscr{P}(N p_1), p_2)+O(N^{1-\sigma})\\
=\ &S_1(\sigma)-S_2(\tau)+S_4(\tau)-S_4(\sigma)+O(N^{1-\sigma}),
\end{split}	
\end{equation}
where we have used a trivial bound similar to \eqref{Buchstab2}.  
Since $\kappa<\sigma$, by \eqref{Goldbachw7} and \eqref{Goldbachw8} we have
		\begin{equation}\label{Goldbachw9}
2D_{1,a}(N)\ge S_1(\kappa)+S_1(\sigma)-2S_2(\tau)-S_4(\kappa)-S_4(\sigma)+O(N^{1-\kappa}).
		\end{equation}

Next we are going to analyze the terms on the right-hand side. In the following Buchstab's identity as well as 
estimates like \eqref{Buchstab1}, \eqref{Buchstab0} and \eqref{Buchstab2} will be used repeatedly, and this will not be 
pointed out at each circumstances.  The above $S_4(\kappa)$ can be written as, since $\kappa<\sigma$,  
\begin{equation*}
\begin{split}
S_4(\kappa)
=\ &\sum_{\substack{N^\kappa \le p_1 \le p_2 \le (N/p_1)^{1/2} \\ (p_1 p_2, N) = 1}}  S(\mathscr{A}_{p_1 p_2}; \mathscr{P}(N p_1), p_2)\\ 
=\ &\sum_{\substack{N^\kappa \le p_1 \le p_2 \le N^\sigma \\ (p_1 p_2, N) = 1}}  S(\mathscr{A}_{p_1 p_2}; \mathscr{P}(N p_1), p_2)\\
&+\sum_{\substack{N^\kappa \le p_1 \le N^\sigma\le p_2 \le (N/p_1)^{1/2} \\ (p_1 p_2, N) = 1}}  S(\mathscr{A}_{p_1 p_2}; \mathscr{P}(N p_1), p_2)\\
&+\sum_{\substack{N^\sigma \le p_1 \le p_2 \le (N/p_1)^{1/2} \\ (p_1 p_2, N) = 1}}  S(\mathscr{A}_{p_1 p_2}; \mathscr{P}(N p_1), p_2)+O(N^{1-\kappa}), 
\end{split}
\end{equation*}
and therefore, 
\begin{equation}\label{Goldbachw10}
S_4(\kappa)
=\sum_{\substack{N^\kappa \le p_1 \le p_2 \le N^\sigma \\ (p_1 p_2, N) = 1}}  S(\mathscr{A}_{p_1 p_2}; \mathscr{P}(N p_1), p_2)+S_5(\kappa,\sigma)+S_4(\sigma)+O(N^{1-\kappa}). 
\end{equation}
The $S_1(\sigma)$ in \eqref{Goldbachw9} can be written as 
\begin{equation*}
\begin{split}
S_1(\sigma)=\ &S(\mathscr{A};\mathscr{P}(N),N^\sigma)\\
= \ & S(\mathscr{A};\mathscr{P}(N),N^\kappa)-\sum\limits_{\substack{N^\kappa\le p\le N^{\sigma} \\(p, N)=1}}S(\mathscr{A}_p ; \mathscr{P}(N), p)\\
= \ & S(\mathscr{A};\mathscr{P}(N),N^\kappa)-\sum\limits_{\substack{N^\kappa\le p\le N^{\sigma} \\(p, N)=1}}S(\mathscr{A}_p ; \mathscr{P}(N), N^\kappa)\\
&+\sum_{\substack{N^\kappa \le p_1 \le p_2 \le N^\sigma \\ (p_1 p_2, N) = 1}}  S(\mathscr{A}_{p_1 p_2}; \mathscr{P}(N), p_1)+O(N^{1-\kappa}). 
\end{split}
\end{equation*}
It follows that 
\begin{equation}\label{Goldbachw11}
S_1(\sigma) =
S_1(\kappa)-S_3(\kappa,\sigma)+\sum_{\substack{N^\kappa \le p_1 \le p_2 \le N^\sigma \\ (p_1 p_2, N) = 1}}  S(\mathscr{A}_{p_1 p_2}; \mathscr{P}(N), p_1)+O(N^{1-\kappa}). 
\end{equation}
To estimate the double sums over $p_1$ and $p_2$ in the last line, we compute 
\begin{equation}\label{Goldbachw12}
\begin{split}
&\sum_{\substack{N^\kappa \le p_1 \le p_2 \le N^\sigma \\ (p_1 p_2, N) = 1}}  S(\mathscr{A}_{p_1 p_2}; \mathscr{P}(N), p_1)-\sum_{\substack{N^\kappa \le p_1 \le p_2 \le N^\sigma \\ (p_1 p_2, N) = 1}}  S(\mathscr{A}_{p_1 p_2}; \mathscr{P}(Np_1), p_2)\\
&= \sum_{\substack{N^\kappa \le p_1 \le p_2 \le p_3 \le N^\sigma \\ (p_1 p_2 p_3, N) = 1}} S(\mathscr{A}_{p_1 p_2 p_3}; \mathscr{P}(N p_1), p_2)+O(N^{1-\kappa})\\
&= S_6(\kappa,\sigma)+O(N^{1-\kappa}). 
\end{split}
\end{equation}
Collecting \eqref{Goldbachw10}, \eqref{Goldbachw11} and \eqref{Goldbachw12} into \eqref{Goldbachw9}, 
we arrive at 
\begin{equation}\label{Goldbachw13}
\begin{split}
		2D_{1,a}(N)\ge\ & 2S_1(\kappa)-2S_2(\tau)-S_3(\kappa,\sigma)-2S_4(\sigma)\\
		&+2S_4(\tau)-S_5(\kappa,\sigma)+S_6(\kappa,\sigma)+O(N^{1-\kappa}). 
\end{split}
\end{equation}

Finally we estimate the $S_4(\sigma)$ in the last line. 
Since $\kappa<\sigma\le \tau$, we have 
		\begin{equation}\label{Goldbachw14}
		\begin{split}
		S_4(\sigma)=\ &\sum_{\substack{N^\sigma \le p_1 \le p_2 \le (N/p_1)^{1/2} \\ (p_1 p_2, N) = 1}}  S(\mathscr{A}_{p_1 p_2}; \mathscr{P}(N p_1), p_2)\\ 
		=\ &\sum_{\substack{N^\sigma \le p_1 \le p_2 \le N^\tau \\ (p_1 p_2, N) = 1}}  S(\mathscr{A}_{p_1 p_2}; \mathscr{P}(N p_1), p_2)\\
		&+\sum_{\substack{N^\sigma \le p_1 \le N^\tau\le p_2 \le (N/p_1)^{1/2} \\ (p_1 p_2, N) = 1}}  S(\mathscr{A}_{p_1 p_2}; \mathscr{P}(N p_1), p_2)\\
		&+\sum_{\substack{N^\tau \le p_1 \le p_2 \le (N/p_1)^{1/2} \\ (p_1 p_2, N) = 1}}  S(\mathscr{A}_{p_1 p_2}; \mathscr{P}(N p_1), p_2)+O(N^{1-\sigma})\\
		=\ &S_7(\sigma,\tau)+S_8(\sigma,\tau)+S_4(\tau)+O(N^{1-\kappa}). 
		\end{split}
		\end{equation}
Inserting this into \eqref{Goldbachw13} proves the desired formula \eqref{Goldbachsmallw2}. 
\end{proof}

The above weighted inequality is designed for the case $a<1.5$. Compared with Proposition \ref{Goldbachbig} for $a>1.5$, its main advantage lies in the additional term $S_4(\tau)$ that appears when applying Proposition \ref{Goldbachbasic}: although this term cannot be bounded from below via the switching principle, it can be canceled with $S_4(\sigma)$, thereby turning $S_4(\sigma)$ into the smaller $S_7(\sigma,\tau)+S_8(\sigma,\tau)$.
Now, we begin by nesting the preceding weighted inequality to obtain a more refined version. Since the structure of our problem when $a<1.5$ has departed far from the classical case of Chen's $a=2$, the weighted inequality below is completely new and fully tailored to the unique features of the new problem.

\begin{proposition}\label{Goldbachweight2}
Let 
$$
1.05+3\varepsilon<a<1.5,\quad \tau=\dfrac{a-1}{a}-\varepsilon, \quad \dfrac{1}{18}<\alpha<\beta<\gamma<\delta< \tau<\dfrac{1}{3} 
$$ 
such that
\begin{equation}\label{Con/bcd}
\beta+\gamma+2\delta< 1. 
\end{equation}
Then 
\begin{equation}\label{Goldbachsmallw4}
		\begin{split}
		4D_{1,a}(N) \geq\ &3H_1 + H_2 -4H_3- H_4- H_5 + H_6 + H_7- 2H_8 \\
		& - 2H_9 -2H_{10}- 2H_{11} - H_{12} - H_{13}-H_{14}-H_{15} \\ 
		& -H_{16}-H_{17} + O(N^{1 - \alpha}),
		\end{split}
		\end{equation}
		where
\[
H_1 := S(\mathscr{A}; \mathscr{P}(N), N^{\alpha}),\quad H_2 := S(\mathscr{A}; \mathscr{P}(N), N^{\beta}),
\]
\[
H_3 :=\sum\limits_{\substack{N^\tau\le p\le N^{1/2} \\(p, N)=1}}S(\mathscr{A}_p ; \mathscr{P}(N), p),
\quad 
H_4 := \sum_{\substack{N^{\alpha} \leq p \le N^{\delta} \\ (p, N) = 1}} S(\mathscr{A}_p; \mathscr{P}(N), N^{\alpha}), 
\]
\[
H_5 := \sum_{\substack{N^{\alpha} \leq p \le N^{\gamma} \\ (p, N) = 1}} S(\mathscr{A}_p; \mathscr{P}(N), N^{\alpha}),
\]
\[ 
H_6 := \sum_{\substack{N^{\alpha} \leq p_1 \le p_2 \le N^{\beta} \\ (p_1 p_2, N) = 1}} S(\mathscr{A}_{p_1 p_2}; \mathscr{P}(N), N^{\alpha}),
\]
\[ 
H_7 := \sum_{\substack{N^{\alpha} \leq p_1 \le N^{\beta} \le p_2 \le N^{\gamma} \\ (p_1 p_2, N) = 1}} S(\mathscr{A}_{p_1 p_2}; \mathscr{P}(N), N^{\alpha}),
\] 
\[
H_8 := \sum_{\substack{N^{\gamma} \leq p_1 \le p_2 \le N^{\tau} \\ (p_1 p_2, N) = 1}} S(\mathscr{A}_{p_1 p_2}; \mathscr{P}(N p_1), p_2),
\]
\[ 
H_{9} := \sum_{\substack{N^{\delta} \leq p_1 \le p_2 \le N^{\tau} \\ (p_1 p_2, N) = 1}} S(\mathscr{A}_{p_1 p_2}; \mathscr{P}(N p_1), p_2),
\]
\[
H_{10} := \sum_{\substack{N^{\gamma} \leq p_1 \le N^{\tau}\le p_2 \le (N/p_1)^{1/2} \\ (p_1 p_2, N) = 1}} S(\mathscr{A}_{p_1 p_2}; \mathscr{P}(N p_1), p_2),
\]
\[
H_{11} := \sum_{\substack{N^{\delta} \leq p_1 \le N^{\tau}\le p_2 \le (N/p_1)^{1/2} \\ (p_1 p_2, N) = 1}} S(\mathscr{A}_{p_1 p_2}; \mathscr{P}(N p_1), p_2),
\] 
\[
H_{12} := \sum_{\substack{N^{\alpha} \leq p_1 \le N^{\delta}\leq p_2 \le (N/p_1)^{1/2} \\ (p_1 p_2, N) = 1}} S(\mathscr{A}_{p_1 p_2}; \mathscr{P}(N p_1), p_2),
\] 
\[
H_{13} := \sum_{\substack{N^{\beta} \leq p_1 \le N^{\gamma} \leq p_2 \le (N/p_1)^{1/2} \\ (p_1 p_2, N) = 1}} S(\mathscr{A}_{p_1 p_2}; \mathscr{P}(N p_1), (N/p_1 p_2)^{1/2}),
\]
\[
H_{14} := \sum_{\substack{N^\beta \le p_1\le N^{\gamma}\\N^\delta \le p_2\le p_3\le (N/p_1p_2)^{1/2} \\ (p_1 p_2 p_3, N) = 1}}  S(\mathscr{A}_{p_1 p_2 p_3}; \mathscr{P}(Np_1), p_2),
\]
\[ 
H_{15} := \sum_{\substack{N^\beta \le p_1\le N^\gamma \le p_2\le N^{\delta}\le p_3\le (N/p_1p_2)^{1/2} \\ (p_1 p_2 p_3, N) = 1}}  S(\mathscr{A}_{p_1 p_2 p_3}; \mathscr{P}(Np_1), p_2),
\]
\[
H_{16} := \sum_{\substack{N^{\alpha} \leq p_1 \le p_2 \le p_3 \le p_4 \le N^{\beta} \\ (p_1 p_2 p_3 p_4, N) = 1}} S(\mathscr{A}_{p_1 p_2 p_3 p_4}; \mathscr{P}(Np_1), p_2),
\]
and 
\[
H_{17} := \sum_{\substack{N^{\alpha} \leq p_1 \le p_2 \le p_3 \le N^{\beta} \leq p_4 \le N^{\gamma}\\ (p_1 p_2 p_3 p_4, N) = 1}} S(\mathscr{A}_{p_1 p_2 p_3 p_4}; \mathscr{P}(Np_1), p_2).
\] 
\end{proposition}
	
\begin{proof}
First, we apply Proposition \ref{Goldbachsmall} with $(\kappa,\sigma,a)=(\alpha,\delta,a)$ to obtain
\begin{equation}\label{GoldW1}
\begin{split}
2D_{1,a}(N)\ge\ & 2S_1(\alpha)-2S_2(\tau)-S_3(\alpha,\delta)-S_5(\alpha,\delta)+S_6(\alpha,\delta)\\
&-2S_7(\delta,\tau)-2S_8(\delta,\tau)+O(N^{1-\alpha})\\
=\ & 2H_1-2H_3-H_4-H_{12}+S_6(\alpha,\delta)\\ 
& -2H_9-2H_{11}+O(N^{1-\alpha}). 
\end{split}
\end{equation}
Another application of Proposition \ref{Goldbachsmall} with $(\kappa,\sigma,a)=(\beta,\gamma,a)$ gives 
\begin{equation}\label{GoldW2}
\begin{split}
2D_{1,a}(N)
\ge\ & 2S_1(\beta)-2S_2(\tau)-S_3(\beta,\gamma)-S_5(\beta,\gamma)+S_6(\beta,\gamma)\\
& -2S_7(\gamma,\tau)-2S_8(\gamma,\tau)+O(N^{1-\beta})\\
=\ & 2H_2-2H_3-S_3(\beta,\gamma)-S_5(\beta,\gamma)+S_6(\beta,\gamma) \\ 
& -2H_8-2H_{10}+O(N^{1-\beta}).
\end{split}
\end{equation}
Adding \eqref{GoldW1} and \eqref{GoldW2}, we get 
\begin{equation}\label{GoldW3}
\begin{split}
4D_{1,a}(N)\ge\ & 2H_1+2H_2-4H_3-H_4-2H_8-2H_9-2H_{10}-2H_{11}\\
& -H_{12} -S_3(\beta,\gamma)-S_5(\beta,\gamma)+S_6(\alpha,\delta)+O(N^{1-\alpha}),
\end{split}
\end{equation}
where we have used the assumption $\alpha<\beta$, and have dropped the term $S_6(\beta,\gamma)$ by non-negativity. 

Next we are going to analyze the terms on the right-hand side. In the following Buchstab's identity as well as 
estimates like \eqref{Buchstab1}, \eqref{Buchstab0} and \eqref{Buchstab2} will be used repeatedly, and this will not be 
pointed out at each circumstances. 
The $H_2$ on the right-hand side of \eqref{GoldW3} can be written as
\begin{equation*}
\begin{split}
H_2=\ &S(\mathscr{A};\mathscr{P}(N),N^\beta)\\
=\ & S(\mathscr{A};\mathscr{P}(N),N^\alpha)-\sum\limits_{\substack{N^\alpha\le p\le N^{\beta} \\(p, N)=1}}S(\mathscr{A}_p ; \mathscr{P}(N), p)\\
=\ & S(\mathscr{A};\mathscr{P}(N),N^\alpha)
-\sum\limits_{\substack{N^\alpha\le p\le N^{\beta} \\(p, N)=1}}S(\mathscr{A}_p; \mathscr{P}(N), N^\alpha) \\ 
&+\sum_{\substack{N^\alpha \le p_1 \le p_2 \le N^\beta \\ (p_1 p_2, N) = 1}}  S(\mathscr{A}_{p_1 p_2}; \mathscr{P}(N), p_1)+O(N^{1-\alpha}), 
\end{split}
\end{equation*}
and this can be further transform as 
\begin{equation}\label{GoldW4}
\begin{split}
H_2=\ & S(\mathscr{A};\mathscr{P}(N),N^\alpha)-\sum\limits_{\substack{N^\alpha\le p\le N^{\beta} \\(p, N)=1}}S(\mathscr{A}_p ; \mathscr{P}(N), N^\alpha)\\
		&+\sum_{\substack{N^\alpha \le p_1 \le p_2 \le N^\beta \\ (p_1 p_2, N) = 1}}  S(\mathscr{A}_{p_1 p_2}; \mathscr{P}(N), N^\alpha)\\
		&-\sum_{\substack{N^\alpha \le p_1 \le p_2 \le p_3\le N^\beta \\ (p_1 p_2 p_3, N) = 1}}  S(\mathscr{A}_{p_1 p_2 p_3}; \mathscr{P}(N), p_1)+O(N^{1-\alpha})\\
=\ &H_1-H_{18}+H_6-H_{19}+O(N^{1-\alpha}),
\end{split}
\end{equation}
where 
\begin{equation}
\begin{split}
H_{18}:=\ &\sum\limits_{\substack{N^\alpha\le p\le N^{\beta} \\(p, N)=1}}S(\mathscr{A}_p ; \mathscr{P}(N), N^\alpha),\\
H_{19}:=\ &\sum_{\substack{N^\alpha \le p_1 \le p_2 \le p_3\le N^\beta \\ (p_1 p_2 p_3, N) = 1}}  S(\mathscr{A}_{p_1 p_2 p_3}; \mathscr{P}(N), p_1). 
\end{split}
\end{equation}
We replace one of the $H_2$'s on the right-hand side of \eqref{GoldW3} by \eqref{GoldW4}, so that 
 \begin{equation}\label{GoldW5}
\begin{split}
4D_{1,a}(N)\ge\ & 3H_1+H_2-4H_3-H_4+H_6-2H_8\\
& -2H_9-2H_{10}-2H_{11} -H_{12} -H_{18}-H_{19} \\ 
& -S_3(\beta,\gamma)-S_5(\beta,\gamma)+S_6(\alpha,\delta)+O(N^{1-\alpha}). 
\end{split}
\end{equation}

For the term $S_3(\beta,\gamma)$ on the right-hand side, 
we have that 
\begin{equation*}
\begin{split}
S_3(\beta,\gamma)=\ &\sum_{\substack{N^\beta \le p \le N^\gamma \\ (p, N) = 1}} S(\mathscr{A}_p; \mathscr{P}(N), N^\beta)\\
=\ & \sum_{\substack{N^\beta \le p \le N^\gamma \\ (p, N) = 1}} S(\mathscr{A}_p; \mathscr{P}(N), N^\alpha)-\sum_{\substack{N^{\alpha} \leq p_1 \le N^{\beta} \le p_2 \le N^{\gamma} \\ (p_1 p_2, N) = 1}} S(\mathscr{A}_{p_1 p_2}; \mathscr{P}(N), p_1), 
\end{split}
\end{equation*}
and hence 
\begin{equation}\label{GoldW6}
\begin{split}
S_3(\beta,\gamma)
=\ & \sum_{\substack{N^\alpha \le p \le N^\gamma \\ (p, N) = 1}} S(\mathscr{A}_p; \mathscr{P}(N), N^\alpha)-\sum_{\substack{N^\alpha \le p \le N^\beta \\ (p, N) = 1}} S(\mathscr{A}_p; \mathscr{P}(N), N^\alpha)\\
		&-\sum_{\substack{N^{\alpha} \leq p_1 \le N^{\beta} \le p_2 \le N^{\gamma} \\ (p_1 p_2, N) = 1}} S(\mathscr{A}_{p_1 p_2}; \mathscr{P}(N), N^{\alpha})\\
		&+\sum_{\substack{N^{\alpha} \leq p_1 \le p_2 \le N^{\beta} \leq p_3 \le N^{\gamma}\\ (p_1 p_2 p_3, N) = 1}} S(\mathscr{A}_{p_1 p_2 p_3}; \mathscr{P}(N), p_1)+O(N^{1-\alpha})\\
		=\ &H_5-H_{18}-H_{7}+H_{20}+O(N^{1-\alpha}),
\end{split}
\end{equation}
where
\begin{equation}
H_{20}:=\sum_{\substack{N^{\alpha} \leq p_1 \le p_2 \le N^{\beta} \leq p_3 \le N^{\gamma}\\ (p_1 p_2 p_3, N) = 1}} 
S(\mathscr{A}_{p_1 p_2 p_3}; \mathscr{P}(N), p_1). 
\end{equation}
Inserting \eqref{GoldW6} into \eqref{GoldW5} yields 
\begin{equation}\label{GoldW7}
\begin{split}
4D_{1,a}(N)
\ge\ & 3H_1+H_2-4H_3-H_4-H_5+H_6+H_7 -2H_8\\ 
& -2H_9-2H_{10}-2H_{11} -H_{12}-H_{19}-H_{20} \\ 
& -S_5(\beta,\gamma)+S_6(\alpha,\delta)+O(N^{1-\alpha}).
\end{split}
\end{equation}

For the term $S_5(\beta,\gamma)$ on the right-hand side, 
we have similarly that 
		\begin{equation}\label{GoldW8}
		\begin{split}
		S_5(\beta,\gamma)=\ &\sum_{\substack{N^\beta \le p_1 \le N^\gamma \le p_2 \le (N/p_1)^{1/2} \\ (p_1 p_2, N) = 1}}  S(\mathscr{A}_{p_1 p_2}; \mathscr{P}(N p_1), p_2)\\
		=\ & \sum_{\substack{N^{\beta} \leq p_1 \le N^{\gamma} \leq p_2 \le (N/p_1)^{1/2} \\ (p_1 p_2, N) = 1}} S(\mathscr{A}_{p_1 p_2}; \mathscr{P}(N p_1), (N/p_1 p_2)^{1/2})\\
		&+\sum_{\substack{N^{\beta} \leq p_1 \le N^{\gamma} \leq p_2 \le p_3\le (N/p_1p_2)^{1/2} \\ (p_1 p_2 p_3, N) = 1}} S(\mathscr{A}_{p_1 p_2}; \mathscr{P}(N p_1), p_3)+O(N^{1-\beta})\\
		=\ &H_{13}+H_{21}+O(N^{1-\alpha}),
		\end{split}
		\end{equation}
		where 
		\begin{equation}
		H_{21}:=\sum_{\substack{N^{\beta} \leq p_1 \le N^{\gamma} \leq p_2 \le p_3\le (N/p_1p_2)^{1/2} \\ (p_1 p_2 p_3, N) = 1}} S(\mathscr{A}_{p_1 p_2}; \mathscr{P}(N p_1), p_3). 
		\end{equation}
Inserting \eqref{GoldW8} into \eqref{GoldW7}, we have
\begin{equation}\label{GoldW9}
\begin{split}
4D_{1,a}(N)
\ge\ & 3H_1+H_2-4H_3-H_4-H_5+H_6+H_7-2H_8\\
&-2H_9-2H_{10}-2H_{11}-H_{12}-H_{13}-H_{19}-H_{20}-H_{21} \\ 
& +S_6(\alpha,\delta)+O(N^{1-\alpha}).
\end{split}
\end{equation}
		
For the term $S_6(\alpha, \delta)$ on the right-hand side, 
we have similarly that 
		\begin{equation}\label{GoldW10}
		\begin{split}
		S_6(\alpha, \delta) =& \sum_{\substack{N^\alpha \le p_1 \le p_2 \le p_3 \le N^{\delta} \\ (p_1 p_2 p_3, N) = 1}} S(\mathscr{A}_{p_1 p_2 p_3}; \mathscr{P}(N p_1), p_2)\\
		=&\sum_{\substack{N^\alpha \le p_1 \le p_2 \le p_3\le N^\beta \\ (p_1 p_2 p_3, N) = 1}}  S(\mathscr{A}_{p_1 p_2 p_3}; \mathscr{P}(Np_1), p_2)\\
		&+\sum_{\substack{N^\alpha \le p_1 \le p_2\le N^\beta \le p_3\le N^{\delta} \\ (p_1 p_2 p_3, N) = 1}}  S(\mathscr{A}_{p_1 p_2 p_3}; \mathscr{P}(Np_1), p_2)\\
		&+\sum_{\substack{N^\alpha \le p_1 \le N^\beta\le p_2\le p_3\le N^{\delta} \\ (p_1 p_2 p_3, N) = 1}}  S(\mathscr{A}_{p_1 p_2 p_3}; \mathscr{P}(Np_1), p_2)\\
		&+\sum_{\substack{N^\beta \le p_1 \le p_2\le p_3\le N^{\delta} \\ (p_1 p_2 p_3, N) = 1}}  S(\mathscr{A}_{p_1 p_2 p_3}; \mathscr{P}(Np_1), p_2). 
\end{split}
\end{equation}
The next-to-last term is dropped by non-negativity. The last term above is, by \eqref{Con/bcd}, 
\begin{equation}\label{GoldW11}
		\begin{split}
		& \ge \sum_{\substack{N^\beta \le p_1\le N^\gamma \le p_2\le p_3\le N^{\delta} \\ (p_1 p_2 p_3, N) = 1}}  S(\mathscr{A}_{p_1 p_2 p_3}; \mathscr{P}(Np_1), p_2)\\
		& = \sum_{\substack{N^\beta \le p_1\le N^\gamma \le p_2\le p_3\le (N/p_1p_2)^{1/2} \\ (p_1 p_2 p_3, N) = 1}}  S(\mathscr{A}_{p_1 p_2 p_3}; \mathscr{P}(Np_1), p_2)\\
		 &\quad  -\sum_{\substack{N^\beta \le p_1\le N^{\gamma}\\N^\delta \le p_2\le p_3\le (N/p_1p_2)^{1/2} \\ (p_1 p_2 p_3, N) = 1}}  S(\mathscr{A}_{p_1 p_2 p_3}; \mathscr{P}(Np_1), p_2)\\
		&\quad - \sum_{\substack{N^\beta \le p_1\le N^\gamma \le p_2\le N^{\delta}\le p_3\le (N/p_1p_2)^{1/2} \\ (p_1 p_2 p_3, N) = 1}}  S(\mathscr{A}_{p_1 p_2 p_3}; \mathscr{P}(Np_1), p_2)\\
		&= H_{22}-H_{14}-H_{15}, 
\end{split}
\end{equation}
where
\begin{equation}
H_{22}:=\sum_{\substack{N^{\beta} \leq p_1 \le N^{\gamma} \leq p_2 \le p_3\le (N/p_1p_2)^{1/2} \\ (p_1 p_2 p_3, N) = 1}} S(\mathscr{A}_{p_1 p_2}; \mathscr{P}(N p_1), p_2). 
\end{equation}
It follows from from \eqref{GoldW10} and \eqref{GoldW11} that 
\begin{equation}\label{GoldW12}
			S_6(\alpha,\delta)\ge H_{23}+H_{24}+H_{22}-H_{14}-H_{15},
		\end{equation}
		where
		\begin{equation}
		\begin{split}
		H_{23}:=&\sum_{\substack{N^\alpha \le p_1 \le p_2 \le p_3\le N^\beta \\ (p_1 p_2 p_3, N) = 1}}  S(\mathscr{A}_{p_1 p_2 p_3}; \mathscr{P}(Np_1), p_2),\\
		H_{24}:=&\sum_{\substack{N^\alpha \le p_1 \le p_2\le N^\beta \le p_3\le N^{\delta} \\ (p_1 p_2 p_3, N) = 1}}  S(\mathscr{A}_{p_1 p_2 p_3}; \mathscr{P}(Np_1), p_2).
		\end{split}
		\end{equation}
Inserting \eqref{GoldW12} into \eqref{GoldW9}, one gets 
		\begin{equation}\label{GoldW13}
		\begin{split}
		4D_{1,a}(N)
		\ge\ 
		& 3H_1+H_2-4H_3-H_4-H_5+H_6+H_7-2H_8\\
		&-2H_9-2H_{10}-2H_{11}-H_{12}-H_{13}-H_{14}-H_{15}\\
		&-H_{19}-H_{20}-H_{21}+H_{22}+H_{23}+H_{24}+O(N^{1-\alpha}).
		\end{split}
		\end{equation}
		
Now we only need to prove  
\begin{equation}\label{GoldW14}
H_{22}+H_{23}+H_{24}\ge H_{19}+H_{20}+H_{21}-H_{16}-H_{17},
\end{equation}
so that \eqref{Goldbachsmallw4} holds by \eqref{GoldW13} and \eqref{GoldW14}.
		
By Buchstab's identity, we have
\begin{equation}\label{GoldW15}
\begin{split}
H_{19}-H_{16}
=\ &\sum_{\substack{N^\alpha \le p_1 \le p_2 \le p_3\le N^\beta \\ (p_1 p_2 p_3, N) = 1}}  S(\mathscr{A}_{p_1 p_2 p_3}; \mathscr{P}(N), p_1)\\&-\sum_{\substack{N^\alpha \le p_1 \le p_2 \le p_3\le p_4\le N^\beta \\ (p_1 p_2 p_3 p_4, N) = 1}}  S(\mathscr{A}_{p_1 p_2 p_3}; \mathscr{P}(Np_1), p_2)\\
=\ &\sum_{\substack{N^\alpha \le p_1 \le p_2 \le p_3\le N^\beta \\ (p_1 p_2 p_3, N) = 1}}  S(\mathscr{A}_{p_1 p_2 p_3}; \mathscr{P}(Np_1), p_2)+O(N^{1-\alpha})\\
=\ &H_{23}+O(N^{1-\alpha}). 
\end{split}
\end{equation}
Similarly, 
\begin{equation}\label{GoldW16}
\begin{split}
H_{20}-H_{17}
=\ &\sum_{\substack{N^{\alpha} \leq p_1 \le p_2 \le N^{\beta} \leq p_3 \le N^{\gamma}\\ (p_1 p_2 p_3, N) = 1}} S(\mathscr{A}_{p_1 p_2 p_3}; \mathscr{P}(N), p_1)\\&-\sum_{\substack{N^{\alpha} \leq p_1 \le p_2 \le p_3 \le N^{\beta} \leq p_4 \le N^{\gamma} \\ (p_1 p_2 p_3 p_4, N) = 1}} S(\mathscr{A}_{p_1 p_2 p_3 p_4}; \mathscr{P}(Np_1), p_2)\\
=\ &\sum_{\substack{N^{\alpha} \leq p_1 \le p_2 \le N^{\beta} \leq p_3 \le N^{\gamma}\\ (p_1 p_2 p_3, N) = 1}} S(\mathscr{A}_{p_1 p_2 p_3}; \mathscr{P}(Np_1), p_2)+O(N^{1-\alpha})\\
\le \ &H_{24}+O(N^{1-\alpha}). 
\end{split}
\end{equation}
Trivially, 
\begin{equation}\label{GoldW17}
\begin{split}
H_{21}=\sum_{\substack{N^{\beta} \leq p_1 \le N^{\gamma} \leq p_2 \le p_3\le (N/p_1p_2)^{1/2} \\ (p_1 p_2 p_3, N) = 1}} S(\mathscr{A}_{p_1 p_2}; \mathscr{P}(N p_1), p_3)\le H_{22}. 
\end{split}
\end{equation}
Now \eqref{GoldW14} follows from \eqref{GoldW15}, \eqref{GoldW16} and \eqref{GoldW17}.
\end{proof}
	
\subsection{The Twin Prime Conjecture}
In this subsection, we will present the weighted sieve method required for Theorem \ref{Twin1}, which correspond to the unconditional 
version of the Twin Prime Conjecture. Recall that, for $1\leq a\leq 2$,   
\begin{equation}\label{pi1a}
\pi_{1,a}(x)=|\{p\le x: p+2=rq,r\le q^{a-1}\}|.
\end{equation}
where $r$ is either a prime or $1$, and $p$ and $q$ are two primes. 
Now we put  
\begin{equation}\label{Def/A=Twin}
\mathscr{A}=\{p+2:\varepsilon x<p<x\}.
\end{equation}
Parallel to Proposition \ref{Goldbachbasic}, we have the following result. 

\begin{proposition}\label{Twinbasic}
Let 
$$
1+3\varepsilon<a<2,\quad \tau=\dfrac{a-1}{a}-\varepsilon, \quad 0<\alpha<\tau<\frac{1}{2}. 
$$
If $N>\varepsilon^{\frac{1-a}{a\varepsilon}}$ then 
\begin{equation}\label{pi1a>}
\pi_{1,a}(x)\ge\sum\limits_{\substack{n\in \mathscr{A}\\P^{-}(n)\ge x^\alpha}}
1-\sum\limits_{\substack{n\in \mathscr{A}\\P^{-}(n)\ge x^\tau\\\Omega(n)\ge 2}}1
+\sum\limits_{\substack{n\in \mathscr{A}\\P^{-}(n)\ge x^\tau\\\Omega(n)\ge 3}}1
-\sum\limits_{\substack{n\in \mathscr{A}\\P^{-}(n)\ge x^\alpha\\\Omega(n)\ge 3}}1. 
\end{equation}
Moreover, the third term on the right vanishes when $a>1.5+3\varepsilon$.
\end{proposition}

\begin{proof}
For \(x > \varepsilon^{\frac{1-a}{a\varepsilon}} \) we define the weight function  
\[
w(n,x)
=\mathbf{1}_{P^{-}(n)\ge x^\alpha}-\mathbf{1}_{\Omega(n)=2}\mathbf{1}_{P^{-}(n)\ge x^\tau}
+\mathbf{1}_{\Omega(n)\ge 3}\mathbf{1}_{P^{-}(n)\ge x^\tau}-\mathbf{1}_{\Omega(n)\ge 3}\mathbf{1}_{P^{-}(n)\ge x^\alpha}. 
\]
It suffices to show  
\begin{equation}\label{pi1a>wn}
\pi_{1,a}(x)\ge\sum\limits_{n\in \mathscr{A}}w(n,x).
\end{equation}

Recall from \eqref{Def/A=Twin} that if $n\in \mathscr{A}$, then $n>\varepsilon x$.  
We examine \( w(n, x) \) according to the size of \( P^{-}(n) \). 

Case 1: $P^{-}(n)<x^{\alpha}$. 
It is easy to see that $w(n,x)=0$ regardless of the value of $\Omega(n)$.
		
Case 2: $x^{\alpha}\le P^{-}(n)<x^{\tau}$.
In this case, $n$ cannot be a prime. When $\Omega(n)=2$, we have $w(n,x)=1-0+0-0=1$; and when $\Omega(n)\ge 3$, we have $w(n,x)=1-0+0-1=0$.
		
Case 3: $P^{-}(n)\ge x^{\tau}$.
If $n$ is a prime, then $w(n,x)=1-0+0-0=1$; if $\Omega(n)=2$, we have $w(n,x)=1-1+0-0=0$; 
and if $\Omega(n)\ge 3$, we have $w(n,x)=1-1+1-1=0$.

Therefore, $w(n,x)=1$ if and only if $n$ is prime or $n=rq$ where $r$ and $q$ are primes satisfying
$$
x^{\alpha}\le r<x^{\tau}, \quad q> \dfrac{\varepsilon x}{x^{\tau}}=\varepsilon x^{1-\tau}. 
$$
In the latter case we have
$r\le q^{a-1}$ when $x>\varepsilon^{\frac{1-a}{a\varepsilon}}$, and 
\eqref{pi1a>wn} holds.
\end{proof}

By replacing Proposition \ref{Goldbachbasic} with Proposition \ref{Twinbasic} and then following the first optimization in Proposition \ref{Goldbachbig} and the second optimization in Proposition \ref{Goldbachweight1}, we obtain the following weighted inequality. The proof is omitted.

\begin{proposition}\label{Twinweight}
Let 
$$
1.5+3\varepsilon<a<2,\quad \tau=\dfrac{a-1}{a}-\varepsilon, \quad \dfrac{1}{18}<\alpha<\beta<\dfrac{1-3\beta}{3}<\gamma< \dfrac{1}{3}<\tau.  
$$ 
Then
\begin{equation}\label{Twinweight4}
\begin{split}
4\pi_{1,a}(x) 
\ge\ & 3T_1 + T_2 -4T_3- T_4- T_5 + T_6 + T_7 \\
&- 2T_8 - T_9 - T_{10} - T_{11} - T_{12} + O(x^{1 - \alpha}),
\end{split}
\end{equation}
where
\[
T_1 := S(\mathscr{A}; \mathscr{P}, x^{\alpha}), \quad 
T_2 := S(\mathscr{A}; \mathscr{P}, x^{\beta}),
\] 
\[
T_3 :=\sum\limits_{x^\tau\le p\le x^{1/2}}S(\mathscr{A}_p ; \mathscr{P}, p), 
\quad 
T_4 := \sum_{x^{\alpha} \leq p \le x^{1/3}} S(\mathscr{A}_p; \mathscr{P}, x^{\alpha}),
\]
\[
T_5 := \sum_{x^{\alpha} \leq p \le x^{\gamma}} S(\mathscr{A}_p; \mathscr{P}, x^{\alpha}), 
\quad 
T_6 := \sum_{x^{\alpha} \leq p_1 \le p_2 \le x^{\beta}} S(\mathscr{A}_{p_1 p_2}; \mathscr{P}, x^{\alpha}), 
\] 
\[
T_7 := \sum_{x^{\alpha} \leq p_1 \le x^{\beta} \le p_2 \le x^{\gamma}} S(\mathscr{A}_{p_1 p_2}; \mathscr{P}, x^{\alpha}), 
\]
\[
T_8 := \sum_{x^{\gamma} \leq p_1 \le p_2 \le (x/p_1)^{1/2}} S(\mathscr{A}_{p_1 p_2}; \mathscr{P}(p_1), p_2),
\]
\[
T_9 := \sum_{x^{\alpha} \leq p_1 \le x^{1/3} \leq p_2 \le (x/p_1)^{1/2}} S(\mathscr{A}_{p_1 p_2}; \mathscr{P}(p_1), p_2),
\] 
\[ 
T_{10} := \sum_{x^{\beta} \leq p_1 \le x^{\gamma} \leq p_2 \le (x/p_1)^{1/2}} S(\mathscr{A}_{p_1 p_2}; \mathscr{P}( p_1), (x/p_1 p_2)^{1/2}), 
\] 
\[
T_{11} := \sum_{x^{\alpha} \leq p_1 \le p_2 \le p_3 \le p_4 \le x^{\beta} } S(\mathscr{A}_{p_1 p_2 p_3 p_4}; \mathscr{P}(p_1), p_2),
\] 
and 
\[
T_{12} := \sum_{x^{\alpha} \leq p_1 \le p_2 \le p_3 \le x^{\beta} \leq p_4 \le x^{\gamma}} S(\mathscr{A}_{p_1 p_2 p_3 p_4}; \mathscr{P}(p_1), p_2).
\] 
\end{proposition}

\section{Proof of Theorem~\ref{Goldbach1}}\label{sec5}

In Proposition \ref{Goldbachweight1} we specify $a=1.9$ and 
\begin{equation}\label{Cond/Goldbig}
\alpha=4/53,\quad \beta=4/33,\quad \gamma=3/11,\quad \tau=9/19-\varepsilon,\quad N>\exp(\exp(10\varepsilon^{-3})). 
\end{equation} 
Thus we deduce from \eqref{Goldbachbigw4} that 
\begin{equation}\label{GoldbachBigw4}
	\begin{split}
	4D_{1,1.9}(N) 
	\geq\ &3G_1 + G_2 -4G_3- G_4- G_5 + G_6 + G_7 \\
	&- 2G_8 - G_9 - G_{10} - G_{11} - G_{12} + O(N^{49/53}),
	\end{split}
	\end{equation}
where $G_1, \ldots, G_{12}$ are as in Proposition \ref{Goldbachweight1}. We are going to estimate $G_1, \ldots, G_{12}$ 
in the subsections below. 
	
\subsection{Evaluation of $G_1, G_2$}
In this subsection we will obtain lower bounds for $G_1$ and $G_2$ respectively, 
where
\[
G_1 = S(\mathscr{A}; \mathscr{P}(N), N^{4/53}), 
\quad 
G_2 = S(\mathscr{A}; \mathscr{P}(N), N^{4/33}).
\] 
	
We now apply Lemma \ref{Sieve-WF} to estimate \(G_1\).  
Recall from \eqref{Def/A=Goldbach} that  
\begin{equation}\label{New/A=}
\mathscr{A}=\{N-p:p<(1-\varepsilon)N\}. 
\end{equation}
Before applying the lemma, we must verify that conditions \eqref{def/SAPz}--\eqref{RIcon3} hold for this \(\mathscr{A}\).  
First, we compute  
\begin{equation}\label{G1X}
X=\Li(N-\varepsilon N)=\dfrac{N}{\log N}\bigg(1-\varepsilon+O\bigg(\dfrac{1}{\log N}\bigg)\bigg), 
\end{equation}
which gives the main term in the sieve. 
We also have to compute the functions $\omega(p)$ and $V(z)$ defined as in \eqref{RIcon1} and \eqref{def/Vz}, respectively,  
associated to this sequence $\mathscr{A}$. 
These are classical in the study of the Goldbach problem; for example, \cite[\S1.3, Example 5]{HR74} shows that  
\begin{equation}\label{New/ome/p}
\omega(p)= \frac{p}{p-1}
\end{equation}
for odd prime $p$, and $\omega(2)=0$. Thus \eqref{RIcon2} is satisfied. 
It is also proved there that 
\begin{equation}\label{Vz}
V(z) = C(\omega) \frac{e^{-\gamma}}{\log z} \bigg(1 + O\bigg(\frac{1}{\log z}\bigg)\bigg), 
\end{equation}
with 
\begin{equation}\label{Come/G1}
\begin{aligned}
C(\omega) &= \prod_p \bigg(1 - \frac{\omega(p)}{p}\bigg) \bigg(1 - \frac{1}{p}\bigg)^{-1} \\
&= 2 \prod_{p > 2} \bigg(1 - \frac{1}{(p-1)^2}\bigg) \prod_{\substack{p \mid N \\ p > 2}} \frac{p-1}{p-2} \\
&= 2C(N).
\end{aligned}
\end{equation}
Moreover, it is easy to check that for this $V(z)$ 
there exists an absolute constant $K > 1$ such that
\begin{equation}\label{Vz1/Vz2/G1}
\frac{V(z_1)}{V(z_2)} \le \frac{\log z_2}{\log z_1} \left( 1 + \frac{K}{\log z_1} \right)
\end{equation}
for $z_2 \ge z_1 \ge 2$. This verifies \eqref{RIcon3}. 

We must also control the error terms 
\begin{equation}\label{ErrTer/LS}
\sum_{1 \leq l \leq L }\sum_{q \mid P(z) } \lambda_l^{\pm} (q) r(\mathscr{A}, q)
\end{equation}
appearing on the right-hand side of Lemma \ref{Sieve-WF}, where $\lambda_q^{\pm}$ are well-factorable coefficients of order $1$ 
and of level $Q$. We note that by assumption \eqref{Lem24/con}, we now have 
$0 < \eta < 1/8, L=\exp(8\eta^{-3}),2 \leq z \leq Q^{1/2}$. 
Since $|\lambda_{l}^{\pm}(q)|\leq 1$ and $L=\exp(8\eta^{-3})$, 
Lemma \ref{BV-New} implies that the absolute value of the quantity in \eqref{ErrTer/LS} is 
\begin{equation}
\begin{aligned}
& \ll L\sum_{q\leq\sqrt{N}/\log^B N}\mu^{2}(q)\max_{y\le N}\max_{(a,q)=1}\bigg|\pi(y;q,a)-\frac{\Li(y)}{\varphi(q)}\bigg|\\
& \ll \frac{\exp(8\eta^{-3})N}{\log^{3} N},
\end{aligned}
\end{equation}
which is under control. 

We can now apply Lemma \ref{Sieve-WF} with
	$$
	z=N^{4/53},\quad Q=\frac{\sqrt{N}}{\log^{B} N}, \quad \eta=\varepsilon 
	$$
to deduce that  
\begin{equation}\label{G1}
\begin{aligned}
G_1\ge\ & \dfrac{2C(N)N}{\log N}\dfrac{e^{-\gamma}}{\log (N^{4/53})}f\bigg(\dfrac{53}{8}\bigg)\bigg(1+O\bigg(\varepsilon+\dfrac{e^K}{\varepsilon^8\log^{1/3} N}+\dfrac{\exp(8\varepsilon^{-3})}{\log  N}\bigg)\bigg)\\
\ge\ &g_1\dfrac{C(N)N}{\log^2 N}(1+O(e^K\varepsilon)),
\end{aligned}
\end{equation}
where 
\begin{equation}\label{g1}
g_1=\dfrac{53}{2}e^{-\gamma}f\bigg(\dfrac{53}{8}\bigg)\ge 14.87710. 
\end{equation}
In the above we have used the facts that 
	\begin{equation}\label{f6,e-gamma}
	f\bigg(\dfrac{53}{8}\bigg)\ge f(6)\ge 0.999895, \quad e^{-\gamma}\ge 0.561459,
	\end{equation}
and that the even integer $N$ satisfies
\begin{equation}\label{G1Error}
N>\exp(\exp(10\varepsilon^{-3}))>\varepsilon^{-\frac{9}{19\varepsilon}},\quad 0<\varepsilon<10^{-10}.
\end{equation}
Combining the above \eqref{G1} and \eqref{g1} gives the desired lower bound for $G_1$. 

\medskip 

By exactly the same idea, we can get that 
\begin{equation}\label{G2}
\begin{aligned}
G_2\ge\ & \dfrac{2C(N)N}{\log N}\dfrac{e^{-\gamma}}{\log (N^{4/33})}f\bigg(\dfrac{33}{8}\bigg)(1+O(e^K\varepsilon))\\
\ge\ & g_2\dfrac{C(N)N}{\log^2 N}(1+O(e^K\varepsilon)),
\end{aligned}
\end{equation}
where the constant $g_2$ can be estimated by Lemma \ref{Ff} as 
\begin{equation}\label{g2}
\begin{aligned}
g_2&=\dfrac{33}{2}e^{-\gamma}f\bigg(\dfrac{33}{8}\bigg)\\
&=8\bigg(\log \bigg(\dfrac{25}{8}\bigg)+\int_3^{25/8} \frac{d t}{t} \int_2^{t-1} \frac{\log (u-1)}{u} d u\bigg)\\
&\ge 9.11587.
\end{aligned}
\end{equation}
Also, combining \eqref{G2} and \eqref{g2} gives the desired lower bound for $G_2$. 

\subsection{Evaluation of $G_4,G_5,G_6,G_7$}
In this subsection we will prove upper bounds for $G_4$ and $G_5$, and lower bounds for $G_6$ and $G_7$, where
\[	
G_4 = \sum_{\substack{N^{4/53} \leq p \le N^{1/3} \\ (p, N) = 1}} S(\mathscr{A}_p; \mathscr{P}(N), N^{4/53}),\\
\]
\[
G_5 = \sum_{\substack{N^{4/53} \leq p \le N^{3/11} \\ (p, N) = 1}} S(\mathscr{A}_p; \mathscr{P}(N), N^{4/53}),
\]
\[
G_6 = \sum_{\substack{N^{4/53} \leq p_1 \le p_2 \le N^{4/33} \\ (p_1 p_2, N) = 1}} S(\mathscr{A}_{p_1 p_2}; \mathscr{P}(N),N^{4/53}),
\]
and 
\[ 
G_7 = \sum_{\substack{N^{4/53} \leq p_1 \le N^{4/33} \le p_2 \le N^{3/11} \\ (p_1 p_2, N) = 1}} S(\mathscr{A}_{p_1 p_2}; \mathscr{P}(N), N^{4/53}).
\]

We start with an upper bound for $G_4$. We want to apply the upper bound part 
in Lemma \ref{Sieve-WF} with
\[
z=N^{4/53},\quad Q=\frac{\sqrt{N}}{p \log^{B} N},\quad \eta=\varepsilon
\] 
to the summand $S(\mathscr{A}_{p};\mathscr{P}(N), N^{4/53})$ of $G_4$. The procedure is 
similar to that for $G_1$ in \S 5.1, except that here we aim at an upper bound instead of a lower bound. 
One checks that the conditions \eqref{def/SAPz}--\eqref{RIcon3} hold for the present sieving function 
$S(\mathscr{A}_{p};\mathscr{P}(N), N^{4/53})$, and this can be done similarly as in 
\eqref{New/A=}--\eqref{Vz1/Vz2/G1}. A simple calculation shows that, 
instead of \eqref{New/A=} and \eqref{G1X}, now we have 
\begin{equation}\label{New/A=/G4}
\mathscr{A}_p=\{N-q: q<(1-\varepsilon)N, \ q\equiv N \bmod{p}\} 
\end{equation}
and 
\begin{equation}\label{G1X/G4}
X=\frac{\Li(N-\varepsilon N)}{\varphi(p)}
= \frac{1}{\varphi(p)} \frac{N}{\log N} \bigg(1-\varepsilon+O\bigg(\dfrac{1}{\log N}\bigg)\bigg), 
\end{equation}
while \eqref{New/ome/p}--\eqref{Vz1/Vz2/G1} remain the same. Now the upper bound part 
in Lemma \ref{Sieve-WF} yields 
\begin{equation}
\begin{split} 
S(\mathscr{A}_{p};\mathscr{P}(N), N^{4/53})
\le\ &  
(1+O(e^K\varepsilon))\dfrac{N}{\log N} 
\frac{V(N^{4/53})}{\varphi(p)} F\bigg(\frac{\log (\sqrt{N} / p)}{\log (N^{4/53})}\bigg) \\ 
& + L \sum\limits_{\substack{q\le\sqrt{N} /p\log^{B} N\\ q\mid P(N^{4/53})}}|r(\mathscr{A}, pq)|. 
\end{split}  
\end{equation}
Inserting this into the expression of $G_4$, we see that the contribution from the last term above 
can be estimated by Lemma \ref{BV-New} as 
\begin{equation}
\begin{split} 
&\ll\sum\limits_{\substack{N^{4/53}\le p\le N^{1/3}\\(p, N)=1}}L
\sum\limits_{\substack{q\le\sqrt{N} /p \log^{B} N\\ q\mid P(N^{4/53})}}|r(\mathscr{A}, pq)|\\ 
&\ll L\sum\limits_{d\le\sqrt{N} /\log^{B} N}\mu^{2}(d)\max\limits _{y\le N}\max\limits _{(a,d)=1} 
\bigg|\pi(y ;d,a)-\frac{\Li(y)}{\varphi(d)}\bigg|\\
&\ll \frac{\exp(8\varepsilon^{-3})N}{\log^3 N},  
\end{split}  
\end{equation}
which is acceptable. It follows that 
\begin{equation}
G_4 \le (1+O(e^K\varepsilon)) \frac{N}{\log N} \sum_{\substack{N^{4/53} \leq p \le N^{1/3} \\ (p, N) = 1}} 
\frac{V(N^{4/53})}{\varphi(p)} F\bigg(\frac{\log (\sqrt{N} / p)}{\log (N^{4/53})}\bigg), 
\end{equation}  
where $V(z)$ is as \eqref{Vz}. 
The condition $(p,N)=1$ can be simply removed. 
Hence we have 
\begin{equation}
G_4 
\le (1+O(e^K\varepsilon))\dfrac{2C(N)N}{\log N}
\frac{e^{-\gamma}}{\log (N^{4/53})}\sum_{N^{4/53} \le p\le N^{1/3}} 
\frac{1}{p-1} F\bigg(\frac{\log (\sqrt{N} / p)}{\log (N^{4/53})}\bigg). 
\end{equation}  
The last sum over $p$ can be estimated by standard procedure of replacing sums over primes by integrals
as 
\begin{equation}\label{G4-pre}
\begin{aligned}
\le\ & \sum_{N^{4/53} \le p \le N^{1/3}} \frac{1}{p}F\bigg(\frac{1/2-\log p/\log N}{4/53}\bigg)\\
\le\ & \int_{N^{4/53}}^{N^{1/3}} \frac{1}{t}F\bigg(\frac{1/2-\log t/\log N}{4/53}\bigg)\dfrac{dt}{\log t}\\
\le\ & \int_{4/53}^{1/3} \frac{1}{u}F\bigg(\frac{53-106u}{8}\bigg)du. 
\end{aligned}
\end{equation}
Hence we conclude that 
\begin{equation}\label{G4}
G_4\le g_4\dfrac{C(N)N}{\log^2 N} (1+O(e^K\varepsilon)),
\end{equation}
where the constant $g_4$ can be calculated by Lemma \ref{Ff} as 
\begin{equation}\label{g4}
\begin{aligned}
g_{4}=\ & \dfrac{53}{2e^{\gamma}}\int_{4/53}^{1/3} \frac{1}{u}F\bigg(\frac{53-106u}{8}\bigg)du\\
=\ &\int_{4/53}^{1/3} \frac{8}{u(1-2u)}du+\int_{4/53}^{29/106} \frac{8du}{u(1-2u)}\int_2^{(45-106u)/8} \frac{\log (v-1)}{v} d v\\
&+\int_{4/53}^{13/106} \frac{8du}{u(1-2u)}\int_2^{(29-106u)/8} \frac{\log (v-1)}{v} d v \int_{v+2}^{(45-106u)/8} \frac{1}{w} \log \frac{w-1}{v+1} dw\\
\le\ & 23.60636. 
\end{aligned}
\end{equation}
Here we have also used \eqref{G1Error} as before.  
The above \eqref{G4} and \eqref{g4} give the desired upper bound for $G_4$. 

An upper bound for $G_5$ can be obtained in exactly the same way. We have 
\begin{equation}\label{G5}
\begin{aligned}
G_5\le\ & (1+O(e^K\varepsilon))\dfrac{2C(N)N}{\log N}
\dfrac{e^{-\gamma}}{\log (N^{4/53})}\sum_{N^{4/53} \le p\le N^{3/11}} \frac{1}{p-1} 
F\bigg(\frac{\log (\sqrt{N} / p)}{\log (N^{4/53})}\bigg)\\
\le\ & g_5\dfrac{C(N)N}{\log^2 N}(1+O(e^K\varepsilon)),
\end{aligned}
\end{equation}
where $g_5>0$ is a constant that can be estimated by Lemma \ref{Ff} as 
\begin{equation}\label{g5}
\begin{aligned}
g_{5}=\ & \dfrac{53}{2e^{\gamma}}\int_{4/53}^{3/11} \frac{1}{u}F\bigg(\frac{53-106u}{8}\bigg)du\\
=\ &\int_{4/53}^{3/11} \frac{8}{u(1-2u)}du+\int_{4/53}^{3/11} \frac{8du}{u(1-2u)}\int_2^{(45-106u)/8} \frac{\log (v-1)}{v} d v\\
&+\int_{4/53}^{13/106} \frac{8du}{u(1-2u)}\int_2^{(29-106u)/8} \frac{\log (v-1)}{v} d v \int_{v+2}^{(45-106u)/8} \frac{1}{w} \log \frac{w-1}{v+1} dw\\
\le\ & 19.51976.
\end{aligned}
\end{equation}
The above two formulas \eqref{G5} and \eqref{g5} give the desired upper bound for $G_5$. 

In the following we work on the lower bounds for $G_6$ and $G_7$. 
Note that the condition $(p_{1}p_{2},N)=1$ can be removed by the estimate 
\begin{equation}
\sum_{\substack{p_{1}\mid N\\ p_{1}\ge N^{\kappa}}}\sum_{N^{\kappa}\le p_{2}<N^{\sigma}}\frac{N}{p_{1}p_{2}}+\sum_{N^{\kappa}\le p_{1}<N^{\sigma}}\sum_{\substack{p_{2}\mid N\\ p_{2}\ge N^{\kappa}}}\frac{N}{p_{1}p_{2}}\ll N^{1-\kappa}\log N,
\end{equation}
where $1/30<\kappa<\sigma<1/2$. As before we can deduce, from Lemmas \ref{Sieve-WF} and \ref{BV-New}, that
\begin{equation*}
\begin{aligned}
G_6
\ge &\   
(1+O(e^K\varepsilon))\dfrac{2C(N)N}{\log N}\dfrac{e^{-\gamma}}{\log (N^{4/53})}\\ 
&\times \sum_{N^{4/53} \leq p_1 \le p_2 \le N^{4/33}} \frac{1}{p_1p_2}f\bigg(\frac{1/2-\log (p_1p_2)/\log N}{4/53}\bigg), 
\end{aligned}
\end{equation*}
while the above double sums over $p_1$ and $p_2$ can be estimated as 
\begin{equation*}
\begin{aligned}
\ge & \int_{N^{4/53}}^{N^{4/33}} \dfrac{dt_1}{t_1\log t_1}\int_{t_1}^{N^{4/33}}f\bigg(\frac{1/2-\log (t_1t_2)/\log N}{4/53}\bigg)\dfrac{dt_2}{t_2\log t_2}\\
\ge & \int_{4/53}^{4/33}\dfrac{du}{u}\int_{u}^{4/33} f\bigg(\frac{53-106u-106v}{8}\bigg)\dfrac{dv}{v}. 
\end{aligned}
\end{equation*}
Hence we obtain 
\begin{equation}\label{G6}
G_6 \ge g_6\dfrac{C(N)N}{\log^2 N} (1+O(e^K\varepsilon)),
\end{equation}
where the constant $g_6$ can be estimated by Lemma \ref{Ff} as 
\begin{equation}\label{g6}
\begin{aligned}
g_{6}=\ & \dfrac{53}{2e^{\gamma}}\int_{4/53}^{4/33}\dfrac{du}{u}\int_{u}^{4/33}f\bigg(\frac{53-106u-106v}{8}\bigg)\dfrac{dv}{v}\\
=\ &\int_{4/53}^{4/33}\dfrac{du}{u}\int_{u}^{4/33} \dfrac{8\log (45-106u-106v)-8\log 8}{v(1-2u-2v)}dv\\
&+\int_{4/53}^{21/212}\dfrac{dv}{v}\int_{4/53}^{v} \dfrac{8du}{u(1-2u-2v)}\int_3^{(45-106u-106v)/8} \frac{dw}{w} \int_{2}^{w-1}\dfrac{\log (x-1)}{x}dx\\
&+\int_{21/212}^{4/33}\dfrac{dv}{v}\int_{4/53}^{21/106-v} \dfrac{8du}{u(1-2u-2v)}\int_3^{(45-106u-106v)/8} \frac{dw}{w} \int_{2}^{w-1}\dfrac{\log (x-1)}{x}dx\\
\ge\ & 1.63357. 
\end{aligned}
\end{equation}
The above \eqref{G6} and \eqref{g6} give the desired lower bound for $G_6$. 

A lower bound for $G_7$ can be obtained in the same way. 
Thus, we have 
\begin{equation}\label{G7}
\begin{aligned}
G_7\ge&\ (1+O(e^K\varepsilon))\dfrac{53}{2e^{\gamma}}\dfrac{C(N)N}{\log ^2N}\int_{4/53}^{4/33}\dfrac{du}{u}\int_{4/33}^{3/11} \dfrac{1}{v}f\bigg(\frac{53-106u-106v}{8}\bigg)\dfrac{dv}{v}\\
\ge &\ g_7\dfrac{C(N)N}{\log^2 N} (1+O(e^K\varepsilon)),
\end{aligned}
\end{equation}
where the constant $g_7$ can be estimated by Lemma \ref{Ff} as 
\begin{equation}\label{g7}
	\begin{aligned}
	g_{7}=\ & \dfrac{53}{2e^{\gamma}}\int_{4/53}^{4/33}\dfrac{du}{u}\int_{4/33}^{3/11}f\bigg(\frac{53-106u-106v}{8}\bigg)\dfrac{dv}{v}\\
	=\ &\int_{4/53}^{4/33}\dfrac{du}{u}\int_{4/33}^{37/106-u} \dfrac{8\log (45-106u-106v)-8\log 8}{v(1-2u-2v)}dv\\
	&+\int_{4/33}^{13/106}\dfrac{dv}{v}\int_{4/53}^{21/106-v} \dfrac{8du}{u(1-2u-2v)}\int_3^{(45-106u-106v)/8} \frac{dw}{w} \int_{2}^{w-1}\dfrac{\log (x-1)}{x}dx\\
\ge\ & 3.79029. 
\end{aligned}
\end{equation}
Note that we have applied the property that $f(s)=0$ for $s\le 2$.  
The desired lower bound for $G_7$ now follows from \eqref{G7} and \eqref{g7}.  

\subsection{Evaluation of $G_3,G_8,G_9,G_{10}$}
In this subsection we will estimate upper bounds for $G_3, G_8, G_9, G_{10}$, where
\[
G_3 =\sum\limits_{\substack{N^{9/19-\varepsilon}\le p\le N^{1/2} \\(p, N)=1}}S(\mathscr{A}_p ; \mathscr{P}(N), p),\\
\]
\[
G_8 = \sum_{\substack{N^{3/11} \leq p_1 \le p_2 \le (N/p_1)^{1/2} \\ (p_1 p_2, N) = 1}} S(\mathscr{A}_{p_1 p_2}; \mathscr{P}(N p_1), p_2),
\] 
\[
G_9 = \sum_{\substack{N^{4/53} \leq p_1 \le N^{1/3} \leq p_2 \le (N/p_1)^{1/2} \\ (p_1 p_2, N) = 1}} S(\mathscr{A}_{p_1 p_2}; \mathscr{P}(N p_1), p_2),
\]
and 
\[
G_{10} = \sum_{\substack{N^{4/33} \leq p_1 \le N^{3/11} \leq p_2 \le (N/p_1)^{1/2} \\ (p_1 p_2, N) = 1}} S(\mathscr{A}_{p_1 p_2}; \mathscr{P}(N p_1), (N/p_1 p_2)^{1/2}). 
\]
	
We shall only majorize  $G_8$, and the others can be treated in the same way. 

First, we verify that \(G_8\) equals the number of quadruples \((p_1,p_2,p_3,p)\) such that  
\(p\) is a prime \(\le N\),  
\(N-p = p_1p_2p_3\), and  
\begin{equation}\label{G8-quad}
N^{3/11} \le p_1 \le p_2 \le (N/p_1)^{1/2},\quad p_3 \ge p_2,\quad p_1p_2p_3\ge \varepsilon N,\quad (p_1p_2p_3, N) = 1.
\end{equation}
Indeed, for fixed \(p_1,p_2\) with \((p_1p_2,N)=1\), the sieve term  
\(S(\mathscr{A}_{p_1p_2}; \mathscr{P}(N p_1), p_2)\) counts those \(a\in\mathscr{A}_{p_1p_2}\) that have no prime factor \(<p_2\) from \(\mathscr{P}(N p_1)\).  
Since \(\mathscr{A}=\{N-p:p<(1-\varepsilon)N\}\), each \(a=N-p\) for some prime \(p\) 
and satisfies \(a \equiv 0 \bmod {p_1p_2}\). Write \(a=p_1p_2m\); then \(N-p=p_1p_2m\).  
The coprimality condition means \(m\) has no prime divisor \(q<p_2\) with \((q,Np_1)=1\).  

Because $p_2\ge p_1\ge N^{3/11}>N^{1/4}$,
any composite $m$ would possess a prime factor
$\le\sqrt{m}\le\sqrt{N/(p_1p_2)}<p_2$,
and this factor would automatically belong to $\mathscr{P}(N p_1)$
unless it divides $N$.
Hence $m$ must be either $1$ or a prime $\ge p_2$.

The case $m=1$ gives $N-p=p_1p_2\ge N^{6/11}$,
which for large $N$ contradicts $p<\varepsilon N$;
thus $m$ is a prime, say $p_3$, with $p_3\ge p_2$.
The inequality $p_1p_2p_3=N-p\ge\varepsilon N$ follows from $p<(1-\varepsilon)N$,
and $(p_1p_2p_3,N)=1$ is inherited from $(p_1p_2,N)=1$ together with
$p_3\in\mathscr{P}(N p_1)$.
Therefore each summand of $G_8$ is exactly the number of admissible $p_3$,
and $G_8$ itself equals the number of quadruples satisfying \eqref{G8-quad}.

To obtain an upper bound suitable for a linear sieve, we define the set 
\begin{equation}
\mathscr{M}:= \{ m : m = p_1p_2, N^{3/11} \le p_1 \le p_2 \le (N/p_1)^{1/2}, (p_1p_2, N) = 1 \},
\end{equation}
and collect the corresponding values $b=N-mp$ where $p$ is an arbitrary prime
$\le N/m$ and the result is still $\le(1-\varepsilon)N$:
\begin{equation}
\mathscr{B}:= \{ b\le (1-\varepsilon) N : b = N - mp, m \in \mathscr{M}, p \text{ prime, } p \le N/m \}.
\end{equation}
Every quadruple counted by $G_8$ yields $b=N-p_1p_2p_3$ which is a prime
$\le (1-\varepsilon)N$ and therefore belongs to $\mathscr{B}$.

Now we examine the sifting function $S(\mathscr{B}; \mathscr{P}(N), N^{1/2})$.
It counts those $b\in\mathscr{B}$ having no prime divisor $< N^{1/2}$ from
$\mathscr{P}(N)$.
All $b\in\mathscr{B}$ are prime numbers, and apart from possibly $O(N^{1/2})$
exceptions they satisfy $b> N^{1/2}$. 
Finally we have
$$
G_8 \le S(\mathscr{B}; \mathscr{P}(N), N^{1/2}) + O(N^{1/2}). 
$$

By applying Lemma \ref{Sieve-WF} with
$$
X = \sum_{m \in \mathscr{M}} \bigg(\Li\bigg( \frac{N}{m} \bigg)-\Li\bigg( \frac{\varepsilon N}{m} \bigg)\bigg), 
\quad Q = \frac{\sqrt{N}}{\log^{B} N},
$$
we obtain
\begin{equation}\label{G8old}
G_8 \leq \frac{8 C(N) X}{\log N} ( 1 + O(e^K\varepsilon) ) + O(\sqrt{N} + R_1 + R_2),
\end{equation}
where
$$ 
R_1 := \sum_{\substack{q \leqslant \sqrt{N}/\log^{B} N \\ (q,N)=1}} \mu^2(q) \biggl| \sum_{\substack{m \in \mathscr{M} \\ (q,m)=1}} \biggl( \sum_{\substack{mp \leqslant N \\ mp \equiv N \bmod q}} 1 - \frac{\Li(N/m)}{\varphi(q)} \biggr) \biggr|, 
$$
and 
$$ 
R_2 := \sum_{\substack{q \leqslant \sqrt{N}/\log^{B} N \\ (q,N)=1}} \frac{\mu^2(q)}{\varphi(q)} 
\sum_{\substack{m \in \mathscr{M} \\ (q,m)>1}}
\Li\biggl(\frac{N}{m}\biggr). 
$$
Let $f(m)$ be the characteristic function of $\mathscr{M}$. Since $m \leqslant N^{2/3}$ for $m \in \mathscr{M}$, Lemma \ref{BV-New} implies
\begin{equation}
R_1 = \sum_{\substack{q \leqslant \sqrt{N}/\log^{B} N \\ (q,N)=1}} \mu^2(q) \biggl| \sum_{\substack{m \leqslant N^{2/3} \\ (q,m)=1}} f(m) \biggl( \sum_{\substack{mp \leqslant N \\ mp \equiv N \bmod q}} 1 - \frac{\Li(N/m)}{\varphi(q)} \biggr) \biggr| \ll \frac{N}{\log^{3} N}.
\end{equation}
Noticing that $(q,m) > 1$ implies $(q,m) \geqslant N^{3/11}$ for 
$m \in \mathscr{M}$ because all the prime factor of $m$ is bigger than $N^{3/11}$, we have
\begin{equation}\label{Chen-R2}
\begin{split}
R_2 &\ll \frac{N}{\log N} \sum_{q \leqslant \sqrt{N}} \frac{\mu^2(q)}{\varphi(q)} 
\sum_{\substack{m \leqslant N^{2/3} \\ (q,m) \geqslant N^{3/11}}} \frac{1}{m} 
\ll N \sum_{q \leqslant \sqrt{N}} \frac{\mu^2(q)}{\varphi(q)} \sum_{d|q, d \geqslant N^{3/11}} \frac{1}{d} \\
&\ll N \sum_{N^{3/11} \le d \leqslant N^{1/2}} \frac{\mu^2(d)}{d\varphi(d)} 
\sum_{l \leqslant \sqrt{N}/d} \frac{\mu^2(l)}{\varphi(l)} \ll N \log N \sum_{N^{3/11} \le d \leqslant N^{1/2}} \frac{\mu^2(d)}{d\varphi(d)}\\
&\ll N^{8/11} \log^2 N.
\end{split}
\end{equation}

Now we give an asymptotic expression for  $X$. By the Prime Number Theorem,
	\begin{equation}\label{G8x}
	\begin{aligned}
	X &= \bigg(1-\varepsilon+O\bigg(\dfrac{1}{\log N}\bigg)\bigg) \sum_{N^{3/11} \leqslant p_1 \le p_2 \le (N/p_1)^{1/2}} \frac{N}{p_1p_2 \log(N/p_1p_2)} \\
	&= \bigg(1-\varepsilon+O\bigg(\dfrac{1}{\log N}\bigg)\bigg) \frac{N}{\log N} \int_{3/11}^{1/3} \frac{\log(1/t - 2)}{t(1 - t)}dt \\
	&= \bigg(1-\varepsilon+O\bigg(\dfrac{1}{\log N}\bigg)\bigg) \frac{N}{\log N} \int_{2}^{8/3} \frac{\log(t - 1)}{t} dt. 
	\end{aligned}
	\end{equation}
Inserting \eqref{G8x} into \eqref{G8old}, we have
\begin{equation}\label{G8}
G_8 \le g_8\dfrac{C(N)N}{\log^2 N}(1+O(e^K\varepsilon)),
\end{equation}
where 
\begin{equation}\label{g8}
g_{8}=8\int_{2}^{8/3} \frac{\log(t - 1)}{t} dt
\le 0.60962.
\end{equation}

Similarly we have
\begin{equation}\label{G3}
G_3 \le g_3 \frac{C(N)N}{\log^2 N} (1+O(e^K\varepsilon)), 
\end{equation}
where 
\begin{equation}\label{g3}
g_{3}
= 8\int_{9/19}^{1/2} \frac{du}{u(1-u)} 
=8\log(10/9)\le 0.84289. 
\end{equation}

For $G_9$, when $p_1 < N^{1/10}$ we use Lemma \ref{Fouvry-new} in place of Lemma \ref{BV-New}.  
The range $p_1 < N^{1/10}$ allows Lemma \ref{Fouvry-new} to supply a much sharper error term than the general estimate of Lemma \ref{BV-New}. When $p_1 \ge N^{1/10}$ we continue to apply Lemma \ref{BV-New}.

For $G_{10}$, note that the sieve level is $(N/p_1p_2)^{1/2}$, which reverses the roles compared with $G_8$.  
Recall that in $G_8$ the condition $p_1 \ge N^{3/11} > N^{1/4}$ forced the remaining factor $m$ to be a single prime: if $m$ were composite, it would contain a prime factor $\le \sqrt{m} \le \sqrt{N/(p_1p_2)} < p_2$, and that factor would be removed by the sieve.  
In $G_{10}$ the sieve limit $(N/p_1p_2)^{1/2}$ plays exactly the same part.  
Suppose the unsieved part of $N-p_1p_2$ were a product of two primes $p_3p_4$; then one of them would satisfy $p_3 \le (N/p_1p_2)^{1/2}$ and would therefore be sieved out.  
Hence only a single prime factor $p_3$ can survive, and we are once again in the situation where the sifted set consists of numbers having exactly one large prime factor, just as in the analysis of $G_8$. 

It follows from the above analysis that for $j=9,10$ we have
\begin{equation}\label{G9-10}
G_j \le g_j \frac{C(N)N}{\log^2 N} (1+O(e^K\varepsilon)), 
\end{equation}
where 
\begin{equation}\label{g9}
\begin{aligned}
g_{9}&=\dfrac{36}{5}\iint\limits_{\substack{4/53\le u\le 1/10\\ 1/3\le v\le (1-u)/2}}\dfrac{dudv}{uv(1-u-v)(1-u)}+8\iint\limits_{\substack{1/10\le u\le 1/3\\ 1/3\le v\le (1-u)/2}}\dfrac{dudv}{uv(1-u-v)}\\
&\le 5.27231,
\end{aligned}
\end{equation}
and 
\begin{equation}\label{g10}
g_{10}=8\iint\limits_{4/33\le u\le 3/11\le v\le (1-u)/2}\dfrac{dudv}{uv(1-u-v)}
\le 5.40996. 
\end{equation}
This completes the analysis for $G_3,G_8,G_9,G_{10}$. 

\subsection{Evaluation of $G_{11},G_{12}$}
The aim of this subsection is upper bounds for $G_{11}$ and $G_{12}$, where
\[
G_{11} = \sum_{\substack{N^{4/53} \leq p_1 \le p_2 \le p_3 \le p_4 \le N^{4/33} \\ (p_1 p_2 p_3 p_4, N) = 1}} S(\mathscr{A}_{p_1 p_2 p_3 p_4}; \mathscr{P}(Np_1), p_2),
\] 
and 
\[
G_{12} = \sum_{\substack{N^{4/53} \leq p_1 \le p_2 \le p_3 \le N^{4/33} \leq p_4 \le N^{3/11}\\ (p_1 p_2 p_3 p_4, N) = 1}} S(\mathscr{A}_{p_1 p_2 p_3 p_4}; \mathscr{P}(Np_1), p_2).
\] 
We use Buchstab's function to estimate $X$ instead of applying the Prime Number Theorem directly. 
Similarly as before, we have, for $j=11,12$, 
\begin{equation}\label{G11-12}
G_j \le  g_j \dfrac{C(N)N}{\log^2 N} (1+O(e^K\varepsilon)), 
\end{equation}
where
\begin{equation*}
\begin{aligned}
g_{11}
=\ & \dfrac{36}{5}\int_{4/53}^{1/10}\dfrac{dt_1}{t_1(1-t_1)}\int_{t_1}^{4/33}\dfrac{dt_2}{t_2^2}\int_{t_2}^{4/33}\dfrac{dt_3}{t_3}\int_{t_3}^{4/33}w\bigg(\dfrac{1-t_1-t_2-t_3-t_4}{t_2}\bigg)\dfrac{dt_4}{t_4}\\
&+8\int_{1/10}^{4/33}\dfrac{dt_1}{t_1}\int_{t_1}^{4/33}\dfrac{dt_2}{t_2^2}\int_{t_2}^{4/33}\dfrac{dt_3}{t_3}\int_{t_3}^{4/33}w\bigg(\dfrac{1-t_1-t_2-t_3-t_4}{t_2}\bigg)\dfrac{dt_4}{t_4}. 
\end{aligned}
\end{equation*}
Since $w(u)\le 0.561522$ for $u \geq 3.5$, we compute that 
\begin{equation}\label{g11}
\begin{aligned}
g_{11}
\le\ & \dfrac{36}{5}\times 0.561522\times \int_{4/53}^{1/10}\dfrac{dt_1}{t_1(1-t_1)}\int_{t_1}^{4/33}\dfrac{dt_2}{t_2^2}\int_{t_2}^{4/33}\dfrac{dt_3}{t_3}\int_{t_3}^{4/33}\dfrac{dt_4}{t_4}\\
&+8\times 0.561522\times\int_{1/10}^{4/33}\dfrac{dt_1}{t_1}\int_{t_1}^{4/33}\dfrac{dt_2}{t_2^2}\int_{t_2}^{4/33}\dfrac{dt_3}{t_3}\int_{t_3}^{4/33}\dfrac{dt_4}{t_4}\\
\le\ & 0.10191. 
\end{aligned}
\end{equation}
This is the desired upper bound for $g_{11}$. 

For $j=12$, we have 
\begin{equation*}
\begin{aligned}
g_{12}
=\ &\dfrac{36}{5}\int_{4/53}^{1/10}\dfrac{dt_1}{t_1(1-t_1)}\int_{t_1}^{4/33}\dfrac{dt_2}{t_2^2}\int_{t_2}^{4/33}\dfrac{dt_3}{t_3}\int_{4/33}^{3/11}w\bigg(\dfrac{1-t_1-t_2-t_3-t_4}{t_2}\bigg)\dfrac{dt_4}{t_4}\\
&+8\int_{1/10}^{4/33}\dfrac{dt_1}{t_1}\int_{t_1}^{4/33}\dfrac{dt_2}{t_2^2}\int_{t_2}^{4/33}\dfrac{dt_3}{t_3}\int_{4/33}^{3/11}w\bigg(\dfrac{1-t_1-t_2-t_3-t_4}{t_2}\bigg)\dfrac{dt_4}{t_4}. 
\end{aligned}
\end{equation*}

Now we claim that 
\begin{equation}\label{w=/w=}
\begin{aligned}
\begin{cases} 
w(u)\le 0.561990, \quad u \geq 3.16, \\ 
w(u)\le 0.564383, \quad u \geq 3, 
\end{cases} 
\end{aligned}
\end{equation}
from which we deduce that 
\begin{equation}\label{g12}
\begin{aligned}
g_{12} 
\le\ & \dfrac{36}{5}\times 0.561990\times \int_{4/53}^{1/10}\dfrac{dt_1}{t_1(1-t_1)}\int_{t_1}^{4/33}\dfrac{dt_2}{t_2^2}\int_{t_2}^{4/33}\dfrac{dt_3}{t_3}\int_{4/33}^{3/11}\dfrac{dt_4}{t_4}\\
&+8\times 0.564383\times\int_{1/10}^{4/33}\dfrac{dt_1}{t_1}\int_{t_1}^{4/33}\dfrac{dt_2}{t_2^2}\int_{t_2}^{4/33}\dfrac{dt_3}{t_3}\int_{4/33}^{3/11}\dfrac{dt_4}{t_4}\\
\le\ & 0.66821. 
\end{aligned}
\end{equation}
This is the desired upper bound for $g_{12}$. 

We must prove our claim \eqref{w=/w=}. Recalling that 
$$
w(u) =
\frac{1+\log(u-1)}{u}+\dfrac{1}{u}\int_{2}^{u-1}\dfrac{\log(t-1)}{t}dt
$$
for $3 \le u \leq 4$, and $w(u)\le 0.561522$ for $u \geq 3.5$,  we only need to treat $u\in (3, 3.5)$. 
For $u\in (3, 3.5)$, the derivative of $w(u)$ is 
$$
w'(u)=\frac{1}{u^2(u-1)}\bigg(1-(u-1)\log (u-1)+u\log (u-2)-(u-1)\int_{2}^{u-1}\frac{\log (t-1)}{t}dt\bigg). 
$$
Writing   
$$
w_2(u)=\dfrac{1}{u-1}-\log (u-1)+\dfrac{u}{u-1}\log (u-2)-\int_{2}^{u-1}\frac{\log (t-1)}{t}dt,
$$
we have $u^2 w'(u)=w_2(u)$, and 
\[
w_2'(u)=\dfrac{u}{(u-1)^2(u-2)}(1-(u-2)\log(u-2)). 
\] 
Since $x\log x$ increase when $x\ge 1$, we have $w_2'(u)>0$ when $3\leq u<2+x_0$, 
where $x_0 \approx 1.76322$ is the solution of $x\log x=1$. 
One computes that 
\[ 
w_2(3.4)=\dfrac{5}{12}-\log (2.4)+\dfrac{17}{12}\log (1.4)-\int_{2}^{2.4}\frac{\log (t-1)}{t}dt\approx -0.01359<0,
\] 
while 
\[
w_2(3.5)=\dfrac{2}{5}-\log (2.5)+\dfrac{7}{5}\log (1.5)-\int_{2}^{2.5}\frac{\log (t-1)}{t}dt\approx 0.00475>0. 
\]	
It follows that there exists an $u_0\in (3.4,3.5)$ such that $w(u)$ decreases for $u\le u_0$, and increases for $u\ge u_0$. 
Direct calculation shows that 
\[
w(3)=\dfrac{1+\log 2}{3}\le 0.564383, 
\] 
and 
\[
w(3.16)=\dfrac{1+\log (2.16)}{3.16}+\dfrac{1}{3.16}\int_{2}^{2.16}\dfrac{\log(t-1)}{t}dt\le 0.561990. 
\] 
This proves the claim \eqref{w=/w=}. 
	
\subsection{Completion of the proof of Theorem \ref{Goldbach1}} 
By \eqref{G1}, \eqref{G2}, \eqref{G4}, \eqref{G5}, \eqref{G6}, \eqref{G7}, \eqref{G8}, \eqref{G3}, \eqref{G9-10}, \eqref{G11-12}, 
we have
\begin{equation}
\begin{split}
4D_{1,1.9}(N)
\ge\ & \dfrac{C(N)N}{\log^2(N)}(3\times 14.87710+9.11587-4\times 0.84289 \\ 
& -23.60636-19.51976+1.63357+3.79029-2\times 0.60962\\ 
& -5.27231-5.40996-0.10191-0.66821)\\
\ge\ & 0.00172\dfrac{C(N)N}{\log^2(N)}. 
\end{split}
\end{equation}
This proves Theorem \ref{Goldbach1}.  

\section{Proof of Theorem~\ref{Twin1}}\label{sec6}
In Proposition \ref{Twinweight} we specify $a=1.75$ and 
\begin{equation}\label{Cond/Twin}
\alpha=1/12,\quad \beta=1/7,\quad \gamma=2/7,\quad \tau=3/7-\varepsilon,\quad x>\exp(\exp(10\varepsilon^{-3})).
\end{equation} 
Thus we deduce from \eqref{Twinweight4} that 
\begin{equation}\label{Twinw4}
\begin{split}
4\pi_{1,1.75}(x) \geq\ &3T_1 + T_2 -4T_3- T_4- T_5 + T_6 + T_7 \\
&- 2T_8 - T_9 - T_{10} - T_{11} - T_{12} + O(x^{11/12}),
\end{split}
\end{equation}
where $T_1, \ldots, T_{12}$ are as in Proposition \ref{Twinweight}. We are going to estimate $T_1, \ldots, T_{12}$ 
in the subsections below.

\subsection{Evaluation of $T_1,T_2$}
In this subsection we will obtain lower bounds for $T_1$ and $T_2$ respectively, where
\[
T_1 = S(\mathscr{A}; \mathscr{P}, x^{1/12}), 
\quad 
T_2 = S(\mathscr{A}; \mathscr{P}, x^{1/7}).
\]

We want to apply Lemma \ref{Sieve-WF} to estimate $T_1$. Recalling \eqref{Def/A=Twin},
we have at present that
\begin{equation}\label{New/A=Twin3}
\mathscr{A}=\{p+2:\varepsilon x<p< x\},
\end{equation}
and, before applying Lemma \ref{Sieve-WF}, we must check that the conditions
\eqref{def/SAPz}--\eqref{RIcon3} hold for the present $\mathscr{A}$. We compute that
\begin{equation}\label{T1X}
X=\Li(x)-\Li(\varepsilon x)=\dfrac{x}{\log x}\bigg(1-\varepsilon+O\bigg(\dfrac{1}{\log x}\bigg)\bigg).
\end{equation}
We also have to compute the functions $\omega(p)$ and $V(z)$ defined as in \eqref{RIcon1} and \eqref{def/Vz}, respectively,
associated to this sequence $\mathscr{A}$, and
this has been done in classical literatures concerning the Twin Prime Conjecture,
see for example \cite[\S1.3, Example 5]{HR74}, where it is proved that
\begin{equation}\label{New/ome/p-Twin}
\omega(p)= \frac{p}{p-1}
\end{equation}
for odd prime $p$, and $\omega(2)=0$. Thus \eqref{RIcon2} is satisfied.
It is also proved there that
\begin{equation}\label{Vz/Twin}
V(z) = C(\omega) \frac{e^{-\gamma}}{\log z} \bigg(1 + O\bigg(\frac{1}{\log z}\bigg)\bigg),
\end{equation}
with
\begin{equation}\label{Come/Twin}
\begin{aligned}
C(\omega) &= \prod_p \bigg(1 - \frac{\omega(p)}{p}\bigg) \bigg(1 - \frac{1}{p}\bigg)^{-1} \\
&= 2 \prod_{p > 2} \bigg(1 - \frac{1}{(p-1)^2}\bigg) = C,  
\end{aligned}
\end{equation}
where $C$ is the constant in the Twin Prime Conjecture as in \eqref{C2}. Moreover, it is easy to check that for this $V(z)$
there exists an absolute constant $K > 1$ such that
\begin{equation}\label{Vz1/Vz2/T1}
\frac{V(z_1)}{V(z_2)} \le \frac{\log z_2}{\log z_1} \left( 1 + \frac{K}{\log z_1} \right)
\end{equation}
for $z_2 \ge z_1 \ge 2$. This verifies \eqref{RIcon3}.

Now we use Lemma \ref{BFI} instead of Lemma \ref{BV-New}, and apply Lemma \ref{Sieve-WF} with
\[
z=x^{1/12},\quad Q=x^{4/7-\varepsilon}, \quad \eta=\varepsilon,
\]
to deduce that
\begin{equation}\label{T1}
\begin{aligned}
T_1\ge\ & \dfrac{Cx}{\log x}\dfrac{e^{-\gamma}}{\log (x^{1/12})}f\bigg(\dfrac{48}{7}\bigg) (1+O(e^K\varepsilon))\\
\ge\ & t_1 \dfrac{Cx}{\log^2 x}(1+O(e^K\varepsilon)),
\end{aligned}
\end{equation}
where
\begin{equation}\label{t1}
t_1=12e^{-\gamma}f\bigg(\dfrac{48}{7}\bigg)\ge 6.73680.
\end{equation}
In the above we have used the facts that
\begin{equation}\label{f8,E-gamma}
f\bigg(\dfrac{48}{7}\bigg)\ge f(6)\ge 0.999895, \quad e^{-\gamma}\ge 0.561459,
\end{equation}
and that the real number $x$ satisfies
\begin{equation}\label{T1Error}
x>\exp(\exp(10\varepsilon^{-3}))>\varepsilon^{-\frac{9}{19\varepsilon}},\quad 0<\varepsilon<10^{-10}.
\end{equation}
The above \eqref{T1} and \eqref{t1} give the desired lower bound for $T_1$.

By exactly the same idea, we can get that
\begin{equation}\label{T2}
\begin{aligned}
T_2\ge\ & \dfrac{Cx}{\log x}\dfrac{e^{-\gamma}}{\log (x^{1/7})}f(4) (1+O(e^K\varepsilon))\\
\ge\ & t_2\dfrac{Cx}{\log^2 x}(1+O(e^K\varepsilon)),
\end{aligned}
\end{equation}
where the constant $t_2$ can be computed by Lemma \ref{Ff} as 
\begin{equation}\label{t2}
t_2=7e^{-\gamma}f(4)
=\frac{7}{2}\log 3 
\ge 3.84514.
\end{equation}
Interting \eqref{t2} into \eqref{T2}, we get the desired lower bound for $T_2$.
	
\subsection{Evaluation of $T_4,T_5,T_6,T_7$}
In this subsection we will estimate upper bounds of $T_4$ and $T_5$, and lower bounds of $T_6$ and $T_7$, 
where
\[
T_4 = \sum_{x^{1/12} \leq p \le x^{1/3}} S(\mathscr{A}_p; \mathscr{P}, x^{1/12}), 
\]
\[
T_5 = \sum_{x^{1/12} \leq p \le x^{2/7}} S(\mathscr{A}_p; \mathscr{P}, x^{1/12}),
\]
\[ 
T_6 = \sum_{x^{1/12} \leq p_1 \le p_2 \le x^{2/7}} S(\mathscr{A}_{p_1 p_2}; \mathscr{P}, x^{1/12}),
\] 
and 
\[
T_7 = \sum_{x^{1/12} \leq p_1 \le x^{1/7} \le p_2 \le x^{2/7}} S(\mathscr{A}_{p_1 p_2}; \mathscr{P}, x^{1/12}).
\] 

First we study \( T_4 \).  
We partition the interval \([x^{1/12}, x^{1/3}]\) into \(O(\log x)\) subintervals of the form \([P, 2P)\), 
and classify them into four types based on which of the following four intervals $P$ belongs to:  
\[
[x^{1/12}, x^{2/7-\varepsilon}],\quad [x^{2/7-\varepsilon}, x^{29/100}],\quad [x^{29/100}, 
x^{1/3-\varepsilon}], \quad [x^{1/3-\varepsilon}, x^{1/3}].
\]
To control the error terms in \( S(\mathscr{A}_p; \mathscr{P}, x^{1/12}) \) for \( p \in [P, 2P) \), 
we apply conditions (C.1), (C.2), and (C.3) of Lemma \ref{Fouvry-level} to the first three types,  
respectively. We have
\begin{equation}
\begin{aligned}
S(\mathscr{A}_{p};\mathscr{P},x^{1/12}) 
\leqslant &\ \frac{(1+O(\varepsilon))x}{\log x}\frac{V(x^{1/12})}{\varphi(p)}F\bigg(\frac{\log(Q/P)}{\log (x^{1/12})}\bigg) \\ 
&\ + \sum_{l<L}\sum_{q|P(x^{1/12})}\lambda_{l}^{+}(q)r(\mathscr{A},pq), 
\end{aligned}
\end{equation}
where \( \lambda_{l}^{+}(q) \) is well-factorable of level \( Q/P \) and order $1$.

For $[P, 2P)$ of the first type, we apply condition (C.1). Let \( \pi_P \) denote the characteristic 
function of primes in \( [P, 2P) \). Since \( P \leqslant x^{2/7-\varepsilon} \) implies \( P \leqslant Q/P \), 
Lemma \ref{WuWF} shows that \( \pi_P * \lambda_{l}^{+} \) is well-factorable of level \( Q \) and order 2. 
Thus, Lemma \ref{Fouvry-level} gives
\begin{equation}
\sum_{P\leqslant p<2P}\sum_{l<L}\sum_{q|P(x^{1/12})} 
\lambda_{l}^{+}(q)r(\mathscr{A},pq) \ll_{\varepsilon} \frac{x}{\log^{4} x}. 
\end{equation}
Similar arguments apply to the second and third types 
using conditions (C.2) and (C.3), respectively. For the last type, we use the trivial estimate
\begin{equation}
\sum_{x^{1/3-\varepsilon} \leq p \le x^{1/3}} S(\mathscr{A}_p; \mathscr{P}, x^{1/12}) 
\ll \dfrac{\varepsilon Cx}{\log ^2 x}. 
\end{equation}
Then, replacing sums over primes by integrals in the standard way yields
\begin{equation}\label{T4}
T_4 \le t_4\dfrac{Cx}{\log^2 x}(1+O(e^K\varepsilon)),
\end{equation}
where, using Lemma \ref{Ff}, we compute
\begin{equation}\label{t4}
\begin{aligned}
t_4=\ &\dfrac{12}{e^\gamma}\int_{1/12}^{2/7}\dfrac{1}{u}F\bigg(\dfrac{4/7-u}{1/12}\bigg)du+\dfrac{12}{e^\gamma}\int_{2/7}^{29/100}\dfrac{1}{u}F\bigg(\dfrac{2-6u}{1/12}\bigg)du\\
	&+\dfrac{12}{e^\gamma}\int_{29/100}^{1/3}\dfrac{1}{u}F\bigg(\dfrac{11/20-u}{1/12}\bigg)du\\
	=\ &\int_{1/12}^{2/7} \frac{14}{u(4-7u)}du+\int_{1/12}^{2/7} \frac{14du}{u(4-7u)}\int_2^{(41-84u)/7} \frac{\log (v-1)}{v} d v\\
	&+\int_{1/12}^{13/84} \frac{14du}{u(4-7u)}\int_2^{(27-84u)/7} \frac{\log (v-1)}{v} d v \int_{v+2}^{(41-84u)/7} \frac{1}{w} \log \frac{w-1}{v+1} dv\\
	&+\int_{2/7}^{29/100}\frac{1}{u(1-3u)}du + \int_{2/7}^{29/100}\frac{du}{u(1-3u)}\int_{2}^{23-72u}\frac{\log(v-1)}{v}dv\\
	&+\int_{29/100}^{1/3}\frac{120}{u(33-60u)}du + \int_{29/100}^{3/10}\frac{120du}{u(33-60u)}\int_{2}^{(26-60u)/5}\frac{\log(v-1)}{v}dv\\
	\le\ & 9.65921.
\end{aligned}
\end{equation}

Similarly, using only condition (C.1), we obtain
\begin{equation}\label{T5}
	T_5 \le t_5\dfrac{Cx}{\log^2 x} (1+O(e^K\varepsilon)),
	\end{equation}
	where 
	\begin{equation}\label{t5}
	\begin{aligned}
	t_5=\ &\dfrac{12}{e^\gamma}\int_{1/12}^{2/7}\dfrac{1}{u}F\bigg(\dfrac{4/7-u}{1/12}\bigg)du\\
	=\ &\int_{1/12}^{2/7} \frac{14}{u(4-7u)}du+\int_{1/12}^{2/7} \frac{14du}{u(4-7u)}\int_2^{(41-84u)/7} \frac{\log (v-1)}{v} d v\\
	&+\int_{1/12}^{13/84} \frac{14du}{u(4-7u)}\int_2^{(27-84u)/7} \frac{\log (v-1)}{v} d v \int_{v+2}^{(41-84u)/7} \frac{1}{w} \log \frac{w-1}{v+1} dv\\
	\le\ & 8.37847.
	\end{aligned}
	\end{equation}

For \( T_6 \) and \( T_7 \), since \( p_1^2 p_2 \le x^{4/7} \), we use Lemma \ref{Sieve-WF} and Lemma \ref{BFI} to get lower bounds. For \( i = 6, 7 \),
\begin{equation}\label{T6-7}
T_i \ge t_i\dfrac{Cx}{\log^2 x} (1+O(e^K\varepsilon)),
\end{equation}
where
\begin{equation}\label{t6}
	\begin{aligned}
	t_6=\ &\dfrac{12}{e^\gamma}\int_{1/12}^{1/7}\dfrac{1}{u}\int_{u}^{1/7}f\bigg(\dfrac{4/7-u-v}{1/12}\bigg)\dfrac{dv}{v}\\
	=\ &\int_{1/12}^{1/7}\dfrac{du}{u}\int_{u}^{1/7} \dfrac{14\log (41-84u-84v)-14\log 7}{v(4-7u-7v)}dv\\
	&+\int_{1/12}^{1/7}\dfrac{du}{u}\int_{u}^{1/7} \dfrac{14dv}{v(4-7u-7v)}\int_3^{(41-84u-84v)/7} \frac{dw}{w} \int_{2}^{w-1}\dfrac{\log (x-1)}{x}dx\\
	\ge\ & 0.96223,
	\end{aligned}
	\end{equation}
	and
	\begin{equation}\label{t7}
	\begin{aligned}
	t_7=\ &\dfrac{12}{e^\gamma}\int_{1/12}^{1/7}\dfrac{1}{u}\int_{1/7}^{2/7}f\bigg(\dfrac{4/7-u-v}{1/12}\bigg)\dfrac{dv}{v}\\
	=\ &\int_{1/12}^{1/7}\dfrac{du}{u}\int_{1/7}^{2/7} \dfrac{14\log (41-84u-84v)-14\log 7}{v(4-7u-7v)}dv\\
	&+\int_{1/12}^{2/21}\dfrac{du}{u}\int_{1/7}^{5/21-u} \dfrac{14dv}{v(4-7u-7v)}\int_3^{(41-84u-84v)/7} \frac{dw}{w} \int_{2}^{w-1}\dfrac{\log (x-1)}{x}dx\\
	\ge\ & 1.90760.
	\end{aligned}
	\end{equation}
The analysis in this subsection gives upper bounds for $T_4$ and $T_5$,  and lower bounds for $T_6$ and $T_7$.

\subsection{Evaluation of $T_3,T_8,T_{10}$}
In this subsection we will estimate the upper bounds of $T_3,T_8,T_{10}$, where
\[
T_3 =\sum\limits_{x^{3/7-\varepsilon}\le p\le x^{1/2}}S(\mathscr{A}_p ; \mathscr{P}, p), 
\]
\[
T_8 = \sum_{x^{2/7} \leq p_1 \le p_2 \le (x/p_1)^{1/2}} S(\mathscr{A}_{p_1 p_2}; \mathscr{P}(p_1), p_2),
\] 
and 
\[	
T_{10} = \sum_{x^{1/7} \leq p_1 \le x^{2/7} \leq p_2 \le (x/p_1)^{1/2}} S(\mathscr{A}_{p_1 p_2}; \mathscr{P}( p_1), (x/p_1 p_2)^{1/2}).
\] 

First we study \(T_8\).  
As in the treatment of \(G_8\), the restriction \(p_2 \ge p_1 \ge x^{2/7} > x^{1/4}\) forces the remaining factor to be a single prime \(p_3\).  
Hence \(T_8\) counts primes \(p \le x\) such that  
\[ 
p + 2 = p_1 p_2 p_3, \quad x^{2/7} \le p_1 \le p_2 \le (x/p_1)^{1/2}, \quad p_3 \ge p_2. 
\] 
Define the set  
\begin{equation}
\mathscr{B} := \{ b - 2 : b = p_1 p_2 p_3 \le x,\ x^{2/7} \le p_1 \le  p_2 \le p_3 \}, 
\end{equation}
so that  
\begin{equation}
T_8 = S(\mathscr{B}; \mathscr{P}, x^{1/2}) + O(x^{1/2}). 
\end{equation}
Then we define cuboid covering of $\mathscr{B}$. Let $\Delta := 1 + \log^{-4} x$. Cover $\mathscr{B}$ by cuboids
\begin{equation}
\mathscr{B}(t_1, t_2, t_3) 
:= \{ b - 2 : 
b = p_1 p_2 p_3 \leqslant x,\ p_i 
\in [ \Delta^{t_i}, \Delta^{t_i + 1} ) \text{ for } i=1, 2, 3\},
\end{equation}
where integers $t_1, t_2, t_3$ satisfy $x^{2/7} \leqslant \Delta^{t_1} \leqslant \Delta^{t_2} \leqslant \Delta^{t_3}$ and $\Delta^{t_1 + t_2 + t_3 + 3} \leqslant x$. 
Since
$$
x^{2/7} \leqslant p_2 \leqslant x^{5/14} \leqslant x^{2/5},
$$ 
we can apply Lemma \ref{Wu-level} with  
\begin{equation}
\alpha_m 
= \begin{cases} 1, & m = p_1 p_3, \\ 
0, & \text{otherwise}, 
\end{cases} 
\quad \beta_n 
= \begin{cases} 1, & n = p_2, \\ 
0, & \text{otherwise}. 
\end{cases}
\end{equation}
This	yields the inequality
	\begin{equation}
	\sum_{\substack{(q, 2) = 1}} \lambda_l^+(q) \bigg( | \mathscr{B}(t_1, t_2, t_3)_q | - \frac{| \mathscr{B}(t_1, t_2, t_3) |}{\varphi(q)} \bigg) \ll \frac{x}{\log^{18} x},
	\end{equation}
where $\lambda_l^+(q)$ is well-factorable of order 1 and level $Q = x^{\theta(t_2)}$ with 
\begin{equation}
\theta(t_2) = \dfrac{2+v}{4}. 
\end{equation}
From Lemma \ref{Sieve-WF} we obtain  
\begin{equation}
S( \mathscr{B}(t_1, t_2, t_3); \mathscr{P}, x^{1/2} ) 
\le \frac{2C(1 + O(\varepsilon))}{\theta(t_2) \log x} | \mathscr{B}(t_1, t_2, t_3) | + O\bigg( \frac{x}{\log^{18} x} \bigg).
\end{equation}
There are \(O(\log^{15} x)\) cuboids. Summing the estimate over all \((t_1, t_2, t_3)\) gives  
\begin{equation}
\begin{aligned}
\sum_{(t_1, t_2, t_3)} \frac{ |\mathscr{B}(t_1, t_2, t_3) |}{\theta(t_2)} 
=&\ \frac{x(1 + O(\varepsilon))}{\log x} \iint\limits_{2/7 \leqslant u \leqslant v \leqslant (1 - u)/2} \frac{dudv}{uv(1 - u - v)\theta(v)} \\ 
&\ + O\bigg( \frac{x}{\log^2 x} \bigg). 
\end{aligned}
\end{equation}
Combining these estimates, we obtain
\begin{equation}\label{T8}
T_8 \le t_8 \dfrac{Cx}{\log^2 x} (1+O(e^K\varepsilon)),
\end{equation}
where 
\begin{equation}\label{t8}
t_{8}=8\iint\limits_{2/7 \leqslant u \leqslant v \leqslant (1 - u)/2} \frac{dudv}{uv(1 - u - v)(2+v)}
\le 0.16037.
\end{equation}

Similarly for $T_{10}$ we have
\begin{equation}\label{T10}
T_{10} \le t_{10}\dfrac{Cx}{\log^2 x} (1+O(e^K\varepsilon)),
\end{equation}
where 
\begin{equation}\label{t10}
\begin{aligned}
	t_{10}=\ &8\iint\limits_{\substack{1/7 \le u \le 1/5\\2/7 \le  v \le 2/5}} \frac{dudv}{uv(1 - u - v)(2+v)}
	+8\iint\limits_{\substack{1/5 \le u \le 2/7\\2/7 \le  v \le (1-u)/2}} \frac{dudv}{uv(1 - u - v)(2+v)}\\
	&+2\iint\limits_{\substack{1/7 \le u \le 1/5\\2/5 \le  v \le (1-u)/2}} \frac{dudv}{uv(1 - u - v)(1-v)}\\
	\le\ & 1.70199.
\end{aligned}
\end{equation}

Finally for $T_3$, we just use $\theta(\log p/\log x)$ instead of $\theta(\log p_2/\log x)$ 
when we apply Lemma \ref{Wu-level}. Note that $3/7>2/5$, and therefore we take $\theta(u)=1-u$, giving 
\begin{equation}\label{T3}
T_{3} \le t_3\dfrac{Cx}{\log^2 x} (1+O(e^K\varepsilon)),
\end{equation}
where 
\begin{equation}\label{t3}
t_{3}=
2\int\limits_{3/7 \le u \le 1/2} \frac{du}{u(1 - u)^2}
=\dfrac{1}{2}+2\log \dfrac{4}{3}
\le 1.07537.
\end{equation}

\subsection{Evaluation of $T_9$}
In this subsection we will estimate an upper bound for $T_9$, where  
\begin{equation}
T_9 = \sum_{x^{1/12} \leq p_1 \le x^{1/3} \leq p_2 \le (x/p_1)^{1/2} } S(\mathscr{A}_{p_1 p_2}; \mathscr{P}(p_1), p_2).
\end{equation}
By the same argument as for $G_8$, the condition $p_3 \ge p_2\ge x^{1/3}\ge p_1\ge x^{1/12}$ 
forces the remaining factor to be a single prime $p_3$, and we obtain  
\begin{equation}
T_{9} = S(\mathscr{B}'; \mathscr{P}, x^{1/2}) + O(x^{1/2}),
\end{equation}
where
\begin{equation}\label{B'}
\mathscr{B}' := \{ b - 2 : b = p_1p_2p_3 \leq x,\ x^{1/12} \leq p_1 \leq x^{1/3} \leq p_2 \leq p_3,\ (p_1p_2p_3,2)=1 \}.
\end{equation}

To apply the switching principle efficiently, we partition $\mathscr{B}'$ into five subsets
$\mathscr{B}_1',\dots,\mathscr{B}_5'$, each corresponding to a different range of the variables
where a suitable level of distribution $\theta$ can be chosen. Concretely, we set
\[
\mathscr{B}' = \mathscr{B}_1' \cup \cdots \cup \mathscr{B}_5',
\]
where each $\mathscr{B}_i'$ is defined by the conditions in \eqref{B'} together with the
following additional restrictions:
\begin{itemize}
	\item $\mathscr{B}_1'$: $p_1 \le x^{1/10}$;
	\item $\mathscr{B}_2'$: $x^{1/10} < p_1 \le x^{1/5}$ and $p_2 \le x^{2/5}$;
	\item $\mathscr{B}_3'$: $x^{1/10} \le p_1 < x^{1/6}$ and $p_2> x^{3/8}p_1^{1/4}$;
	\item $\mathscr{B}_4'$: $p_2 > x^{2/5}$ and $p_2 \le x^{3/8}p_1^{1/4}$;
	\item $\mathscr{B}_5'$: $x^{1/5} \le p_1 < x^{1/3}$.
\end{itemize}
The five cases are mutually exclusive and cover the full range of $p_1,p_2$ appearing
in \eqref{B'}.

\begin{figure}[htbp]
	\centering
	\includegraphics[width=0.6\textwidth]{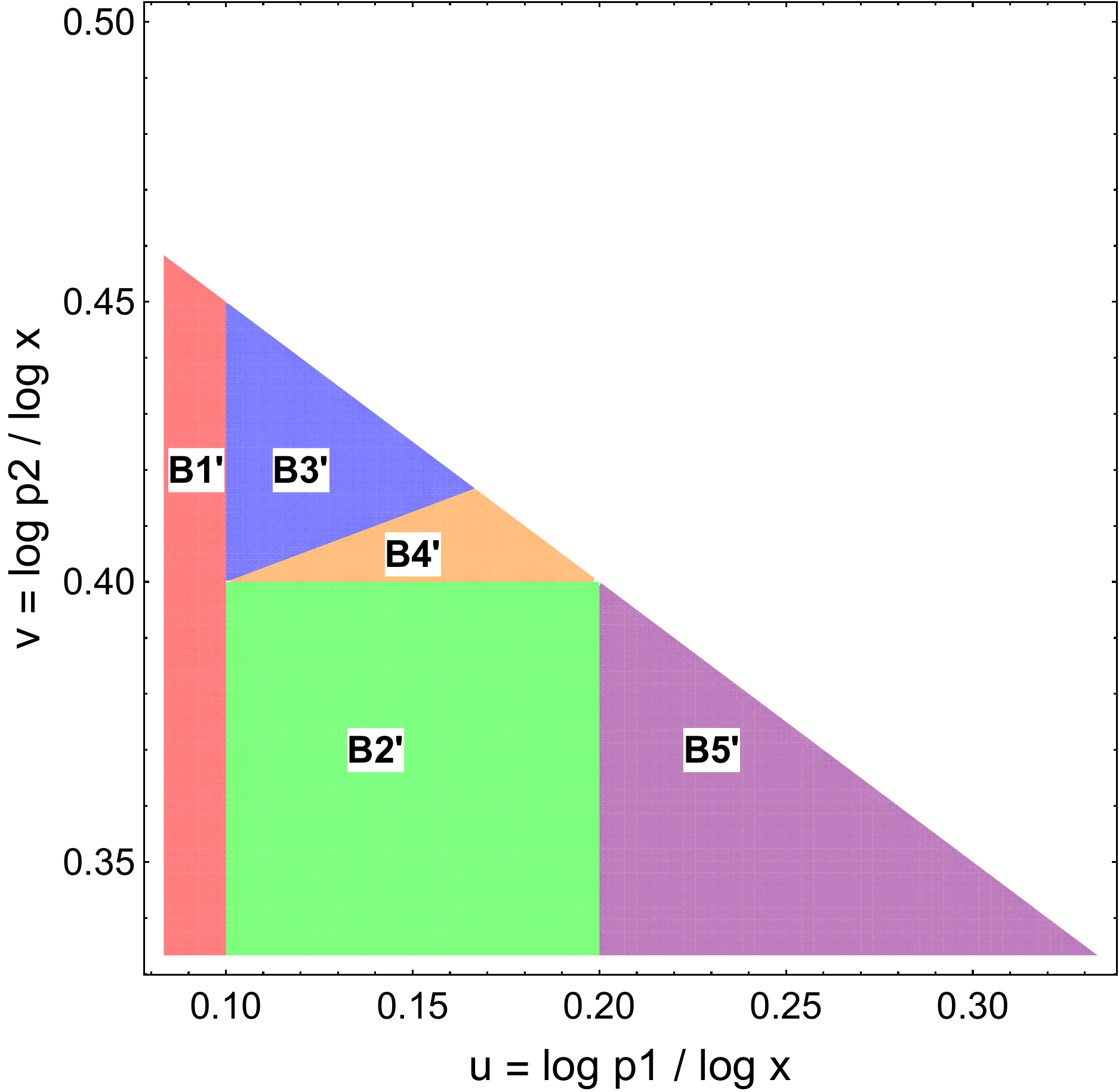}
	\caption{Partition of the set $\mathscr{B}'$ into $\mathscr{B}_1',\dots,\mathscr{B}_5'$ in the $(u,v)$-plane, where $u=\frac{\log p_1}{\log x}$ and $v=\frac{\log p_2}{\log x}$.}
	\label{fig:T9_split}
\end{figure}

Now we apply Lemma \ref{Wu-level} to each subset with the corresponding level function $\theta$:
\begin{equation}
\begin{aligned}
\theta\bigg( \frac{\log p_1}{\log x} \bigg) &= \frac{1 + 2\log p_1 / \log x}{2} && \text{for } \mathscr{B}_1', \\[2pt]
\theta\bigg( \frac{\log p_1}{\log x} \bigg) &= \frac{5 - 2\log p_1 / \log x}{8} && \text{for } \mathscr{B}_2' \text{ and } \mathscr{B}_3', \\[2pt]
\theta\bigg( \frac{\log p_2}{\log x} \bigg) &= 1 - \frac{\log p_2}{\log x} && \text{for } \mathscr{B}_4', \\[2pt]
\theta\bigg( \frac{\log p_2}{\log x} \bigg) &= \frac{2 + \log p_2 / \log x}{4} && \text{for } \mathscr{B}_5'.
\end{aligned}
\end{equation}
Summing the contributions of the five subsets, we obtain
\begin{equation}\label{T9}
T_9 \le  t_9 \dfrac{Cx}{\log^2 x} (1+O(e^K\varepsilon)),
\end{equation}
where the constant $t_9$ is the sum of the integrals corresponding to
$\mathscr{B}_1', \ldots, \mathscr{B}_5'$:
\begin{equation}\label{t9}
\begin{aligned}
t_9
&= 4\iint\limits_{\substack{1/12 \le u \le 1/10\\1/3 \le v \le (1-u)/2}} \frac{du\,dv}{uv(1-u-v)(1+2u)} 
+ 16\iint\limits_{\substack{1/10 \le u \le 1/5\\1/3 \le v \le 2/5}} \frac{du\,dv}{uv(1-u-v)(5-2u)} \\
&\quad + 16\iint\limits_{\substack{1/10 \le u \le 1/6\\(3+2u)/8 \le v \le (1-u)/2}} \frac{du\,dv}{uv(1-u-v)(5-2u)} 
+ 2\iint\limits_{\substack{2/5 \le v \le 5/12\\(8v-3)/2 \le u \le 1-2v}} \frac{du\,dv}{uv(1-u-v)(1-v)} \\ 
& \quad + 8\iint\limits_{\substack{1/5 \le u \le 1/3\\1/3 \le v \le (1-u)/2}} \frac{du\,dv}{uv(1-u-v)(2+v)} \\
&\le 2.03512.
\end{aligned}
\end{equation}
	
\subsection{Evaluation of $T_{11},T_{12}$}
The aim of this subsection is upper bounds for $T_{11}$ and $T_{12}$, where
\[
T_{11} = \sum_{x^{1/12} \leq p_1 \le p_2 \le p_3 \le p_4 \le x^{1/7} } S(\mathscr{A}_{p_1 p_2 p_3 p_4}; \mathscr{P}(p_1), p_2),
\]
and 
\[
T_{12} = \sum_{x^{1/12} \leq p_1 \le p_2 \le p_3 \le x^{1/7} \leq p_4 \le x^{2/7}} S(\mathscr{A}_{p_1 p_2 p_3 p_4}; \mathscr{P}(p_1), p_2).
\]
We use Buchstab's function to estimate $X$ instead of applying the Prime Number Theorem directly.
Similarly as before, we have, for $i=11,12$,
\begin{equation}\label{T11-12}
T_i \le  t_i \dfrac{Cx}{\log^2 x} (1+O(e^K\varepsilon)),
\end{equation}
where
\begin{equation}\label{t11}
\begin{aligned}
t_{11}
= &\ 4\int_{1/12}^{1/10}\frac{dt_1}{t_1(1+2t_1)}\int_{t_1}^{1/7}\frac{dt_2}{t_2^2}\int_{t_2}^{1/7}\frac{dt_3}{t_3}\int_{t_3}^{1/7} w\bigg(\frac{1-t_1-t_2-t_3-t_4}{t_2}\bigg)\frac{dt_4}{t_4}\\
&+16\int_{1/10}^{1/7}\frac{dt_1}{t_1(5-2t_1)}\int_{t_1}^{1/7}\frac{dt_2}{t_2^2}\int_{t_2}^{1/7}\frac{dt_3}{t_3}\int_{t_3}^{1/7} w\bigg(\frac{1-t_1-t_2-t_3-t_4}{t_2}\bigg)\frac{dt_4}{t_4}.
\end{aligned}
\end{equation}
To obtain a numerical bound we use the known estimates for $w(u)$.
For the range of $u$ appearing above we have 
\begin{equation}\label{w-Twin}
w(u)\le \begin{cases} 0.561522, & u \ge 3.3, \\ 
0.564383, & u \ge 3.  
\end{cases} 
\end{equation}
Indeed, 
\[
w(u)=\frac{1+\log(u-1)}{u}+\frac1u\int_2^{u-1}\frac{\log(t-1)}{t}dt
\] 
for $3\le u\le 4$, and $w(u)\le 0.561522$ for $u\ge 3.5$; the region $3\le u\le 3.5$ is handled by the monotonicity
argument as before, with a maximum at $u=3$ giving $w(3)\le0.564383$, and $w(3.3)\le0.561097<0.561522$.
Thus \eqref{w-Twin} holds. Inserting these bounds into $t_{11}$ yields
\[
\begin{aligned}
t_{11}
\le &\ 4\times 0.561522\times \int_{1/12}^{1/10}\frac{dt_1}{t_1(1+2t_1)}\int_{t_1}^{1/7}\frac{dt_2}{t_2^2}\int_{t_2}^{1/7}\frac{dt_3}{t_3}\int_{t_3}^{1/7}\frac{dt_4}{t_4}\\
&+16\times 0.564383\times\int_{1/10}^{1/7}\frac{dt_1}{t_1(5-2t_1)}\int_{t_1}^{1/7}\frac{dt_2}{t_2^2}\int_{t_2}^{1/7}\frac{dt_3}{t_3}\int_{t_3}^{1/7}\frac{dt_4}{t_4}\\
\le &\ 0.06507.
\end{aligned}
\]
This is the desired upper bound for $t_{11}$.

Set $W(u)=w(u)$ for $2\le u<3$, and $W(u)=0.564383\ge w(u)$ for $u\ge 3$.
For $i=12$, we have
\begin{equation*}
\begin{aligned}
t_{12}
= \ &4\int_{1/12}^{1/10}\frac{dt_1}{t_1(1+2t_1)}\int_{t_1}^{1/7}\frac{dt_2}{t_2^2}\int_{t_2}^{1/7}\frac{dt_3}{t_3}\int_{1/7}^{2/7} W\bigg(\frac{1-t_1-t_2-t_3-t_4}{t_2}\bigg)\frac{dt_4}{t_4}\\
&+16\int_{1/10}^{1/7}\frac{dt_1}{t_1(5-2t_1)}\int_{t_1}^{1/7}\frac{dt_2}{t_2^2}\int_{t_2}^{1/7}\frac{dt_3}{t_3}\int_{1/7}^{2/7} W\bigg(\frac{1-t_1-t_2-t_3-t_4}{t_2}\bigg)\frac{dt_4}{t_4}. 
\end{aligned}
\end{equation*}
Using the same numerical bounds we obtain
\begin{equation}\label{t12}
t_{12} \le 0.31841.
\end{equation}
This completes the estimation of $T_{11}$ and $T_{12}$.

\subsection{Completion of the proof of Theorem \ref{Twin1}} 
By \eqref{T1}, \eqref{T2}, \eqref{T4}, \eqref{T5}, \eqref{T6-7}, \eqref{T8}, \eqref{T10}, \eqref{T3}, \eqref{T9}, \eqref{T11-12}, we have
\begin{equation}
\begin{split}
4\pi_{1,1.75}(x)
\ge &\ \dfrac{Cx}{\log^2 x} (3\times 6.73680 + 3.84514 - 4\times 1.07537 - 9.65921  \\
& - 8.37847 + 0.96223 + 1.90760 - 2\times 0.16037 -2.03512 \\
& - 1.70199  - 0.06507 - 0.31841) \\ 
\ge &\ 0.14488\,\dfrac{Cx}{\log^2 x}.
\end{split}
\end{equation}
This proves Theorem \ref{Twin1}. 	
	
\section{Proof of Theorem~\ref{Goldbach2}}\label{sec7}

In Proposition \ref{Goldbachweight2} we specify $a=1.4$ and 
\begin{equation}\label{Cond/Goldsmall}
\alpha=1/12,\quad \beta=1/8,\quad \gamma=1/5,\quad \delta=1/4,\quad \tau=2/7-\varepsilon,\quad N>\exp(\exp(10\varepsilon^{-3})).
\end{equation} 
Thus we deduce from \eqref{Goldbachsmallw4} that 
\begin{equation*}
		\begin{split}
		4D_{1, 1.4}(N) \geq\ & 3H_1 + H_2 -4H_3- H_4- H_5 + H_6 + H_7- 2H_8 \\
		& - 2H_9 -2H_{10}- 2H_{11} - H_{12} - H_{13}-H_{14}-H_{15} \\ 
		& -H_{16}-H_{17} + O(N^{11/12}), 
		\end{split}
\end{equation*}
where $H_1, \ldots, H_{17}$ are as in Proposition \ref{Goldbachweight2}. We are going to estimate $H_1, \ldots, H_{17}$ 
in the subsections below.

\subsection{Evaluation of $H_1,H_2$}
In this subsection we will obtain lower bounds for $H_1$ and $H_2$ respectively, where
\[
H_1 = S(\mathscr{A}; \mathscr{P}(N), N^{1/12}), \quad 
H_2 = S(\mathscr{A}; \mathscr{P}(N), N^{1/8}).
\]

We want to apply Lemma \ref{Sieve-standard} to estimate $H_1$. Recalling \eqref{Def/A=Goldbach2}, 
we have at present that
\begin{equation}\label{New/A=Goldbach4}
\mathscr{A}=\{N-p:p<(1-\varepsilon)N\},
\end{equation}
and, before applying Lemma \ref{Sieve-standard}, we must check that the conditions
\eqref{def/SAPz}--\eqref{def/Vz} and \eqref{Iwacon} hold for the present $\mathscr{A}$. We compute that
\begin{equation}\label{H1X}
X=\Li(N-\varepsilon N)=\dfrac{N}{\log N}\bigg(1-\varepsilon+O\bigg(\dfrac{1}{\log N}\bigg)\bigg).
\end{equation}
We also have to compute the functions $\omega(p)$ and $V(z)$ defined as in \eqref{RIcon1} and \eqref{def/Vz}, respectively,
associated to this sequence $\mathscr{A}$, and
this has been done in classical literatures concerning the Goldbach Conjecture,
see for example \cite[\S1.3, Example 5]{HR74}, where it is proved that
\begin{equation}\label{New/ome/p-small}
\omega(p)= \frac{p}{p-1}
\end{equation}
for odd prime $p$, and $\omega(2)=0$. Thus \eqref{RIcon2} is satisfied.
It is also proved there that
\begin{equation}\label{Vz/Small}
V(z) = 2C(N) \frac{e^{-\gamma}}{\log z} \bigg(1 + O\bigg(\frac{1}{\log z}\bigg)\bigg),
\end{equation}
with
\begin{equation}\label{Come/Small}
\begin{aligned}
C(\omega) 
=&\  \prod_p \bigg(1 - \frac{\omega(p)}{p}\bigg) \bigg(1 - \frac{1}{p}\bigg)^{-1} \\ 
=&\ 2 \prod_{p > 2} \bigg(1 - \frac{1}{(p-1)^2}\bigg) \prod_{\substack{p \mid N \\ p > 2}} \frac{p-1}{p-2}
= 2C(N).
\end{aligned}
\end{equation}
Moreover, it is easy to check that for this $\omega(p)$,
\begin{equation}\label{Vz1/Vz2/H1}
\sum_{\substack{z_1 \le p < z_2}} \frac{\omega(p)}{p} = \log \frac{\log z_2}{\log z_1} 
+ O\left(\frac{1}{\log z_1}\right)
\end{equation}
holds for $z_2 \ge z_1 \ge 2$. This verifies \eqref{Iwacon}.

Now we assume Conjecture $\mathrm{WEH}(0.999)$, which provides the necessary control on the remainder terms.
We apply Lemma \ref{Sieve-standard} with
\[
z=N^{1/12},\quad Q=N^{0.999-\varepsilon},
\]
to deduce that
\begin{equation}\label{H1}
\begin{aligned}
H_1\ge\ & \bigg(1+O\bigg(\dfrac{1}{\log^{1/3}N}\bigg)\bigg)\dfrac{2C(N)N}{\log N}\dfrac{e^{-\gamma}}{\log (N^{1/12})}f(11.988)\\
\ge\ & h_1\dfrac{C(N)N}{\log^2 N} \bigg(1+O\bigg(\dfrac{1}{\log^{1/3}N}\bigg)\bigg),
\end{aligned}
\end{equation}
where
\begin{equation}\label{h1}
h_1=24e^{-\gamma}f(11.988)\ge 13.47361. 
\end{equation}
In the above we have used the facts that
\begin{equation}\label{f6,E-gamma}
f(11.988)\ge f(6)\ge 0.999895, \quad e^{-\gamma}\ge 0.561459.
\end{equation}
The above \eqref{H1} and \eqref{h1} give the desired lower bound for $H_1$.

By exactly the same idea, we can get that
\begin{equation}\label{H2}
\begin{aligned}
H_2\ge\ & \bigg(1+O\bigg(\dfrac{1}{\log^{1/3}N}\bigg)\bigg)\dfrac{2C(N)N}{\log N}\dfrac{e^{-\gamma}}{\log (N^{1/8})}f(7.992)\\
\ge\ & h_2\dfrac{C(N)N}{\log^2 N} \bigg(1+O\bigg(\dfrac{1}{\log^{1/3}N}\bigg)\bigg),
\end{aligned}
\end{equation}
where the constant $h_2$ can be computed by Lemma \ref{Ff} as 
\begin{equation}\label{h2}
h_2=16e^{-\gamma}f(7.992)\ge 8.98240.
\end{equation}
Also, \eqref{H2} and \eqref{h2} give the desired lower bound for $H_2$.

\subsection{Evaluation of $H_4,H_5,H_6,H_7$}
In this subsection we will prove upper bounds for $H_4$ and $H_5$, and lower bounds for $H_6$ and $H_7$, where
\[
H_4 = \sum_{\substack{N^{1/12} \leq p \le N^{1/4} \\ (p, N) = 1}} S(\mathscr{A}_p; \mathscr{P}(N), N^{1/12}),
\]
\[
H_5 = \sum_{\substack{N^{1/12} \leq p \le N^{1/5} \\ (p, N) = 1}} S(\mathscr{A}_p; \mathscr{P}(N), N^{1/12}),
\]
\[
H_6 = \sum_{\substack{N^{1/12} \leq p_1 \le p_2 \le N^{1/8} \\ (p_1 p_2, N) = 1}} S(\mathscr{A}_{p_1 p_2}; \mathscr{P}(N), N^{1/12}),
\]
and 
\[
H_7 = \sum_{\substack{N^{1/12} \leq p_1 \le N^{1/8} \le p_2 \le N^{1/5} \\ (p_1 p_2, N) = 1}} S(\mathscr{A}_{p_1 p_2}; \mathscr{P}(N), N^{1/12}).
\]

We start with an upper bound for $H_4$. We want to apply the upper bound part
in Lemma \ref{Sieve-standard} with
\[
z=N^{1/12},\quad Q=N^{0.999-\varepsilon}
\]
to the summand $S(\mathscr{A}_{p};\mathscr{P}(N), N^{1/12})$. The procedure is
similar to that for $H_1$, except that here we aim at an upper bound instead of a lower bound.
One checks that the conditions \eqref{def/SAPz}--\eqref{def/Vz} hold for the present sieving function
$S(\mathscr{A}_{p};\mathscr{P}(N), N^{1/12})$, and this can be done similarly as in
\eqref{New/A=Goldbach4}--\eqref{Vz1/Vz2/H1}. A simple calculation shows that,
instead of \eqref{New/A=Goldbach4} and \eqref{H1X}, now we have
\begin{equation}\label{New/A=/H4}
\mathscr{A}_p=\{N-q: q<(1-\varepsilon)N, \ q\equiv N \bmod{p}\}
\end{equation}
and
\begin{equation}\label{H1X/H4}
X=\frac{\Li(N-\varepsilon N)}{\varphi(p)}= \frac{1}{\varphi(p)} \frac{N}{\log N} \bigg(1-\varepsilon+O\bigg(\dfrac{1}{\log N}\bigg)\bigg),
\end{equation}
while \eqref{New/ome/p-small}--\eqref{Vz1/Vz2/H1} remain the same. Now the upper bound part
in Lemma \ref{Sieve-standard} yields
\begin{equation}
\begin{split}
& S(\mathscr{A}_{p};\mathscr{P}(N), N^{1/12}) \\ 
& \le 
\bigg(1+O\bigg(\dfrac{1}{\log^{1/3}N}\bigg)\bigg)\dfrac{N}{\log N}
\frac{V(N^{1/12})}{\varphi(p)} F\bigg(\frac{\log (N^{0.999-\varepsilon} / p)}{\log (N^{1/12})}\bigg) \\
&\quad  + \sum_{\substack{pq < D \\ pq \mid P(N^{1/12})}} |r_{pq}|.
\end{split}
\end{equation}
Inserting this into the expression of $H_4$, we see that the contribution from the last term above
can be estimated by Conjecture $\mathrm{WEH}(0.999)$ as
\begin{equation}
\begin{split}
&\ll \sum_{\substack{N^{1/12}\le p\le N^{1/4}\\(p, N)=1}} \sum_{\substack{pq < D \\ pq \mid P(N^{1/12})}} |r_{pq}|\\
&\ll \sum_{d\le N^{0.999-\varepsilon}} \mu^{2}(d) \max_{y\le N} \max_{(a,d)=1} \bigg|\pi(y; d, a) - \frac{\Li(y)}{\varphi(d)}\bigg| \\
&\ll \frac{N}{\log^3 N},
\end{split}
\end{equation}
which is acceptable. It follows that
\begin{equation}
H_4 \le \bigg(1+O\bigg(\dfrac{1}{\log^{1/3}N}\bigg)\bigg) \frac{N}{\log N} \sum_{\substack{N^{1/12} \le p \le N^{1/4} \\ (p, N) = 1}}
\frac{V(N^{1/12})}{\varphi(p)} F\bigg(\frac{\log (N^{0.999-\varepsilon} / p)}{\log (N^{1/12})}\bigg),
\end{equation}
where $V(z)$ is as \eqref{Vz/Small}. The condition $(p,N)=1$ can be simply removed. 
Hence we have
\begin{equation}
H_4 \le \bigg(1+O\bigg(\dfrac{1}{\log^{1/3}N}\bigg)\bigg)\dfrac{2C(N)N}{\log N}
\frac{e^{-\gamma}}{\log (N^{1/12})}\sum_{N^{1/12} \le p \le N^{1/4}}
\frac{1}{p-1} F\bigg(\frac{\log (N^{0.999-\varepsilon} / p)}{\log (N^{1/12})}\bigg).
\end{equation}
The last sum over $p$ can be estimated by the standard procedure of replacing sums over primes by integrals, as
\begin{equation*}
\begin{aligned}
\le\ & \int_{N^{1/12}}^{N^{1/4}} \frac{1}{t}F\bigg(\frac{0.999-\log t/\log N}{1/12}\bigg)\dfrac{dt}{\log t}\\
\le\ & \int_{1/12}^{1/4} \frac{1}{u}F(11.988-12u)du. 
\end{aligned}
\end{equation*}
It follows that 
\begin{equation}\label{H4}
\begin{aligned}
H_4\le h_4\dfrac{C(N)N}{\log^2 N}\bigg(1+O\bigg(\dfrac{1}{\log^{1/3}N}\bigg)\bigg),
\end{aligned}
\end{equation}
where the constant $h_4$ can be estimated by Lemma \ref{Ff} as
\begin{equation}\label{h4}
\begin{aligned}
h_{4}=\ & \dfrac{24}{e^{\gamma}}\int_{1/12}^{1/4} \frac{1}{u}F(11.988-12u)du
\le\dfrac{24}{e^{\gamma}}\cdot F(8.988)\cdot\int_{1/12}^{1/4} \frac{1}{u}du\\
\le\ & 14.80391. 
\end{aligned}
\end{equation}
In the above we have used the facts that 
\begin{equation}\label{F7,E-gamma}
F(8.988)\le F(7)\le 1.000005, \quad e^{-\gamma}\le 0.561460.
\end{equation}
The inequalities \eqref{H4} and \eqref{h4} give the desired upper bound for $H_4$.

An upper bound for $H_5$ can be obtained in exactly the same way. We have
\begin{equation}\label{H5}
\begin{aligned}
H_5\le\ & \bigg(1+O\bigg(\dfrac{1}{\log^{1/3}N}\bigg)\bigg)\dfrac{2C(N)N}{\log N}\dfrac{e^{-\gamma}}{\log (N^{1/12})}\sum_{N^{1/12} \le p \le N^{1/5}} \frac{F\bigg(\frac{\log (N^{0.999-\varepsilon}) / p)}{\log (N^{1/12})}\bigg)}{p-1} \\
\le\ & \bigg(1+O\bigg(\dfrac{1}{\log^{1/3}N}\bigg)\bigg)\dfrac{24}{e^{\gamma}}\dfrac{C(N)N}{\log ^2N}\int_{1/12}^{1/5} \frac{1}{u}F(11.988-12u)du\\
\le\ &h_5\dfrac{C(N)N}{\log^2 N} \bigg(1+O\bigg(\dfrac{1}{\log^{1/3}N}\bigg)\bigg),
\end{aligned}
\end{equation}
where the constant $h_5$ satisfies 
\begin{equation}\label{h5}
\begin{aligned}
h_{5}=&\ \dfrac{24}{e^{\gamma}}\int_{1/12}^{1/5} \frac{1}{u}F(11.988-12u)du
\le\dfrac{24}{e^{\gamma}} \cdot F(9.588) \cdot \int_{1/12}^{1/5} \frac{1}{u}du\\
\le&\ 11.79703,
\end{aligned}
\end{equation}
on using \eqref{F7,E-gamma} again. 
The above two formulas \eqref{H5} and \eqref{h5} give the desired upper bound for $H_5$.

In the following we work on the lower bounds for $H_6$ and $H_7$.
Note that the condition $(p_{1}p_{2},N)=1$ can be removed by the estimate
\begin{equation}
\sum_{\substack{p_{1}\mid N\\ p_{1}\ge N^{\kappa}}}\sum_{N^{\kappa}\le p_{2}<N^{\sigma}}\frac{N}{p_{1}p_{2}}+\sum_{N^{\kappa}\le p_{1}<N^{\sigma}}\sum_{\substack{p_{2}\mid N\\ p_{2}\ge N^{\kappa}}}\frac{N}{p_{1}p_{2}}\ll N^{1-\kappa}\log N,
\end{equation}
where $1/30<\kappa<\sigma<1/2$. As before we can deduce, from Lemma \ref{Sieve-standard} and Conjecture $\mathrm{WEH}(0.999)$, that
\begin{equation*}
\begin{aligned}
H_6 
\ge&\ 
\bigg(1+O\bigg(\dfrac{1}{\log^{1/3}N}\bigg)\bigg)\dfrac{2C(N)N}{\log N} \dfrac{e^{-\gamma}}{\log (N^{1/12})} \\ 
&\ \times \sum_{N^{1/12} \leq p_1 \le p_2 \le N^{1/8}} \frac{1}{p_1p_2}f\bigg(\frac{0.999-\log (p_1p_2)/\log N}{1/12}\bigg). 
\end{aligned}
\end{equation*}
The double sums over $p_1$ and $p_2$ is 
\begin{equation*}
\begin{aligned}
\ge&\ 
\int_{N^{1/12}}^{N^{1/8}} \dfrac{dt_1}{t_1\log t_1}\int_{t_1}^{N^{1/8}}f\bigg(\frac{0.999-\log (t_1t_2)/\log N}{1/12}\bigg)\dfrac{dt_2}{t_2\log t_2}\\
\ge&\ \int_{1/12}^{1/8}\dfrac{du}{u}\int_{u}^{1/8} f(11.988-12u-12v)\dfrac{dv}{v}. 
\end{aligned}
\end{equation*}
It follows that 
\begin{equation}\label{H6}
\begin{aligned}
H_6\ge h_6\dfrac{C(N)N}{\log^2 N}\bigg(1+O\bigg(\dfrac{1}{\log^{1/3}N}\bigg)\bigg),
\end{aligned}
\end{equation}
with 
\begin{equation}\label{h6}
\begin{aligned}
h_{6}=\ & \dfrac{24}{e^{\gamma}} 
\int_{1/12}^{1/8}\dfrac{du}{u}\int_{u}^{1/8}f(11.988-12u-12v)\dfrac{dv}{v}\\
\ge\ & \dfrac{24}{e^{\gamma}} \cdot f(8.988) \cdot \int_{1/12}^{1/8}\dfrac{du}{u}\int_{u}^{1/8}\dfrac{dv}{v}\\
\ge\ & 1.10765.
\end{aligned}
\end{equation}
Here we have used \eqref{f6,E-gamma} again. 
The above \eqref{H6} and \eqref{h6} give the desired lower bound for $H_6$.

A lower bound for $H_7$ can be obtained in the same way.
Thus, we have
\begin{equation}\label{H7}
\begin{aligned}
H_7
\ge\ & \bigg(1+O\bigg(\dfrac{1}{\log^{1/3}N}\bigg)\bigg)\dfrac{2C(N)N}{\log N}\dfrac{e^{-\gamma}}{\log (N^{1/12})}\\
&\times \sum_{N^{1/12} \leq p_1 \le N^{1/8} \le p_2 \le N^{1/5} } \frac{1}{p_1p_2}f\bigg(\frac{0.999-\log (p_1p_2)/\log N}{1/12}\bigg)\\
\ge\ &h_7\dfrac{C(N)N}{\log^2 N} \bigg(1+O\bigg(\dfrac{1}{\log^{1/3}N}\bigg)\bigg),
\end{aligned}
\end{equation}
where the constant $h_7$ is bounded by \eqref{f6,E-gamma} as 
\begin{equation}\label{h7}
\begin{aligned}
h_{7}=\ & \dfrac{24}{e^{\gamma}}\int_{1/12}^{1/8}\dfrac{du}{u}\int_{1/8}^{1/5}f(11.988-12u-12v)\dfrac{dv}{v}\\
\ge\ & \dfrac{24}{e^{\gamma}}\cdot f(8.088)\cdot \int_{1/12}^{1/8}\dfrac{du}{u}\int_{1/8}^{1/5}\dfrac{dv}{v}\\
\ge\ & 2.56793.
\end{aligned}
\end{equation}
The desired lower bound for $H_7$ now follows from \eqref{H7} and \eqref{h7}.

\subsection{Evaluation of $H_{9},H_{10},H_{11},H_{13}$}
In this subsection we will estimate upper bounds for $H_9,H_{10},H_{11},H_{13}$, where
\[
H_{9}  = \sum_{\substack{N^{1/4} \leq p_1 \le p_2 \le N^{2/7-\varepsilon} \\ (p_1 p_2, N) = 1}} S(\mathscr{A}_{p_1 p_2}; \mathscr{P}(N p_1), p_2),
\]
\[
H_{10} = \sum_{\substack{N^{1/5} \leq p_1 \le N^{2/7-\varepsilon}\le p_2 \le (N/p_1)^{1/2} \\ (p_1 p_2, N) = 1}} S(\mathscr{A}_{p_1 p_2}; \mathscr{P}(N p_1), p_2),
\]
\[
H_{11} = \sum_{\substack{N^{1/4} \leq p_1 \le N^{2/7-\varepsilon}\le p_2 \le (N/p_1)^{1/2} \\ (p_1 p_2, N) = 1}} S(\mathscr{A}_{p_1 p_2}; \mathscr{P}(N p_1), p_2),
\]
and 
\[
H_{13} = \sum_{\substack{N^{1/8} \leq p_1 \le N^{1/5} \leq p_2 \le (N/p_1)^{1/2} \\ (p_1 p_2, N) = 1}} S(\mathscr{A}_{p_1 p_2}; \mathscr{P}(N p_1), (N/p_1 p_2)^{1/2}).
\]

We shall only majorize $H_9$, and the others can be treated similarly.
Since the exponent of the least prime factor is $1/4$, the quantity $H_9$ equals the number of primes $p \leqslant N$ such that $N - p = p_1p_2p_3$ with
\begin{equation}
N^{1/4} \le p_1 \le p_2 \le (N/p_1)^{1/2},\quad p_3 \ge p_2,\quad p_1p_2p_3\ge \varepsilon N,\quad (p_1p_2p_3, N) = 1.
\end{equation}
Define the sets
\begin{equation}
\begin{aligned}
\mathscr{M} &:= \{ m : m = p_1p_2,\ N^{1/4} \le p_1 \le p_2 \le (N/p_1)^{1/2},\ (p_1p_2, N) = 1 \}, \\
\mathscr{B} &:= \{ b\le (1-\varepsilon) N : b = N - mp,\ m \in \mathscr{M},\ p \leqslant N/m \}.
\end{aligned}
\end{equation}
It is clear that
\[
H_9 \le S(\mathscr{B}; \mathscr{P}(N), N^{1/2}) + O(N^{1/2}).
\]

Applying Lemma \ref{Sieve-standard} with
\[
X = \sum_{m \in \mathscr{M}} \bigg(\Li\bigg( \frac{N}{m} \bigg)-\Li\bigg( \frac{\varepsilon N}{m} \bigg)\bigg),\quad Q = N^{0.999-\varepsilon},
\]
we obtain
\begin{equation}\label{H9old}
H_9 \leqslant 4.00401\frac{C(N) X}{\log N} (1 + O(\varepsilon)) + O(\sqrt{N} + R_1 +R_2),
\end{equation}
where
$$ 
R_1 := \sum_{\substack{q \leqslant N^{0.999-\varepsilon} \\ (q,N)=1}} \mu^2(q) \biggl| \sum_{\substack{m \in \mathscr{M} \\ (q,m)=1}} \biggl( \sum_{\substack{mp \leqslant N \\ mp \equiv N \bmod q}} 1 - \frac{\Li(N/m)}{\varphi(q)} \biggr) \biggr|, 
$$
and 
$$ 
R_2 := \sum_{\substack{q \leqslant N^{0.999-\varepsilon}\\ (q,N)=1}} \frac{\mu^2(q)}{\varphi(q)} \sum_{\substack{m \in \mathscr{M}\\ (q,m)>1}} \Li\biggl(\frac{N}{m}\biggr). $$
	Let $f(m)$ be the characteristic function of $\mathscr{M}$. Since $m \leqslant N^{3/4}$ for $m \in \mathscr{M}$, $\mathrm{WEH}(0.999)$ implies
\begin{equation}
R_1 = \sum_{\substack{q \leqslant N^{0.999-\varepsilon} \\ (q,N)=1}} \mu^2(q) \biggl| \sum_{\substack{m \leqslant N^{3/4} \\ (q,m)=1}} f(m) \biggl( \sum_{\substack{mp \leqslant N \\ mp \equiv N \bmod q}} 1 - \frac{\Li(N/m)}{\varphi(q)} \biggr) \biggr| \ll \frac{N}{\log^{3} N}.
\end{equation}
And, similarly as \eqref{Chen-R2}, 
	\begin{equation}
	R_2 \ll N^{3/4} \log^{2} N.
	\end{equation}

The asymptotic expression for $X$ follows from the Prime Number Theorem. Thus, 
\begin{equation}\label{H9x}
\begin{aligned}
X &= \bigg(1-\varepsilon+O\bigg(\dfrac{1}{\log N}\bigg)\bigg) 
\sum_{N^{1/4} \leqslant p_1 \le p_2 \le N^{2/7-\varepsilon}} \frac{N}{p_1p_2 \log(N/p_1p_2)} \\
&= \bigg(1-\varepsilon+O\bigg(\dfrac{1}{\log N}\bigg)\bigg)\frac{N}{\log N} \int_{1/4}^{2/7}\dfrac{dt}{t}\int_{t}^{2/7}\dfrac{du}{u(1-t-u)}.
\end{aligned}
\end{equation}
Collecting \eqref{H9x} into \eqref{H9old}, we have
\begin{equation}\label{H9}
H_9 \le h_9\dfrac{C(N)N}{\log^2 N} \bigg(1+O\bigg(\dfrac{1}{\log^{1/3}N}\bigg)\bigg),
\end{equation}
where
\begin{equation}\label{h9}
h_{9}=4.00401\int_{1/4}^{2/7}\dfrac{dt}{t}\int_{t}^{2/7}\dfrac{du}{u(1-t-u)}
\le 0.07684.
\end{equation}
Similarly we have for $i=10,11,13$ that 
\begin{equation}\label{H10-13}
H_i \le h_i\dfrac{C(N)N}{\log^2 N}\bigg(1+O\bigg(\dfrac{1}{\log^{1/3}N}\bigg)\bigg),
\end{equation}
where
\begin{equation}\label{h10}
h_{10}=4.00401\int_{1/5}^{2/7}\dfrac{dt}{t}\int_{2/7}^{(1-t)/2}\dfrac{du}{u(1-t-u)}
\le 0.94673,
\end{equation}
\begin{equation}\label{h11}
h_{11}=4.00401\int_{1/4}^{2/7}\dfrac{dt}{t}\int_{2/7}^{(1-t)/2}\dfrac{du}{u(1-t-u)}
\le 0.32602,
\end{equation}
and 
\begin{equation}\label{h13}
h_{13}=4.00401\int_{1/8}^{1/5}\dfrac{dt}{t}\int_{1/5}^{(1-t)/2}\dfrac{du}{u(1-t-u)}
\le 2.60454. 
\end{equation}
In the above we note that the sieve size in $H_{13}$ is $(N/p_1p_2)^{1/2}$, 
which restricts the cofactor to a single prime just as in the analysis of $G_8$. 

\subsection{Evaluation of $H_{3},H_{8},H_{12},H_{14},H_{15},H_{16},H_{17}$}
In this subsection we will estimate upper bounds for $H_{3},H_{8},H_{12},H_{14},H_{15},H_{16},H_{17}$, where
\[
H_3 =\sum_{\substack{N^{2/7-\varepsilon}\le p\le N^{1/2} \\(p, N)=1}}S(\mathscr{A}_p ; \mathscr{P}(N), p),
\]
\[
H_8 = \sum_{\substack{N^{1/5} \leq p_1 \le p_2 \le N^{2/7-\varepsilon} \\ (p_1 p_2, N) = 1}} S(\mathscr{A}_{p_1 p_2}; \mathscr{P}(N p_1), p_2),
\]
\[
H_{12} = \sum_{\substack{N^{1/12} \leq p_1 \le N^{1/4}\leq p_2 \le (N/p_1)^{1/2} \\ (p_1 p_2, N) = 1}} S(\mathscr{A}_{p_1 p_2}; \mathscr{P}(N p_1), p_2),
\]
\[
H_{14} = \sum_{\substack{N^{1/8} \le p_1\le N^{1/5}\\N^{1/4} \le p_2\le p_3\le (N/p_1p_2)^{1/2} \\ (p_1 p_2 p_3, N) = 1}}  S(\mathscr{A}_{p_1 p_2 p_3}; \mathscr{P}(Np_1), p_2),
\]
\[
H_{15} = \sum_{\substack{N^{1/8} \le p_1\le N^{1/5} \le p_2\le N^{1/4}\le p_3\le (N/p_1p_2)^{1/2} \\ (p_1 p_2 p_3, N) = 1}}  S(\mathscr{A}_{p_1 p_2 p_3}; \mathscr{P}(Np_1), p_2),
\]
\[
H_{16} = \sum_{\substack{N^{1/12} \leq p_1 \le p_2 \le p_3 \le p_4 \le N^{1/8} \\ (p_1 p_2 p_3 p_4, N) = 1}} S(\mathscr{A}_{p_1 p_2 p_3 p_4}; \mathscr{P}(Np_1), p_2),
\]
and 
\[
H_{17} = \sum_{\substack{N^{1/12} \leq p_1 \le p_2 \le p_3 \le N^{1/8} \leq p_4 \le N^{1/5}\\ (p_1 p_2 p_3 p_4, N) = 1}} S(\mathscr{A}_{p_1 p_2 p_3 p_4}; \mathscr{P}(Np_1), p_2).
\]

We use Buchstab’s function to estimate the integrals instead of applying the Prime Number Theorem directly.  
Let $w(u)$ be the Buchstab function, and set $W(u)=w(u)$ for $1\le u<3$, and $W(u)=0.564383\ge w(u)$ for $u\ge 3$.
The estimates $w(u)\le 0.561522$ for $u\ge 3.4$ follow from the monotonicity of $w$ as before:
since 
$$
w(u)=\frac{1+\log(u-1)}{u}+\frac{1}{u}\int_{2}^{u-1}\frac{\log(t-1)}{t}dt
$$ 
for $3\le u\le 4$, and $w(u)\le 0.561522$ for $u\ge 3.5$,
it suffices to check $w(3.4)\le 0.560862<0.561522$, giving the required bound for $u\ge 3.4$. The bound $W(u)\le 0.564383$ for $u\ge 3$ comes from $w(3)\le 0.564383$. These numerical bounds are used throughout. 

As before, for \(i = 8, 12, 14, 15, 16, 17\),
\begin{equation}\label{H8-17}
H_i \le h_i\dfrac{C(N)N}{\log^2 N}\bigg(1+O\bigg(\dfrac{1}{\log^{1/3}N}\bigg)\bigg),
\end{equation}
where
\begin{equation}\label{h8}
\begin{aligned}
h_{8}=\ &4.00401\int_{1/5}^{2/7}\dfrac{dt_1}{t_1}\int_{t_1}^{2/7}W\bigg(\dfrac{1-t_1-t_2}{t_2}\bigg)\dfrac{dt_2}{t_2^2}\\
\le\ & 0.55225,
\end{aligned}
\end{equation}
\begin{equation}\label{h12}
\begin{aligned}
h_{12}=\ &4.00401\int_{1/12}^{1/4}\dfrac{dt_1}{t_1}\int_{1/4}^{(1-t_1)/2}W\bigg(\dfrac{1-t_1-t_2}{t_2}\bigg)\dfrac{dt_2}{t_2^2}\\
\le \ & 4.67442,
\end{aligned}
\end{equation}
\begin{equation}\label{h14}
\begin{aligned}
h_{14}=\ &4.00401\int_{1/8}^{1/5}\dfrac{dt_1}{t_1}\int_{1/4}^{(1-t_1)/3}\dfrac{dt_2}{t_2^2}\int_{t_2}^{(1-t_1-t_2)/2}W\bigg(\dfrac{1-t_1-t_2-t_3}{t_2}\bigg)\dfrac{dt_3}{t_3}\\
\le \ &0.06040,
\end{aligned}
\end{equation}
\begin{equation}\label{h15}
\begin{aligned}
h_{15}=\ &4.00401\int_{1/8}^{1/5}\dfrac{dt_1}{t_1}\int_{1/5}^{1/4}\dfrac{dt_2}{t_2^2}\int_{1/4}^{(1-t_1-t_2)/2}W\bigg(\dfrac{1-t_1-t_2-t_3}{t_2}\bigg)\dfrac{dt_3}{t_3}\\
\le \ &0.25706,
\end{aligned}
\end{equation}
\begin{equation}\label{h16}
\begin{aligned}
h_{16}=\ &4.00401\times 0.561522\times\int_{1/12}^{1/8}\dfrac{dt_1}{t_1}\int_{t_1}^{1/8}\dfrac{dt_2}{t_2^2}\int_{t_2}^{1/8}\dfrac{dt_3}{t_3}\int_{t_3}^{1/8}\dfrac{dt_4}{t_4}\\
\le \ &0.02592,
\end{aligned}
\end{equation}
and 
\begin{equation}\label{h17}
\begin{aligned}
h_{17}=\ &4.00401\times 0.561522\times\int_{1/12}^{1/8}\dfrac{dt_1}{t_1}\int_{t_1}^{1/8}\dfrac{dt_2}{t_2^2}\int_{t_2}^{1/8}\dfrac{dt_3}{t_3}\int_{1/8}^{1/5}\dfrac{dt_4}{t_4}\\
\le\ & 0.11551.
\end{aligned}
\end{equation}

Finally we treat $H_3$. We want to apply the upper bound part
in Lemma \ref{Sieve-standard} with
\[
z=p,\quad Q=N^{0.999-\varepsilon}
\]
to the summand $S(\mathscr{A}_{p};\mathscr{P}(N), p)$. The procedure is 
similar to that for $H_4$, and therefore 
\begin{equation}
H_3 \le \bigg(1+O\bigg(\dfrac{1}{\log^{1/3}N}\bigg)\bigg) \frac{N}{\log N} \sum_{\substack{N^{2/7} \le p \le N^{1/2} \\ (p, N) = 1}}
\frac{V(p)}{\varphi(p)} F\bigg(\frac{\log (N^{0.999-\varepsilon} / p)}{\log p}\bigg),
\end{equation}
where $V(p)$ is as \eqref{Vz/Small}. The condition $(p,N)=1$ can be simply removed. 
Hence we have
\begin{equation}
H_3 \le \bigg(1+O\bigg(\dfrac{1}{\log^{1/3}N}\bigg)\bigg)\dfrac{2C(N)N}{\log N}
\sum_{N^{2/7} \le p \le N^{1/2}}\frac{e^{-\gamma}}{\log p}
\frac{1}{p-1} F\bigg(\frac{\log (N^{0.999-\varepsilon} / p)}{\log p}\bigg).
\end{equation}
By standard procedure of replacing sum over primes by integral, 
the last sum over $p$ can be estimated as 
\begin{equation*}
\begin{aligned}
\le\ & \int_{N^{2/7}}^{N^{1/2}} \frac{e^{-\gamma}}{t\log  t}F\bigg(\frac{0.999-\log t/\log N}{\log t/\log N}\bigg)\dfrac{dt}{\log t}\\
\le\ & \int_{2/7}^{1/2} \frac{e^{-\gamma}}{u^2}F\bigg(\frac{0.999-u}{u}\bigg)du.
\end{aligned}
\end{equation*}
It follows that 
\begin{equation}\label{H3}
\begin{aligned}
H_3\le h_3\dfrac{C(N)N}{\log^2 N}\bigg(1+O\bigg(\dfrac{1}{\log^{1/3}N}\bigg)\bigg),
\end{aligned}
\end{equation}
where the constant $h_3$ can be estimated by Lemma \ref{Ff} as
\begin{equation}\label{h3}
\begin{aligned}
h_{3}=\ & \dfrac{2}{e^{\gamma}}\int_{2/7}^{1/2} \frac{1}{u^2}F\bigg(\frac{0.999-u}{u}\bigg)du
\le 4\int_{2/7}^{1/2} \frac{1}{u(0.999-u)}du\\
\le\ & 3.67124. 
\end{aligned}
\end{equation}
The inequalities \eqref{H3} and \eqref{h3} give the desired upper bound for $H_3$.	

\subsection{Completion of the proof}
By \eqref{H1}, \eqref{H2}, \eqref{H4}, \eqref{H5}, \eqref{H6}, \eqref{H7}, \eqref{H9}, \eqref{H10-13}, 
\eqref{H8-17}, \eqref{H3}, we have
\begin{equation}
\begin{split}
4D_{1,1.4}(N)\ge\ & \dfrac{C(N)N}{\log^2(N)} 
(3\times 13.47361+8.98240-4\times 3.67124-14.80391\\
&-11.79703 +1.10765+2.56793-2\times 0.55225-2\times 0.07684\\
&-2\times 0.94673 -2\times 0.32602-4.67442-2.60454-0.06040\\ 
& -0.25706-0.02592-0.11551)\\ 
\ge\ & 0.25138\dfrac{C(N)N}{\log^2 N},
\end{split}
\end{equation}
from which Theorem \ref{Goldbach2} follows. 
	
\section{Further remarks}\label{sec8}
\subsection{What is the difference between Proposition $(1+a)$ and Proposition $(1+2)$?} 
In the proof of our Proposition $(1+a)$, a term of the form
\[
G:=\sum_{\substack{N^{\tau}\le p\le N^{1/2} \\ (p,N)=1}} S(\mathscr{A}_p;\mathscr{P}(N),p)
\]
always appears. 
A difficulty arises when the upper limit of the summation approaches $N^{1/2}$.
If we try to estimate an upper bound without using the switching principle, we obtain
\begin{equation}
\begin{aligned}
G\le\ & (1+O(E))\dfrac{2C(N)N}{\log N}\sum_{N^{\tau} \le p\le N^{1/2}}\dfrac{e^{-\gamma}}{\log p} \frac{F\big(\frac{\log (\sqrt{N} / p)}{\log p}\big)}{p-1} \\
\le\ & (1+O(E))\dfrac{2C(N)N}{e^{\gamma}\log ^2N}\sum_{N^{\tau} \le p\le N^{1/2}}\dfrac{\log N}{p\log p}F\bigg(\frac{1/2-\log p/\log N}{\log p/\log N}\bigg)\\
\le\ & (1+O(E))\dfrac{2C(N)N}{e^\gamma\log ^2N}\int_{N^{\tau}}^{N^{1/2}} \dfrac{\log N}{t\log t}F\bigg(\frac{1/2-\log t/\log N}{\log t/\log N}\bigg)\dfrac{dt}{\log t}\\
\le\ & (1+O(E))\dfrac{2C(N)N}{e^\gamma\log ^2N}\int_{\tau}^{1/2} \dfrac{1}{u^2}F\bigg(\frac{1/2-u}{u}\bigg)du\\
\le\ &  g \dfrac{C(N)N}{\log^2 N}(1+O(E)),
\end{aligned}
\end{equation}
where
\begin{equation}
g:= \frac{2}{e^\gamma}\int_{\tau}^{1/2} \dfrac{1}{u^2}F\bigg(\frac{1/2-u}{u}\bigg)du=+\infty,
\end{equation}
where $E$ represents the corresponding error term. 
This shows that when the summation range includes large primes, Chen's switching principle plays an indispensable role; without it, we would not be able to make the upper bound of this term have the correct order of growth of its expected value. 
In the other terms of Proposition $(1+2)$ a more suitable choice of weights avoids this phenomenon; at least each term can be bounded above with the correct order of magnitude. 
However, because the constants in those bounds are poor, avoiding the switching principle still fails to give a nontrivial lower bound for $D_{1,2}(N)$.

Thus, in the proof of Proposition $(1+2)$ the switching principle is decisive, whereas it is not 
necessary for Proposition $(1-2)$. 
To this day we still have no proof of Proposition $(1+2)$ that completely circumvents the switching principle. 
Greaves\cite{Gr82} showed that when the level of distribution reaches $1/1.936$ the switching principle can be avoided; this observation also underlies the work of Fouvry and Grupp\cite{FG86}, who, after applying the Bombieri--Friedlander--Iwaniec $4/7$‑level \cite{BFI86}, managed to prove Proposition $(1-2)$ without using the switching principle.

\subsection{When analytic tools advance, how do combinatorial tools keep pace with them?} 
Since R\'enyi~\cite{Ren48}, the nearly eighty-year process of improving the level of distribution has deeply influenced the weighted sieve. 
The step-by-step approach to the Goldbach Conjecture is inseparable from the synergy between analytic and combinatorial methods.
We explain what the weighted sieve method does, and why it is necessary to design different weighted sieves. 

\subsubsection{The switching principle can only be used to estimate upper bounds}
In the weighted sieve method, the fundamental reason why we employ multiple Buchstab iterations to obtain extremely complicated weights is that applying the switching principle yields only upper bounds, never lower bounds.

We illustrate this with one term from the proof:
\begin{equation}
G_8 = \sum_{\substack{N^{3/11} \le p_1 \le p_2 \le (N/p_1)^{1/2} \\ (p_1 p_2, N) = 1}} S(\mathscr{A}_{p_1 p_2}; \mathscr{P}(N p_1), p_2).
\end{equation}
Because the exponent of the least prime factor satisfies $3/11>1/4$, the quantity $G_8$ equals the number of primes $p\le N$ such that $N-p = p_1p_2p_3$ with
\begin{equation}
N^{3/11} \le p_1 \le p_2 \le (N/p_1)^{1/2},\quad p_3 \ge p_2,\quad p_1p_2p_3\ge \varepsilon N,\quad (p_1p_2p_3, N) = 1.
\end{equation}
Define the sets
\begin{equation}
\begin{aligned}
\mathscr{M} &:= \{ m\ge \varepsilon N : m = p_1p_2,\ N^{3/11} \le p_1 \le p_2 \le (N/p_1)^{1/2},\ (p_1p_2, N) = 1 \}, \\
\mathscr{B} &:= \{ b\le (1-\varepsilon) N : b = N - mp,\ m \in \mathscr{M},\ p \le N/m \}.
\end{aligned}
\end{equation}
Clearly
\[
G_8 = S(\mathscr{B}; \mathscr{P}(N), N^{1/2}) + O(N^{1/2}).
\]
In our paper we obtained an upper bound for $G_8$, whereas its lower bound is trivial: $G_8\ge 0$, because $f(s)=0$ for $s\le 2$.

In retrospect, in our proof of Proposition $(1+1.9)$ we used a weighted inequality with twelve terms instead of the simpler Proposition~\ref{Goldbachbig} which states that 

\begin{equation}\label{Exp/2D}
\begin{split}
2D_{1,a}(N)
\ge\ & 2S_1(\kappa)-2S_2(\tau)-S_3(\kappa,\sigma)-2S_4(\sigma)\\
&-S_5(\kappa,\sigma)+S_6(\kappa,\sigma)+O(N^{1-\kappa}) 
\end{split}
\end{equation}
with 
\[
S_6(\kappa,\sigma)= \sum_{\substack{N^\kappa \le p_1 \le p_2 \le p_3 \le N^\sigma \\ (p_1 p_2 p_3, N) = 1}} S(\mathscr{A}_{p_1 p_2 p_3}; \mathscr{P}(N p_1), p_2). 
\]
We made this choice since we were unable to derive a lower bound for \(S_6\) by switching principle. 
If we estimated \eqref{Exp/2D} directly, we would get a lower bound worse than the trivial one. 

\subsubsection{At different scales of prime size, the level of prime distribution varies} 
Looking back at Lemmas \ref{Wu-level}, \ref{Fouvry-new}, and \ref{Fouvry-level}, we see that the level of distribution varies depending on the range in which they are applied. 
Consequently, when using the weighted sieve method, we often carefully select parameters to increase the level of distribution for those terms that are most sensitive to such an increase. 
To illustrate this point, we present two different terms that exhibit different sensitivities as the level of distribution improves.

The first term is  
\[
G_1 =S(\mathscr{A},\mathscr{P}(N),N^{4/53}).
\]  
In our proof we have  
\begin{equation}
G_1\ge g_1\dfrac{C(N)N}{\log^2 N}(1+O(E)),
\end{equation}
where  
\begin{equation}
g_1=\dfrac{53}{2}e^{-\gamma}f\bigg(\dfrac{53}{8}\bigg)\ge 14.87710,
\end{equation}
and we have used
\begin{equation}
f\bigg(\dfrac{53}{8}\bigg)\ge f(6)\ge 0.999895, \quad e^{-\gamma}\ge 0.561459. 
\end{equation}
If we replace the Bombieri--Vinogradov theorem by the Elliott--Halberstam Conjecture, we obtain  
\begin{equation}
G_1\ge g_1'\dfrac{C(N)N}{\log^2 N}(1+O(E)),
\end{equation}
with  
\begin{equation}
g_1'=\dfrac{53}{2}e^{-\gamma}f\bigg(\dfrac{53}{4}-\varepsilon\bigg)\ge 14.87851,
\end{equation}
where  
$$
f\bigg(\dfrac{53}{4}-\varepsilon\bigg)\ge 0.99999. 
$$
Thus the gain from the Elliott--Halberstam Conjecture in this term is tiny.  

The second example is  
\begin{equation}
G_3 =\sum\limits_{\substack{N^{9/19}\le p\le N^{1/2} \\(p, N)=1}}S(\mathscr{A}_p ; \mathscr{P}(N), p),
\end{equation} 
In our proof we have  
\begin{equation}
G_3 \le g_3\dfrac{C(N)N}{\log^2 N}(1+O(E)),
\end{equation}
where  
\begin{equation}
g_{3}=8\log(10/9)\le 0.84289. 
\end{equation}
If we instead use the weighted Elliott--Halberstam Conjecture, we obtain  
\begin{equation}
G_3 \le g_3'\dfrac{C(N)N}{\log^2 N}(1+O(E)),
\end{equation}
with  
\begin{equation}
g_{3}'=4\log(10/9)\le 0.42145. 
\end{equation}
Here the benefit of the weighted Elliott--Halberstam Conjecture is substantial.

These two examples clearly show that the effect of improving the level of distribution is highly uneven across different sieve terms. 
A careful design of the weighted sieve is therefore essential to exploit the available analytic information optimally. 

\subsubsection{The calculus of variations cannot help us find the most suitable sieve function and parameters}   
The calculus of variations has achieved great success in mathematics and physics, providing a convenient framework for deriving the functional form that maximizes a given target parameter. 
Yet, in the context of the weighted sieve method presented here, it cannot be applied—for two principal reasons.

First, the operations leading to weighted inequalities are inherently discrete, such as the Buchstab identity and the switching principle. 
Moreover, the resulting expressions depend critically on the ranges of parameters. 
For instance, in the setting of Proposition \ref{Goldbachweight2}, if the condition \(\beta + \gamma + 2\delta < 1\) is changed to \(\beta + \gamma + 2\delta > 1\), the form of the weighted inequality changes dramatically.

Second, even after an appropriate range is chosen and the form of the weighted inequality is fixed, the calculus of variations remains impractical because the relevant auxiliary functions \(F\), \(f\), and \(w\) are defined by differential-difference equations. 
The calculus of variations is naturally adapted to differential equations, not to difference equations. 

\medskip 
\noindent 
{\bf Acknowledgements}. The authors record their thanks to many people, in particular Sergey Zheleznov, for 
checking the first version of the manuscript. The authors are supported by the National Key Research and Development 
Program of China (No. 2021YFA1000700), and the National Natural Science Foundation of China (Nos. 12031008 \& 123B2001).

\end{document}